\documentclass[11pt,letterpaper]{article}
\usepackage[margin=1in]{geometry}
\usepackage{natbib}

 \bibpunct[, ]{[}{]}{,}{n}{}{,}%

\usepackage[colorlinks=true,breaklinks=true,bookmarks=true,urlcolor=blue,
     citecolor=blue,linkcolor=blue,bookmarksopen=false,draft=false]{hyperref}

\usepackage{graphics,color} 
\usepackage{epsfig} 
\usepackage{amsmath} 
\usepackage{amssymb}  
\usepackage{wrapfig}

\usepackage{subfig}
\usepackage{multirow}
\usepackage{pifont}
\usepackage{tcolorbox}

\usepackage{algorithm}
\usepackage{algpseudocode}
\algnewcommand{\Initialize}[1]{%
	\State \textbf{Initialization:}
	\Statex {\raggedright #1}
}

\usepackage[colorinlistoftodos]{todonotes}

\usepackage{comment}

\usepackage{multirow}
\newcommand{\us}[1]{\begin{color}{black}#1\end{color}}
\newcommand{\uvs}[1]{\begin{color}{black}#1\end{color}}
\newcommand{\fy}[1]{\begin{color}{black}#1\end{color}}

\newcommand{\yqro}[1]{\begin{color}{black}#1\end{color}}
\newcommand{\fyy}[1]{\begin{color}{black}#1\end{color}}

\usepackage{accents}
\newcommand\munderbar[1]{%
  \underaccent{\bar}{#1}}

\let\proof\relax

\usepackage[para,online,flushleft]{threeparttable}
 \usepackage{amsthm}
 
 \newtheorem{assumption}{Assumption}
\newtheorem{theorem}{Theorem}
\newtheorem{lemma}{Lemma}
\newtheorem{proposition}{Proposition}
\newtheorem{definition}{Definition}
\newtheorem{corollary}{Corollary}
\newtheorem{remark}{Remark}

\title{\LARGE \bf
Zeroth-Order Federated Methods for Stochastic MPECs and Nondifferentiable Nonconvex Hierarchical Optimization
}
\author{Yuyang Qiu\thanks{Postdoctoral Scholar in the Department of Electrical and Computer Engineering,  University of California, Santa Barbara, CA 93106,  USA ({\tt\small yuyang\_qiu@ucsb.edu}).}
\and Uday V.  Shanbhag\thanks{Professor in the Department of Industrial and Manufacturing Engineering, Pennsylvania State University, University Park, PA 16802, USA  ({\tt\small{udaybag@psu.edu}}).}
\and Farzad Yousefian\thanks{Associate Professor in the Department of Industrial and Systems Engineering, Rutgers University, Piscataway, NJ 08854, USA
  ({\tt\small{farzad.yousefian@rutgers.edu}}).
 We acknowledge funding from the U.S. Department of Energy under grant \#DE-SC0023303, Air Force Office of Scientific Research (AFOSR) 
Grant FA9550-24-1-0259, and the U.S. Office of Naval Research under grants \#N00014-22-1-2757 and \#N00014-22-1-2589.}
}

\begin{document}
\sloppy
\maketitle
\thispagestyle{empty}
\pagestyle{plain}

\maketitle
\begin{abstract}{Motivated by the emergence of federated learning (FL), we design and analyze communication-efficient decentralized computational methods for addressing four broadly applicable but hitherto unstudied problem classes
in FL: (i) Nondifferentiable nonconvex optimization; (ii) Bilevel optimization;
(iii) Minimax problems; {and} (iv) Two-stage stochastic mathematical programs
with equilibrium constraints ({\bf 2s-SMPEC}). Notably, in an implicit sense, (ii),
(iii), and (iv) are instances of (i). However, these hierarchical problems are
often complicated by the absence of a closed-form expression for the implicit
objective function. Research on these problems has been limited and afflicted
by reliance on strong assumptions, including the need for differentiability of
the implicit function and {the} absence of constraints in the lower-level
problem, among others. We address these shortcomings by making the following
contributions. In (i), by leveraging convolution-based smoothing and Clarke’s
subdifferential calculus, we devise a randomized smoothing-enabled zeroth-order
FL method and derive communication and iteration complexity guarantees for
computing an approximate Clarke stationary point. To contend with (ii) and
(iii), we devise a unifying randomized implicit zeroth-order FL framework,
equipped with explicit communication and iteration complexities. Importantly,
our method utilizes delays during local steps to skip calls to the
inexact lower-level FL oracle,  resulting in significant reduction in
communication overhead. In (iv), we devise an inexact implicit variant of
the method in (i).  Remarkably,  this method achieves a total communication
complexity {matching} that of single-level nonsmooth nonconvex optimization in
FL. This appears to be the first time that {\bf 2s-SMPEC}s are provably addressed in
FL. 
We empirically validate the theoretical \uvs{claims} on instances of
federated nonsmooth and hierarchical problems including training of
ReLU neural networks, hyperparameter learning, fair classification, and
Stackelberg-Nash-Cournot equilibrium seeking.\footnote{A
preliminary version of this work has been accepted at The 37th Annual
Conference on Neural Information Processing Systems (NeurIPS
2023)~\cite{qiu2023zerothorder}.}}
\end{abstract}

\section{Introduction}
Federated learning (FL) has emerged as a promising enabling framework
for learning predictive models from a multitude of distributed,
privacy-sensitive \fyy{local datasets, which may exhibit varying data distributions (i.e., non-iid)}. This is accomplished
through the use of efficiently devised periodic communications between a
central server and a collection of clients. The FL algorithmic framework
allows for addressing several key obstacles in the development and
implementation of standard machine learning methods in a distributed and
parallel manner. For instance, the conventional parallel\us{ized} stochastic gradient
descent (SGD) method requires the exchange of information among the computing
nodes at every single time step, resulting in excessive communication
overhead \cite{lian2017can}. In contrast, FL methods including FedAvg~\cite{FedAvg17} and Local
SGD~\cite{Stich19a} overcome this onerous communication bottleneck by
provably attaining the linear speedup of parallel SGD by using 
significantly fewer communication
rounds~\cite{haddadpour2019local,FedAc,MimeMVR,das2022faster}. These guarantees have
been further complemented by recent
efforts~\cite{khaled2020tighter,FedBuff} where the presence of both \yqro{non-iid local datasets (i.e., the distribution of local data among clients is different, which may lead to client drift \cite{karimireddy2020scaffold})} and device heterogeneity
(i.e., variability of edge devices in computational power, memory, and
bandwidth) have been addressed. \fyy{It has been shown that the vanilla FedAvg method can be provably equipped with convergence guarantees under both iid and non-iid data settings, although convergence under non-iid data may require more communication rounds to achieve the same level of accuracy~\cite{khaled2020tighter}.}

Despite recent advances, much needs to be understood about
designing communication-efficient decentralized methods  for resolving {problems afflicted by nonconvexity, nonsmoothness, and the presence of hierarchy; to this end, we consider the} following broadly applicable problem classes (see Section~\ref{sec:num} for details of motivating examples for each class):

\smallskip

\noindent (a) {\em Nondifferentiable nonconvex locally constrained FL.} Consider the prototypical FL setting given as
\begin{align}\tag{\bf FL$_{nn}$}
    \min_x \quad \left\{ \, f(x) \triangleq  \frac{1}{m}\sum_{i=1}^m  \mathbb{E}_{\xi_i \in\mathcal{D}_i} [\, {\tilde{f}_i}(x,\xi_i) \, ] \, {\large \mid} \, 
    x \, \in \, X \, \triangleq \, \bigcap_{i=1}^mX_i \, \right\},\label{prob:main}
\end{align}
where $f$ is a nonsmooth nonconvex global loss function {defined over} a group of $m$ clients, indexed by $i \, \in \, [\, m \, ] \, \triangleq \, \{1,\ldots,m\}$,
$\mathcal{D}_i$ denotes the local dataset, ${\tilde{f}_i} :\mathbb{R}^n\times
\mathcal{D}_i \to \mathbb{R}$ denotes the local loss function, and $X_i \subseteq
\mathbb{R}^n$ is an easy-to-project local constraint set. Notably, 
local datasets may vary across clients, allowing for \yqro{non-iid data}.  We also consider client-specific local sets to induce
personalization. 

\smallskip

\noindent (b) {\em Constrained bilevel FL.} Overlaying a bilevel term {in \eqref{prob:main}} leads to 
\begin{align} \tag{\bf FL$_{bl}$}\label{prob:main_BL}
    &    \min_{x  \, \in \, X \, \triangleq \, \bigcap_{i=1}^mX_i}\left\{ \frac{1}{m}\sum_{i=1}^m  \mathbb{E}_{\xi_i \in\mathcal{D}_i}  [\, {\tilde{f}_i}(x,y(x),\xi_i) \, ] \mid  
    y(x) {=} {\displaystyle \operatornamewithlimits{\mbox{argmin}}_{y\in Y(x) \subseteq \mathbb{R}^{\tilde{n}}}}\,   \frac{1}{m}\sum_{i=1}^m  \mathbb{E}_{\zeta_i \in\tilde{\mathcal{D}}_i}[\tilde{h}_i(x,y,\zeta_i) ]\right\},\notag
\end{align}
where $y(\bullet):\mathbb{R}^n \to \mathbb{R}^{\tilde{n}}$ denotes a single-valued map returning the unique solution to the lower-level problem at $x$ and $Y(\bullet):\mathbb{R}^n \to \mathbb{R}^{\tilde{n}}$ denotes a {closed and convex-valued} set-valued map at any $x$.
 Let $f(\bullet) \triangleq  \frac{1}{m}\sum_{i=1}^m  \mathbb{E}_{\xi_i \in\mathcal{D}_i}  [\, {\tilde{f}_i}(\bullet,y(\bullet),\xi_i) \, ]$ denote the implicit objective function. 

\smallskip

\noindent (c) {\em Constrained minimax FL.} Consider the minimax setting, 
 defined as  
\begin{align}\tag{\bf FL$_{mm}$}
    &\min_{x \in X \triangleq \cap_{i=1}^m X_i} \max_{y\in {Y(x) \subseteq }\mathbb{R}^{\tilde{n}}} \quad  \frac{1}{m}\sum_{i=1}^m  \mathbb{E}_{\xi_i \in\mathcal{D}_i} [\, \tilde{f}_i(x,y,\xi_i)\,  ]  ,\label{prob:main_minimax}
\end{align}
where we assume that $y(x) = \mbox{arg}\max_{y\in Y(x)  }
\frac{1}{m}\sum_{i=1}^m  \mathbb{E}_{\zeta_i
\in\tilde{\mathcal{D}}_i}[\tilde{f}_i(x,y,\xi_i) ]$ is unique for all $x$. Let
$f(\bullet) \triangleq  \frac{1}{m}\sum_{i=1}^m  \mathbb{E}_{\xi_i
\in\mathcal{D}_i}  [\, {\tilde{f}_i}(\bullet,y(\bullet),\xi_i) \, ]$ denote the
implicit objective function. Indeed, problem~\eqref{prob:main_BL} subsumes this
minimax formulation when we choose $\tilde{h}_i:=-\tilde{f}_i$, {$\zeta_i \fyy{:=} \xi_i$}, and
$\tilde{\mathcal{D}}_i:=\mathcal{D}_i$.

\smallskip
\noindent (d) {\em Two-stage stochastic mathematical program with equilibrium constraints (2s-SMPEC) in FL.}
Consider a decentralized variant of a two-stage stochastic mathematical program~\cite{cui2022complexity} defined as
\begin{align}\tag{{\bf FL}$_{2s}$}\label{prob:fed2smpecs}
  \hspace{-0.4in}  \min_{x  \, \in \, X \, \triangleq \, \bigcap_{i=1}^mX_i}\left\{  \tfrac{1}{m}\sum_{i=1}^m  \mathbb{E}_{\xi_i \in\mathcal{D}_i}  [\, {\tilde{f}_i}(x,y(x,\xi_i),\xi_i) \, ]  \mid  
    y(x,\xi_i) =\text{SOL}(Y(x, \xi_i), G(x, \bullet , \xi_i)), \ \hbox{for a.e. $\xi_i$}\right\},\notag
\end{align}
where $G:\mathbb{R}^{n}\times\mathbb{R}^{\tilde{n}}\times \mathcal{D}_i \to
\mathbb{R}^{\tilde{n}}$ denotes a parametrized mapping, $Y(x, \xi_i)$ is closed
convex for any $x$ and $\xi_i$, and $\text{SOL}(Y(x, \xi_i), G(x, \bullet ,
\xi_i))$ denotes the solution set of the parametrized variational inequality
problem $\text{VI}(Y(x, \xi_i), G(x, \bullet , \xi_i))$ given an upper-level
decision $x$ and a random sample $\xi_i$. In this setting, we assume that $
y(x,\xi_i) $ is a unique solution to the lower-level VI problem and we let
$f(\bullet) \triangleq  \frac{1}{m}\sum_{i=1}^m  \mathbb{E}_{\xi_i
\in\mathcal{D}_i}  [\, {\tilde{f}_i}(\bullet,y(\bullet,\xi_i),\xi_i) \, ]$ denote
the implicit objective function.

\smallskip

{{\bf Key challenges.} Notably, in an implicit sense, problems (b), (c), and
(d) are all instances of problem (a). This is because these problems can be
cast as $\min_{x \in X} \ f(x)$ where $f(\bullet)$ denotes the implicit
objective function which we defined earlier for each of the problems (b), (c),
and (d) accordingly. \us{However}, the implicit function is afflicted by several challenges.
First, we note that generally $f(\bullet)$ is nondifferentiable and nonconvex.
Second, these hierarchical problems are often complicated by the absence of a
closed-form expression for the implicit objective.  As a consequence,  both
zeroth and first-order information of $f(\bullet)$ in problems (b), (c), and
(d) become unavailable,  making the design and analysis of computational
methods for these problems more challenging than that for (a).  Third,  to
construct an approximate zeroth-order evaluation of $f$, one may consider
employing inexact evaluations of the unique solution to the lower-level problem
using a suitable solver for the lower-level problem. This approach, however,
may result in \yqro{the} emergence of a bias \fyy{(caused by the inexact approximations) in} the zeroth-order gradients
which may then undesirably propagate throughout the implementation of the
method.  Fourth,  we assume that the lower-level problem in (b), (c), and (d)
can be characterized by a constraint set parametrized by the upper-level
decision $x$.  Such constraints create a challenge in addressing hierarchical
problems.  To elaborate,  consider, for example,  problem $\min_{x \in [-1,1]}\
\max_{y \in [-1,1],\ x+y\leq 0}\ x^2+y$ that admits the unique solution
$(x^*,y^*)=(0.5,-0.5)$. Now consider this problem with a \us{reversal in the} order of min
and max, i.e., $\max_{y \in [-1,1]}\ \min_{x \in [-1,1], \ x+y\leq 0}\  x^2+y,$
admitting the unique solution $(x^*,y^*)=(-1,1)$. \yqro{Consequently}, the
well-known primal-dual gradient methods, extensively employed
for addressing minimax problems with independent constraint sets, may fail to
converge to a saddle-point in minimax problems with coupling constrain\us{t}s. Last
but not least,  all the aforementioned challenges are exacerbated when {faced by 
a need for a} decentralized communication-efficient (i.e.,
federated) algorithmic framework. 

\smallskip

{\bf Related work.}}
\textit{\bf (i) Nondifferentiable nonconvex optimization.} Nonsmooth and nonconvex optimization has been studied extensively with convergence guarantees to Clarke-stationary points via gradient sampling~\cite{burke02approximating,burke05robust} and difference-of-convex approaches~\cite{cui22non}. Most complexity and rate guarantees necessitate smoothness of the nonconvex term \cite{ghadimi2016mini,xu2017globally,davis19stochastic,davis19proximally,lei2020asynchronous} or convexity of the nonsmooth term \cite{gasnikov2022power}, while only a few results truly consider nonsmooth nonconvex {objectives}~\cite{lin2022gradient,cui2022complexity,shanbhag2021zeroth}.
\textit{\bf (ii) Nondifferentiable nonconvex federated learning.} The research on FL
was initially motivated by decentralized neural networks where local
functions are nondifferentiable and nonconvex~\cite{FedAvg17}. Nevertheless,
    theoretical guarantees that emerged after FedAvg {required} either nonsmooth convex or smooth nonconvex local costs,
    under {either} iid ~\cite{Stich19a,K_AVG,coSGD21,Stich19b} or non-iid
    datasets~\cite{nonIIDFL20,khaled2020tighter}, while provable guarantees for FL methods under nonconvexity~\cite{K_AVG,ParSGDnonconvexFL19,nonconvexFL_haddadpour19} require $L$-smoothness of local functions. 
Unfortunately, these assumptions do not hold either for ReLU neural networks or risk-averse learning and necessitate the use of Clarke
calculus~\cite{clarke1975generalized}.
Moreover, existing work on zeroth-order FL methods in convex \cite{li2021communication} and nonconvex settings \cite{fang2022communication} rely on the smoothness properties of the objective function. \emph{However, there appear to be no provably convergent FL schemes with
    complexity guarantees for computing approximate Clarke stationary points
    of nondifferentiable nonconvex problems}. 
    \textit{\bf (iii) Federated bilevel optimization.}  Hyperparameter tuning~\cite{Hazan2018} {and its federated counterpart~\cite{FLreview21}} is a crucial, and yet, computationally complex integrant of machine learning (ML) pipeline. 
    Bilevel models where the lower-level is a parameterized training model while the upper-level requires selecting the best configuration for the unknown hyperparameters~\cite{GhW2018,JiYangLiang2021,tibshirani2005sparsity}. 
Solving such hierarchical problems is challenging because of nondifferentiable nonconvex terms and {an} absence of an analytical form for the implicit objective. These challenges exacerbate the development of provable guarantees for \yqro{data-}privacy-aware and communication-efficient schemes. \yqro{Existing methods in bilevel FL \fyy{often} consider the estimation of \uvs{the} hypergradient \cite{tarzanagh2022fednest}, \fyy{efficient} ways of computing \uvs{the} hypergradient \cite{yang2024simfbo,li2024communication}, and Lagrangian reformulations~\cite{yang2025first}, \fyy{among others}. Such schemes often rely on strong assumptions, such as the smoothness of the implicit function. Table \ref{tbl:flbl} summarizes and compares the assumptions made in recent works on bilevel FL.} {\emph{Few schemes, if any, exist for contending with bilevel problems with lower-level constraints.}} 
\textit{\bf (iv) Federated minimax optimization.} Minimax optimization {has assumed relevance} in adversarial learning
\cite{GAN14,NEURIPS2018_5a9d8bf5,sinha2018certifiable} and fairness in ML
\cite{fairgan18}, amongst other efforts. Recently, FL was extended to
distributed minimax
problems~\cite{AgnosticFL19,distrobustfl20,sharma2022federated}, but relatively
little exists in nonconvex-strongly concave settings.  In particular,  there
has been recent progress in addressing bilevel and minimax problems in FL,
including Local SGDA and
FedNest~\cite{tarzanagh2022fednest,sharma2022federated}. Our {proposed} work
{enjoys} two {key} distinctions with existing FL methods, described as
follows. (1) We do not require {either} the differentiability {or}
$L$-smoothness of the implicit objective, {assumptions that} may fail to
hold; e.g., in constrained hierarchical FL problems. (2) The existing FL
methods for bilevel and minimax problems assume that the lower-level problem is
unconstrained. In fact,  as emphasized earlier, even in centralized regimes
addressing hierarchical problems where the lower-level constraints depend on
$x$ have remained challenging.  Our proposed algorithmic framework allows for
accommodating these challenging problems in FL.  {\bf (v) SMPECs.}
Mathematical programs with equilibrium constraints (MPECs) provide a powerful
modeling framework for capturing {a range of hierarchical} decision-making problems in economics (such as the Stackelberg equilibrium problem), statistics, power markets, and
traffic networks (cf.~\cite{luo1996mathematical,outrata2013nonsmooth}). In the deterministic
regime,  {prior efforts have derived stationary conditions that cope with the ill-posedness of such problems}~\cite{scheel2000mathematical} {and extended} sequential quadratic programming schemes~\cite{jiang2000smooth} and
interior-point methods~\cite{leyffer2006interior} {to MPEC settings}.  In {uncertain}
regimes,  stochastic variants of MPECs have found
prime relevance in capturing leader-follower games. In two-stage
leader-follower games, the follower makes a second-stage decision $y$
contingent on the leader's decision $x$ and the realization of uncertainty
denoted by $\xi$.  In addressing 2s-SMPECs,  sample average approximation (SAA)
schemes have been studied~\cite{chen2015regularized,liu2011convergence}.  A key
shortcoming of the SAA approach is that the SAA problem becomes
increasingly difficult to solve as the number of scenarios of the second-stage
problem grows.  To resolve this issue and provide complexity
guarantees for the resolution of SMPECs, in a recent
work~\cite{cui2022complexity}, we devise {amongst the first known schemes with complexity guarantees for resolving SMPECs, by introducing} a class of {centralized} smoothing-enabled inexact zeroth-order methods. \yqro{We also note that a distributed zeroth-order gradient tracking method was developed~\cite{ebrahimi2024distributed} to \fyy{address} single-stage SMPECs.}
\begin{table}[htb]
\centering 
\scriptsize
\caption{Assumptions and methodologies of recent works on bilevel FL.}
\yqro{\begin{tabular}{|l|c|c|c|}
\hline
Algorithm & \begin{tabular}[c]{@{}c@{}}Lipschitz \\ $\nabla f_i$, $\nabla^2 h_i$\end{tabular} & \begin{tabular}[c]{@{}c@{}}Hypergradient \\ estimation \end{tabular} & Unconstrained $x$, $y$ \\ \hline
\begin{tabular}[c]{@{}c@{}}FedNest \cite{tarzanagh2022fednest} \end{tabular} & \ding{51} & \ding{51} & \ding{51} \\ \hline
\begin{tabular}[c]{@{}c@{}}SimFBO \cite{yang2024simfbo} \end{tabular} & \begin{tabular}[c]{@{}c@{}}\ding{51} \end{tabular} & \ding{55} & \ding{51} \\ \hline
\begin{tabular}[c]{@{}c@{}}FedBiOAcc \cite{li2024communication} \end{tabular} & \ding{51} & \ding{55} & \ding{51} \\ \hline
\begin{tabular}[c]{@{}c@{}}MemFBO \cite{yang2025first} \end{tabular} & \ding{51} & \ding{55} & \ding{51} \\ \hline
FedRZO$_{\text{bl}}$ (this work) & \ding{55} & \ding{55} & \ding{55} \\ \hline
\end{tabular}}
\label{tbl:flbl}
\end{table}

\begin{tcolorbox}
{\bf Gaps \& goal.} To the best of our knowledge, there are no known efficient FL algorithms that can contend with nondifferentiability and nonconvexity in an unstructured sense (as in \eqref{prob:main}). Generalizations {accommodating} a \uvs{constrained} bilevel \uvs{formulation}~\eqref{prob:main_BL}, minimax \uvs{formulations}~\eqref{prob:main_minimax}, or a two-stage program of the form~\eqref{prob:fed2smpecs} remain unaddressed in FL.  Motivated by these gaps, our work here is aimed at the development of a unified FL framework accommodating nondifferentiable nonconvex settings with extensions allowing for bilevel, minimax, and SMPEC {formulations}.
\end{tcolorbox}
We now  introduce the key ideas utilized in our work, postponing details until the subsequent sections.
 
{\bf A smoothed sampled zeroth-order framework.} We consider a smoothed framework for contending with constrained, nonsmooth, and nonconvex regimes. Specifically, given that $f$ is an expectation-valued function and $X$ is a closed and convex set, both of which are defined in \eqref{prob:main}, \fyy{an $\eta$-}smoothed unconstrained approximation is given as follows.
\begin{align}\label{fl_nn_eta}
    &\left\{ \begin{aligned}
\min & \ \tfrac{1}{m} \textstyle\sum_{i=1}^m \mathbb{E}_{\xi_i} [\, \tilde{f}_i(x,\xi_i) \,  ] \\
\mbox{subject to }
 & \ x \, \in \, \cap_{i=1}^m X_i.
    \end{aligned} \right\}  \,  \equiv  
\Big\{ \begin{aligned}
\min & \ \tfrac{1}{m} \textstyle\sum_{i=1}^m  \uvs{\big[}\, \mathbb{E}_{\xi_i} \, [\, \tilde{f}_i(x,\xi_i) \, ] + \underbrace{\mathbb{I}_{X_i}(x)}_{\tiny{\mbox{Indicator function}}}\,  \uvs{\big]}
\end{aligned} \Big\} \, \nonumber\\ 
 \tag{\bf FL$^{\eta}_{nn}$}
    & \overset{\mbox{\tiny Smoothing}}{\approx} 
\left\{ \begin{aligned}
\min & \  \tfrac{1}{m} \textstyle\sum_{i=1}^m \, \uvs{\bigg[} \underbrace{\mathbb{E}_{u_i \, \in \, \mathbb{B}}  [ \, \mathbb{E}_{\xi_i} [\, \tilde{f}_i(x+\eta u_i,\xi_i) \,  ] ]}_{\tiny  \mbox{Convolution smoothing}}\,  + \underbrace{\mathbb{I}^{\eta}_{X_i}(x)}_{\tiny\mbox{Moreau smoothing}} \uvs{\bigg]}\,  \end{aligned} \right\}.
\end{align}
If $f$ is {as} defined in \eqref{prob:main} and $d(x) \, \triangleq \, \tfrac{1}{m}\sum_{i=1}^m \mathbb{I}_{X_i}(x)$, then  ${\bf f}$ and its smoothing ${\bf f}^{\eta}$ are defined as 
\begin{align} 
    {\bf f}(x) \, & \triangleq \, f(x) + d(x) \mbox{ and } {\bf f}^{\eta}(x) \, \triangleq \, f^{\eta}(x) + {d^{\eta}(x)},\label{def:f_bold} \\
    \mbox{ where } f^{\eta} (x) & \triangleq  
    \tfrac{1}{m} \sum_{i=1}^m  [\, \mathbb{E}_{u_i \, \in \, \mathbb{B}}  [ \, \mathbb{E}_{\xi_i} [\, \tilde{f}_i(x+\eta u_i,\xi_i) \,  ]] ] \mbox{ and } d^{\eta}(x) \, \triangleq \, \tfrac{1}{m}\sum_{i=1}^m \mathbb{I}^{\eta}_{X_i}(x),\nonumber
\end{align}
where $\eta>0$ is a smoothing constant (see {\it Notation}  for details).

\noindent (i) {\em Clarke-stationarity}. Consider the original problem
\eqref{prob:main}. Under the assumption that the objective of \eqref{prob:main}
is Lipschitz continuous, then Clarke-stationarity of $x$ {with respect
to~\eqref{prob:main}} requires that $x$ {satisfies} $0 \,\in \, \partial f(x) +
{\cal N}_{\cap_{i=1}^m X_i}(x)$, where $\partial f(x)$ represents the Clarke
generalized gradient {\cite{clarke1975generalized}} of $f$ at $x$. However, a
negative result has been provided regarding the efficient computability of
$\epsilon$-stationarity in nonsmooth nonconvex
regimes~\cite{zhang20complexity}. Consequently, we focus on the \fyy{$\eta$-}smoothed
counterpart \eqref{fl_nn_eta}, a smooth nonconvex problem. In fact, under
suitable conditions \fyy{(see Proposition \ref{proposition:1})}, it can be shown that any stationary point of
\eqref{fl_nn_eta} is a $\yqro{\eta}$-Clarke stationary point of the original problem,
i.e.  $
    [ \, 0 \, \in \, \partial {\bf f}^{\eta}(x) \,  ] \, \implies \,  [\, 0 \, \in \, \partial_{\yqro{\eta}} {\bf f}(x) \,  ], 
$
where $\partial_{\yqro{\eta}} {\bf f}(x)$ represents the $\yqro{\eta}$-Clarke generalized gradient of ${\bf f}$ at $x$. 

\noindent (ii) {\em Meta-scheme for efficient resolution of \eqref{fl_nn_eta}.} We
develop {zeroth-order} stochastic gradient schemes for resolving \eqref{fl_nn_eta}.
This requires a zeroth-order gradient estimator for {$f^{\eta} (x)$},  denoted by
$\tfrac{1}{m}\sum_{i=1}^m {g}^\eta_i(x,{\xi_i},v_i)$ where $v_i \, \in \,
{\eta}\mathbb{S}$ for $i \, \in \, [ m ]$ and $\mathbb{S}$ denotes the surface of the
unit ball. Note that the Moreau smoothing of the
indicator function of $X_i$, denoted by $\mathbb{I}_{X_i}^{\eta}(x)$, admits a gradient, defined as $\nabla_x
\mathbb{I}_{X_i}^{\eta}(x) = \tfrac{1}{\eta}(x-\mathcal{P}_{X_i}(x))$, where
$\mathcal{P}_{X_i}(x) \triangleq \mbox{arg}\displaystyle{\min_{y \in X_i}} \|y-{x}\|^2$. The resulting {\em
meta-scheme} is defined next. 
 \begin{align}\tag{\bf Meta-ZO}
    x_{k+1} \, := \,   x_k - \gamma \left( \tfrac{1}{m}\textstyle\sum_{i=1}^m \left(g_{i}^{\eta}(x_k,{\xi_{i,k}},v_{i,k}) + \tfrac{1}{\eta}\left(x_k - \mathcal{P}_{X_i}(x_k)\right)\right)\right)  ,  \, k \geq 0.
\end{align}
{While the above scheme does not enforce feasibility at every iterate, we can show that the scheme generates an iterate that is $\mathcal{O}(\eta)$-feasible (see \fyy{Proposition \ref{proposition:1}~(ii)}).}  
 {Our first goal lies in leveraging this meta-scheme to design a randomized smoothing variant of the FedAvg method~\cite{FedAvg17} in which the clients {make} local steps of the form  $$x_{i,k+1}: = x_{i,k} -\gamma \left(  g_{i}^{\eta}(x_{i,k},{\xi_{i,k}},v_{i,k}) + \tfrac{1}{\eta}\left(x_{i,k} - \mathcal{P}_{X_i}(x_{i,k})\right)\right), \qquad \hbox{for } T_r \leq k \leq T_{r+1} -1,$$ 
where $x_{i,k} \in \mathbb{R}^{n}$ is a local copy of the decision $x$ maintained by client $i$ at local time index $k$, $r=0,1,\ldots$ denotes the communication round index, and $T_r$ is the global time index of communication round $r$ at which the server receives $x_{i,T_r}$ from clients,  does an aggregation step, and then sends the aggregated vector to clients.  This framework allows for preserving privacy \fyy{in terms of the local constraint sets} $X_i$. \yqro{We note that the term ``privacy" \fyy{carries} a different meaning \us{from that implied in} the field of differential privacy, where a typical approach \us{involves adding} random noise to the local information \cite{wei2020federated}.}
 }

\noindent {(iii) {\em Inexact implicit generalizations for \eqref{prob:main_BL}, \eqref{prob:main_minimax}, and \eqref{prob:fed2smpecs}.}} In addressing
the bilevel problem \eqref{prob:main_BL}, unlike in \eqref{prob:main},  we consider settings where the clients do not have access to the exact evaluation of the
implicit local objective $\tilde{f}_i(\bullet,y(\bullet),\xi_i)$,  
 a challenge addressed by developing {\em inexact}
implicit variants of the zeroth-order scheme, where clients compute an
$\varepsilon$-approximation of $y(x)$, denoted by $y_\varepsilon(x)$, in a
federated fashion.  Let $
{f}_{i}^{\eta}$ denote the \fyy{$\eta$-}smoothed  implicit local function.  We approximate $\nabla_{x} { f}_{i}^{\eta}{(x)}$ as follows.
\begin{align*}
    \nabla_{x} {f}_{i}^{\eta}(x)  = \underbrace{\mathbb{E}_{\xi_i,v} \left[g_{i}^{\eta}(x,\xi_i,v)\right]}_{\scriptsize \mbox{cannot be tractably evaluated}}\ \overset{(a.1)}{\approx}  \underbrace{g_{i}^{\eta}(x,\xi_{i,k},v_{T_r})}_{\scriptsize \mbox{intractable since $y(x)$ is unavailable}}\   \overset{(a.2)}{\approx} \  g_{i}^{\eta,\varepsilon}(x,\xi_{i,k},v_{T_r}),
\end{align*}
where $ g_{i}^{\eta,\varepsilon}(x,\xi,v)
\triangleq  \tfrac{n}{\yqro{2}\eta}({\tilde f}_i(x+v, y_{\varepsilon}(x+v),{\xi}) -
{\tilde f}_i(x\yqro{-v}, y_{\varepsilon}(x\yqro{-v}),{\xi}))\tfrac{v}{\|v\|}$ denotes an inexact
implicit zeroth-order gradient. Here, we contend with approximation of $\nabla_{x} { f}_{i}^{\eta}{(x)}$ through approximating the expectation by
sampling in $(a.1)$, while in $(a.2)$, we replace $y(x)$ by an inexact form $y_{\varepsilon}(x)$. This leads to the introduction of $g_i^{\eta,\varepsilon}(x,\xi,v)$.  Note that at each round of communication at the
upper level, $y_\varepsilon(x)$ can be computed using calls to a standard FL
method, e.g., FedAvg, in the lower level. {Notably,} such
calls to an FL oracle {should} be made only at the {global} step to preserve the
communication efficiency of the scheme. It follows that $
g_i^{\eta,\varepsilon}(x,\xi_{i,k},v_{T_r}) = \nabla_{x} {f}_{i}^{\eta}(x) +
\tilde{e}_{i,\varepsilon}$ where the approximation error
$\tilde{e}_{i,\varepsilon}$ is a {random variable,  possibly with a nonzero mean}. This bias \fyy{(caused by the inexact approximations)} can
then be mitigated by carefully updating the inexactness level $\varepsilon$ at
each communication round, as {done} in this work.  In addressing \eqref{prob:main_minimax} and \eqref{prob:fed2smpecs},  we employ a similar approach, with some distinctions {to be presented later.}

\yqro{\noindent (iv) {\it The need for zeroth-order methods in federated learning.} (a) In hierarchical \fyy{optimization} where the lower-level problem \fyy{is constrained}, the implicit function \fyy{might be} nondifferentiable (\fyy{see, for example,} problem \eqref{eg:bl} and Figure~\ref{fig:NDNCBiLv}). \fyy{Existing} methods \fyy{for} bilevel FL \fyy{often assume that the implicit function has Lipschitz continuous gradients}. This work appears to be the first to provably address bilevel FL problems under nondifferentiable and nonconvex settings. \fyy{(b)} In single-level FL, computing gradients in a high-dimensional neural network model (i.e., backpropagation) may be challenging for some clients, as high computing and memory requirements are needed. Zeroth-order methods are naturally backpropagation-free \cite{malladi2023fine}, relying solely on loss function evaluations, requiring less computation and memory at the client side \cite{zhang2024revisiting}, making them a promising alternative for resource-constrained settings.}

{\bf Contributions.}  In this work, we devise federated optimization methods with iteration and communication-complexity bounds for computing an approximate Clarke stationary point to the nonsmooth problem~\eqref{prob:main} and the hierarchical problem formulations given by \eqref{prob:main_BL}, \eqref{prob:main_minimax}, and \eqref{prob:fed2smpecs}.  In the hierarchical problems, we do not require differentiability of the implicit function. Instead, we consider a weakening in that the implicit function is assumed to be Lipschitz continuous and consider settings where the lower-level problem has a unique solution.  Our main contributions are summarized as follows.

\begin{table} [!htb]
 \centering
\newcommand{\tabincell}[2]{\begin{tabular}{@{}#1@{}}#2\end{tabular}}
  \caption{Comparison of our scheme with standard FL schemes for the nonconvex setting (see Proposition~\ref{thm:main_bound})}\label{TAB-lit}
 \centering
     \begin{tabular}{|l|c|c|c|c|c|}
        \hline
    Reference & Metric  &Rate &  Communication complexity & Assumption \\ \hline
    ~\cite{ParSGDnonconvexFL19}  & $\|\nabla_x f(x)\|^2$ & $ \tfrac{G^2}{\sqrt{mK}} $ & $ m^{3/4} K^{3/4}$ & \tabincell{c}{Bounded gradients,\\ $L$-smooth nonconvex functions}\\\hline  
    ~\cite{coSGD21}  & $\|\nabla_x f(x)\|^2$ & $ \tfrac{1}{\sqrt{mK}} $ & $ m^{3/2} K^{1/2} $ &  \tabincell{c}{$L$-smooth nonconvex functions} \\ \hline 
    ~\cite{haddadpour2019local} & $f(x)-f^*$ & $ \tfrac{1}{{mK}} $ & $ m^{1/3} K^{1/3} $ &  \tabincell{c}{$L$-smooth functions, PL-condition} \\ \hline 
    {\bf This work} & ${\|\nabla_x {\bf f}^{\eta}(x)\|^2}$ & $ \tfrac{1}{\sqrt{mK}} $ & $ m^{3/4} K^{3/4} $ &  \tabincell{c}{Lipschitz functions} \\ \hline 
      \end{tabular}
        \centering 
\end{table}

\begin{table}[ht]
\centering
 
\caption{Communication complexity for bilevel and minimax FL ($K$, $\tilde{K}$: maximum iteration in upper and lower level, resp.)}
\label{Table:FLbl}
\begin{center}
    \begin{tabular}{| c | c | c | c |}
    \hline
      Comm. complexity (upper level) & \multicolumn{2}{c|}{Lower level (standard FL schemes)}   &  Total comm. complexity \\ \hline
     \multirow{3}{*}{$m^{3/4} K^{3/4}$ ({\bf This work}, Theorem~\ref{thm:bilevel})}  & iid (Local SGD \cite{khaled2020tighter}) &$ m  $  & $ m^{\frac{7}{4}}K^{\frac{3}{4}} $ ({\bf This work})\\ \cline{2-4}  
      &iid (FedAC \cite{FedAc}) &  $ m^{\frac{1}{3}} $ &$ m^{\frac{{13}}{12}}K^{\frac{3}{4}} $( {\bf This work})\\ \cline{2-4}
      &\yqro{non-iid} (LFD \cite{haddadpour2019convergence}) &$\ m^{\frac{1}{3}}\tilde{K}^{\frac{1}{3}} $ & \fy{$  m^{\frac{11}{12}} K^{\frac{11}{12}} $} ({\bf This work})\\ \cline{2-4}  \hline
    \end{tabular}
\end{center}
\end{table}
 
\begin{table} [!htb]
 \centering
\newcommand{\tabincell}[2]{\begin{tabular}{@{}#1@{}}#2\end{tabular}}
  \caption{Rate and complexity results for federated two-stage SMPECs (see Theorem~\ref{thm:fed2smpec})}\label{TAB-mpec}
 \centering
     \begin{tabular}{|c|c|c|c|}
        \hline
  & Total comm.  complexity  & Rate (upper level) & Lower-level projections  \\ \hline

    {\bf This work} & \fy{$ m^{3/4} K^{3/4}$}& \fy{$\tfrac{1}{\sqrt{mK}} \ \ $  (for $K\geq m$)}&  \fy{ $mK/\text{ln}\left(\tfrac{1}{1- \kappa^{-2}_F}\right)$ }  \\ \hline 
      \end{tabular}
        \centering 
\end{table}

\noindent (i) {\it FL for nondifferentiable nonconvex problems.} To address
\eqref{prob:main} with \yqro{non-iid} datasets, we develop a Randomized
Zeroth-Order Locally-Projected Federated Averaging method ({FedRZO$_{\texttt
{nn}}$}). We derive iteration complexity of $\mathcal{O}\left(\tfrac{1}{
{m\epsilon^2}}\right)$ and communication complexity of
$\mathcal{O}\left(m^{3/4} K^{3/4}\right)$ for computing an
$\epsilon$-approximate Clarke stationary point to the smoothed problem where
$K$ denotes the termination time index. Such guarantees {appear} to be new in
the context of resolving nondifferentiable nonconvex FL problems, e.g. in
training of ReLU neural networks {(see Table~\ref{TAB-lit})}. {This is distinct
from existing zeroth-order methods, including
FedZO~\cite{fang2022communication}, that rely on differentiability and
$L$-smoothness of the local loss functions.}

\smallskip 
\noindent (ii) {\it Federated bilevel optimization.} In addressing
\eqref{prob:main_BL}, we develop FedRZO$_{\texttt {bl}}$, an inexact implicit
extension of FedRZO$_{\texttt {nn}}$, \yqro{where we employ an implicit programming approach and deal with the implicit function $f(\bullet)=\tfrac{1}{m}\sum_{i=1}^m \mathbb{E}_{\xi_i\in\mathcal{D}_i}  [\, {\tilde{f}_i}(\bullet,y(\bullet),\xi_i) \, ]$}. {A {\bf novelty} in the design of
FedRZO$_{\texttt {bl}}$ lies in the use of delays that facilitates {skipping}
calls to the inexact lower-level FL oracle during the local steps. As a
result, FedRZO$_{\texttt {bl}}$ is a} communication-efficient FL scheme with
single-timescale local steps, resulting in significant reduction in
communication overhead. Table~\ref{Table:FLbl} summarizes the  communication
complexity of this scheme. In all cases, we consider \yqro{non-iid} data in the
upper level. In the lower level, depending on what conventional FL scheme is
employed, we report the communication complexity accordingly. 

\smallskip 
\noindent (iii) {\it Federated minimax optimization.} FedRZO$_{\texttt {bl}}$ can be employed for addressing \eqref{prob:main_minimax} where $\tilde{h}_i: = -\tilde{f}_i $. As such, the complexity results in (ii) hold for solving (nondifferentiable nonconvex)-(strongly concave) FL minimax problems. Such results are new for this class of FL problems.

\smallskip 
\noindent (iv) {\it Federated two-stage SMPECs.}  Unlike in the bilevel
framework in \eqref{prob:main_BL},  the lower-level problem in
\eqref{prob:fed2smpecs} is not characterized by an expectation operator.
Instead, the equilibrium constraints hold for almost any realization of the
local random variables.  To contend with this stochastic hierarchical
structure, we propose a federated method called FedRZO$_{\texttt {2s}}$ in that
client $i$ employs a gradient method to approximate $y(\bullet,\xi_i)$
inexactly.  Analyzing this method, we derive explicit iteration and
communication complexity bounds for computing an $\epsilon$-approximate Clarke
stationary point to the smoothed implicit problem.  Importantly, we show that
FedRZO$_{\texttt {2s}}$ achieves a total communication complexity of $ m^{3/4}
K^{3/4}$, matching that of single-level nonsmooth nonconvex optimization
in FL. Table~\ref{TAB-mpec} presents a summary of the rate and complexity
results for FedRZO$_{\texttt {2s}}$. This is the first time that these
guarantees are provided for two-stage SMPECs in FL.  

\smallskip

{\bf Outline of the paper.} The paper is organized as follows.  After introducing notation, in Section~\ref{sec:nn} we consider nonsmooth nonconvex FL and present the convergence theory for the FedRZO$_{\texttt {nn}}$ method.  Then, in Section~\ref{sec:bl}, we {outline} FedRZO$_{\texttt {bl}}$ and derive guarantees for addressing \eqref{prob:main_BL}.  In Section~\ref{sec:mm},  we present the main convergence theory of FedRZO$_{\texttt {bl}}$ when applied to the minimax FL problem given by \eqref{prob:main_minimax}.  In Section~\ref{sec:2s}, we present FedRZO$_{\texttt {2s}}$ along with convergence guarantees for solving two-stage SMPECs.  To validate the theoretical results,  in Section~\ref{sec:num}, we provide a set of numerical experiments for FL problems including training ReLU neural networks,  decentralized hyperparameter tuning, and fair classification.  We provide some concluding remarks in Section~\ref{sec:conc}. For ease of exposition, we present all the proofs in Section~\ref{sec:proofs}. 

{\bf Notation and Preliminaries.} Throughout,  vectors are denoted by
lower-case letters and assumed to be column vectors.  Given two vectors $x,y
\in \mathbb{R}^n$, their Euclidean inner product is denoted by $x^{\fyy{\top}}y$.  We let
$\|x\|$ denote the Euclidean norm of vector $x$.  A mapping $F:X \to
\mathbb{R}^n$ is said to be merely monotone on a convex set $X \subseteq
\mathbb{R}^n$ if $(F(x)-F(y))^{\fyy{\top}}(x-y)\geq 0$ for all $x,y \in X$. The mapping
$F$ is said to be $\mu$-strongly monotone on $X$ if $\mu>0$ and
$(F(x)-F(y))^{\fyy{\top}}(x-y)\geq \mu\|x-y\|^2$ for any $x,y \in X$. Also, $F$ is said to
be Lipschitz continuous with parameter $L>0$ on the set $X$ if
$\|F(x)-F(y)\|\leq L\|x-y\|$ for all $x,y \in X$. The Euclidean projection of
vector $x$ onto a closed convex set $X$ is denoted by $ \mathcal{P}_{X}(x) $,
where $\mathcal{P}_{X}(x)\triangleq  \mbox{arg}\min_{ y \in  X}\| x- y\|$. We
let $\mbox{dist}(x,X)\triangleq \|x-\mathcal{P}_{X}(x)\|$ denote the distance
of $x$ from the set $X$.   Given a closed convex set $X$, we let
$\mathbb{I}_X:\mathbb{R}^n \to \{0,\fyy{\infty}\}$ denote the indicator function of
$X$.   Given a set $X \subseteq \mathbb{R}^n$ and a vector $x \in X$,  the
normal cone of $X$ at $x$ is defined as $\mathcal{N}_X(x) = \{y \in
\mathbb{R}^n \mid y^{\fyy{\top}}(z-x) \leq 0, \ \hbox{for any }z \in X\}$.  We let
$\mathbb{B}$ and $\mathbb{S}$ denote the $n$-dimensional unit ball and its
surface, respectively, i.e., $\mathbb{B} = \{u \in \mathbb{R}^n \mid \|u\| \leq
1\}$ and $\mathbb{S} = \{v \in \mathbb{R}^n \mid \|v\| = 1\}$. We use {\it
a.s.} and {\it a.e.} to  abbreviate ``almost surely" and ``almost everywhere",
respectively.  Throughout, we use $\mathbb{E}[\bullet]$ to denote the
expectation of a random variable.  The {\em Clarke generalized gradient} at $x$
is defined as $\partial h(x) \triangleq \left\{ \zeta \in \mathbb{R}^n \mid
h^{\circ}(x,v) \geq \langle \zeta, v\rangle, \  \forall v \in
\mathbb{R}^n\right\},$ where $h^{\circ}(x,v)$ is the {\em generalized
directional derivative}~\cite{clarke2008nonsmooth} of $h$ at $x$ in a direction
$v$, defined as $h^{\circ}(x,v)\triangleq {\textstyle \limsup_{y \to x, t
\downarrow 0}} \left({h(y+tv)-h(y)}\right)/{t}$. We denote the big O notation
by $\mathcal{O}\left(\bullet\right)$. When we ignore logarithmic and constant
numerical factors, we use $\tilde{\mathcal{O}}\left(\bullet\right)$ instead.

\section{Nonsmooth nonconvex federated optimization}\label{sec:nn}
In this section, we introduce a smoothing\us{-enabled} federated averaging framework for addressing \eqref{prob:main} and outline our main assumptions next.
\begin{assumption}\label{assum:Heterogeneous}\em 
Consider problem \eqref{prob:main}. The following hold. 

\noindent (i) \yqro{The function ${\tilde{f}_i}(\bullet,\xi_i)$ is {Lipschitz continuous with parameter $L_0(\xi_i) >0$} for all $i \in [m]$, where $L_0\triangleq\max_{i=1,\dots,m} \{\sqrt{\mathbb E_{\xi_i\in \mathcal D_i}[(L_0(\xi_i))^2]}\}<\fyy{\infty}$.}

\noindent \yqro{(ii)} For any $i \in [ m ]$, 
    client $i$ has access to a zeroth-order oracle ${\tilde{f}_i}(x,\xi_i)$ satisfying the following for every $x$ in an {almost-sure} sense: $\mathbb{E}_{\yqro{\xi_i\in \mathcal D_i}}[{{\tilde{f}_i}(x,\xi_{i})} \mid x] = f_i(x)$.  

\noindent (iii) The set $X_i$ is nonempty, closed, and convex for all $i \in [m]$. In addition, \us{for two positive scalars $B_1$ and $B_2$}, the following {\it bounded set-dissimilarity} condition holds for all $x\in \mathbb{R}^n$.
\begin{align}\label{assum:diss_set}
\tfrac{1}{m}\textstyle\sum_{i=1}^{m}\mbox{dist}^2(x,X_i) \leq B_1^2 + B_2^2\ {\left\|x- \tfrac{1}{m}\textstyle\sum_{i=1}^m\mathcal{P}_{X_i}(x)\right\|^2}.
\end{align} 

\end{assumption}
\yqro{Notably, Assumption~\ref{assum:Heterogeneous}~(ii) allows for non-iid data among the clients. \fyy{Such an assumption has been employed in prior works to study FL in settings with non-iid data~\cite{khaled2020tighter,karimireddy2020scaffold}.}} We note that the bounded set-dissimilarity {condition
\eqref{assum:diss_set}} is naturally similar to the so-called {\it bounded
gradient-dissimilarity} condition that has been employed in the literature,
e.g., in~\cite{karimireddy2020scaffold}. In particular, when the bounded
gradient-dissimilarity condition is stated for the Moreau smoothing
of the indicator function of $X_i$, denoted by $\mathbb{I}_{X_i}^{\eta}(x)$, we \us{arrive at} \eqref{assum:diss_set}. Notably, condition~\eqref{assum:diss_set}
holds for the generated iterate by the algorithm when, for example, the iterate
remains bounded. {As such, \eqref{assum:diss_set} can be viewed as a
weakening of the assumption that the generated iterate by the FL method is
bounded. }

\noindent {\bf Nonsmooth unconstrained reformulation.} Consider an unconstrained reformulation of \eqref{prob:main} given by ${\displaystyle \min_{x \in \mathbb{R}^n}}\ {\bf f}(x)$,  where ${\bf f}$ is defined as in~\eqref{def:f_bold}. 
Notably, the nonsmoothness of ${\bf f}$ arises from that of  $f$
and the local indicator functions {$\mathbb{I}_{X_i}$}.  The minimization of
${\bf f}$ is  challenging, {as noted by} recent findings \us{in nonsmooth nonconvex optimization}, where it \us{has been} shown~\cite{zhang20complexity} that for a suitable class
of nonsmooth functions, computing an $\mathbf{\epsilon}$-stationary point,
i.e., a point $\bar x$ for which $\mbox{dist}(0_n,\partial {\bf f}(\bar x))
\leq \epsilon$, is impossible in finite time. 

\noindent {\bf Approximate Clarke stationarity.} To circumvent this challenge, as a weakening of $\epsilon$-stationarity, a notion of $(\delta,\epsilon)$-stationarity is introduced~\cite{zhang20complexity} for a vector $\bar x$ when $\mbox{dist}(0_n,{\partial_{\delta}} {\bf f}(\bar x)) \leq \epsilon$, where the set 
$$\partial_{\delta} {\bf f}(x) \triangleq \mbox{conv}\left\{ \zeta: \zeta \in \partial {\bf f}(y), \|x-y\| \leq \delta\right\}$$ 
denotes the $\delta$-Clarke generalized gradient of ${\bf f}$ at $x$~\cite{goldstein77}; i.e. if $x$ is $(\delta,\epsilon)$-stationary, then there exists a convex combination of gradients in a $\delta$-neighborhood of $x$ that has a norm of at most $\epsilon$~\cite{shamir21nearapproximatelystationary}. 

This discussion naturally leads to the following key question: {\em
Can we devise provably convergent FL methods for computing approximate Clarke
stationary points of minimization of ${\bf f}$?} The aim of this section is to
provide an answer to this question by proposing a
{zeroth-order FL} method that employs smoothing. To contend with the nonsmoothness,
we employ the Moreau-smoothed variant of $\mathbb{I}_X(x)$, {where $ X  \triangleq  \bigcap_{i=1}^mX_i$,} and a randomized
smoothed variant of $f$.  These techniques are described next.

{\noindent {\bf Moreau smoothing for constraints.} Given a smoothing parameter $\eta>0$, the Moreau-smoothed variant of {$\mathbb{I}_{X_i}$} is given by {$\mathbb{I}^\eta_{X_i}$}, defined as $\mathbb{I}^{\eta}(x) \triangleq \frac{1}{2\eta}\mbox{dist}^2(x,X_i)$.  Note that $\mathbb{I}^\eta_{X_i}(x)$ is $\frac{1}{\eta}$-smooth. }

\noindent {\bf Randomized smoothing of loss function. } For {smoothing} the loss function $f$,  we employ a {\em randomized
smoothing} approach where the smoothing parameter is maintained as sufficiently small. This framework is rooted in the seminal work
by Steklov~\cite{steklov1}, leading to progress in both
convex~\cite{DeFarias08,YNS_Automatica12,Duchi12} and nonconvex~\cite{nesterov17}
regimes.  We consider a smoothing of $f$, given by {$f^{\eta}$}, defined as $
f^{\eta}(x) \triangleq  \mathbb{E}_{u \in \mathbb{B}}[f(x+\eta u)]$, where $u$ is a random vector in the unit ball $\mathbb{B}$. 
Further, $\mathbb{S}$ {denotes} the surface of the ball $\mathbb{B}$
and $\eta\mathbb{B}$ and $\eta\mathbb{S}$ denote a {ball} with radius $\eta$ and its
surface, respectively. \yqro{In this work, we employ spherical smoothing rather than Gaussian smoothing \fyy{employed in some prior work on zeroth-order methods~\cite{nesterov17}}. \fyy{A key reason is as follows.} Let us use $g^\eta$ to temporarily denote the zeroth-order gradient of an $L_0$-Lipschitz continuous function. \fyy{When} Gaussian smoothing \fyy{is employed}, it is shown that $\mathbb E[\|g^\eta\|^2]\leq \mathcal O(L_0^2n^2)$ \cite{nesterov17}. \fyy{However, when spherical smoothing is employed, this bound improves to} $\mathcal O(L_0^2n)$~\fyy{\cite{lin2022gradient}}. 
}

\begin{lemma}[Randomized spherical smoothing]\label{SphericalSmooth}\em      
Let $h:\mathbb{R}^n\to \mathbb{R}$ be a given continuous function and define $h^{\eta}(x)\triangleq \mathbb{E}_{u \in \mathbb{B}} \left[ h(x+ \eta u)  \right].$ Then, the following hold. 

\noindent (i) $h^\eta$ is continuously differentiable and for any $x\in \mathbb{R}^n$, 
$$
\nabla h^\eta(x) = \left(\tfrac{n}{\eta}\right)  \mathbb{E}_{v \in \mathbb{\eta S}}  \left[  h(x+ v)  \tfrac{v}{\|v\|}  \right] = \fyy{\left(\tfrac{n}{2\eta}\right)  \mathbb{E}_{v \in \mathbb{\eta S}}  \left[  (h(x+ v)-  h(x- v)) \tfrac{v}{\|v\|}  \right]}.$$  

Suppose $h$ is Lipschitz continuous with parameter $L_0>0$. Then, the following statements hold.

\noindent (ii) $| h^{\eta}(x) -  h^{\eta}(y)| \leq L_0 \|x-y \|$ for all $x,y\in \mathbb{R}^n$; 

\noindent (iii) $| h^{\eta}(x) -  h(x)| \leq L_0\eta$ for all $x\in \mathbb{R}^n$; 

\noindent (iv) \yqro{$\| \nabla h^{\eta}(x) - \nabla h^{\eta} (y)\| \leq \frac{L_0 \sqrt{n}}{\eta}  \|x-y \|$ for all \fyy{$x,y\in \mathbb{R}^n$.}}

%
%
%

\end{lemma}
The discussion leads to the consideration of the following smoothed federated problem. 
\begin{definition}[\underline{Unconstrained smoothed approximate problem}]
Given {$\eta > 0$}, consider an unconstrained smoothed problem given as
\begin{align}
\min_{x \in \mathbb{R}^n} \ {\bf f}^{\eta}(x) \left\{\triangleq \tfrac{1}{m}\textstyle\sum_{i=1}^m {\bf f}_i^{\eta}(x)\right\},  \hbox{where }  {\bf f}_i^{\eta}(x) \triangleq  \mathbb{E}_{\xi_i,u_i\in\mathbb{B}}[{\tilde{f}_i}\left(x+\eta u_i,\xi_i\right)] +\tfrac{1}{2\eta}\mbox{dist}^2(x,X_i).\label{prob:main_unconstrained}
\end{align}
\end{definition}

To solve the smoothed problem~\eqref{prob:main_unconstrained}, we propose FedRZO$_{\texttt {nn}}$ (Algorithm~\ref{alg:PZO_LSGD}). Here, client $i$ employs a zeroth-order stochastic gradient of the form $g^{\eta}_{i,k} \triangleq \frac{n}{\yqro{2}\eta^2} \left({\tilde{f}_i}(x_{i,k}+v_{i,k},\xi_{i,k}) - {\tilde{f}_i}(x_{i,k}\yqro{-v_{i,k}},\xi_{i,k})\right)v_{i,k}$, augmented by the gradient of the Moreau smoothed function. The random sample $v_{i,k} \in \eta\mathbb{S}$ is locally generated by each client $i$, allowing for randomized smoothing. \us{In view of Lemma~\ref{SphericalSmooth} (i), observe that this indeed} {facilitates} the development of a randomized zeroth-order gradient.
\begin{algorithm}[htb]                                                                                                                                 
	\caption{Randomized Zeroth-Order Locally-Projected Federated Averaging (FedRZO$_{\texttt {nn}}$)}\label{alg:PZO_LSGD}
{\begin{algorithmic}[1] 
	 \State \textbf{input:} Server chooses a random initial point $\hat x_0 \in X$, stepsize $\gamma$, smoothing parameter $\eta$, synchronization indices $T_0:=0$ and $T_r \geq 1$, where $r\geq 1$ is the communication round index
	 \For {$r = 0,1, \dots$}
	  \State Server broadcasts $\hat x_{r}$ to all clients: $x_{i,T_r}:= \hat x_r, \ \forall i \in [m]$
	 \For {$k = T_r,  \ldots,T_{r+1}-1$} {\bf in parallel by clients}
	 \State Client $i$ generates the random \yqro{sample} $\xi_{i,k} \in \mathcal{D}_i$ and $v_{i,k} \in \eta\mathbb{S}$
	 \State  $g^{\eta}_{i,k}: =\frac{n}{\yqro{2}\eta^2} \left({\tilde{f}_i}(x_{i,k}+v_{i,k},\xi_{i,k}) - {\tilde{f}_i}(x_{i,k}\yqro{-v_{i,k}},\xi_{i,k})\right)v_{i,k}$
	 \State  {Client $i$ does a local update as $x_{i,k+1} :=  x_{i,k} - \gamma \left( g^{\eta}_{i,k} + \frac{1}{\eta}\left(x_{i,k}-\mathcal{P}_{X_i}(x_{i,k}) \right)\right)  $}
	 
	 \EndFor
	\State Server receives $x_{i,T_{r+1}}$ from all clients and aggregates, i.e., $\hat x_{r+1} := \frac{1}{m}\sum_{i=1}^m x_{i,T_{r+1}}$

	 \EndFor
	\end{algorithmic}}
\end{algorithm}
Let $\bar{x}_k   \triangleq \frac{\sum_{i=1}^m x_{i,k}}{m}$ denote the averaged iterates of the clients. {Next, we define a communication frequency bound and present the main convergence result and complexity bounds for FedRZO$_{\texttt {nn}}$. }
\begin{definition}[\underline{Communication frequency bound}]\label{def:H} 
Consider Algorithm~\ref{alg:PZO_LSGD}. Let $H>0$ denote {a finite} upper bound on the communication frequency, i.e., $H \geq \max_{r =0,1,\ldots}{|T_{r+1} -T_r|}$.
\end{definition}

\begin{proposition}[Complexity bounds for FedRZO$_{\texttt {nn}}$]\label{thm:main_bound}\em \hspace{-0.15in}             
Consider Algorithm~\ref{alg:PZO_LSGD}. Let Assumption~\ref{assum:Heterogeneous} hold.  

\noindent (i) {\bf [Error bound]} Suppose $\gamma \leq \min\left\{\yqro{\tfrac{\eta}{4L_0\sqrt{n}}},\tfrac{1}{4H},\yqro{\tfrac{\eta}{12\sqrt{3} B_2(L_0\sqrt{n}+1)H}}\right\}$. Let $k^*$ denote \us{a random} integer drawn uniformly from $\{0,\ldots,K\}$, \yqro{where $K$ denotes the maximum iteration number} and ${\bf f}^{\eta,*}\triangleq \inf_{x}{\bf f}^{\eta}(x)$. Then \us{the following bound holds.} 
 \begin{align*}
  \hspace{-0.2in} \mathbb{E}\left[\|\nabla {\bf f}^\eta(\bar x_{k^*})\|^2\right]  &\leq  \tfrac{8( \mathbb{E}\left[ {\bf f}^{\eta}(\bar x_{0})\right] -  {\bf f}^{\eta,*} )}{\gamma(K+1)}   + \yqro{\tfrac{\gamma 64\sqrt{2\pi} L_0^3n^{1.5}}{\eta m}  +\tfrac{72H^2\gamma^2(L_0\sqrt{n}+1)^2}{\eta^2}\left(\tfrac{B_1^2}{\eta^2}+8\sqrt{2\pi}(4B_2^2+1)L_0^2n  \right) }  .
\end{align*}

\noindent (ii) {\bf [Iteration complexity]} Let $\yqro{\gamma:=\sqrt{\tfrac{\eta m}{Kn^{1.5}L_0^3}}}$ and $H:=   \left\lceil\sqrt[4]{\tfrac{K}{m^3}}\right\rceil$ where $\eta>0$. Let $\epsilon>0$ be an arbitrary scalar and $K_\epsilon$ denote the number of iterations such that $\mathbb{E}\left[\|\nabla f^\eta(\bar x_{k^*})\|^2\right] \leq \epsilon$. Then, the iteration complexity \us{$K_{\epsilon}$ for computing an $\epsilon$-solution is  as follows.} \yqro{$$K_{\epsilon}=\mathcal{O}\left(\left(\tfrac{L_0^{1.5}n^{0.75}}{\eta^{0.5}}+\tfrac{B_2^2L_0 \sqrt n}{\eta} \right)^2\tfrac{1}{m\epsilon^2}\right).$$}

\noindent (iii) {\bf [Communication complexity]} Suppose $K_\epsilon \geq m^3$. Then, the number of communication rounds to achieve the accuracy level in (ii) is $R:=\mathcal{O}\left((mK_\epsilon)^{3/4}\right)$.
\end{proposition}
\begin{remark}\label{rem:fedrzonn}
FedRZO$_{\texttt {nn}}$ achieves an iteration complexity of
$\mathcal{O}\left(\tfrac{1}{ {m\epsilon^2}}\right)$ and a communication
complexity of $\mathcal{O}\left(m^{3/4} K^{3/4}\right)$ in computing an
$\epsilon$-Clarke stationary point to the smoothed problem.  Both of these
bounds match  the standard existing bounds for smooth nonconvex federated
optimization~\cite{ParSGDnonconvexFL19}. \yqro{We note that the iteration complexity in (ii) in terms of the problem dimension $n$ is $\mathcal{\tilde O}\left(n^{1.5}\right)$, which matches the best known \us{dependence obtained for centralized schemes for resolving} nonsmooth nonconvex optimization~\cite{lin2022gradient}. Furthermore, in terms of smoothing parameter $\eta$, the iteration complexity is given as $\mathcal{\tilde O}\left(\eta^{-2}\right)$, combined with Lemma \ref{SphericalSmooth} (iii), we note a trade-off between \fyy{approximation accuracy and convergence speed}. Therefore, \fyy{selecting an appropriate} $\eta$ \fyy{may be crucial} for the performance of our method. We demonstrate the effect of different choices of $\eta$ in the numerical experiments.}
$\hfill$ $\Box$ 
\end{remark}
\yqro{Note that in Algorithm \ref{alg:PZO_LSGD}, we efficiently deal with the double expectations in \eqref{fl_nn_eta} through Monte Carlo sampling. Specifically, as shown in Algorithm \ref{alg:PZO_LSGD}, by leveraging only one sample of the smoothing vector $v_{i,k}$ per iteration, the computational cost remains reasonable in practice.}

We now formally relate the original nonsmooth problem and its smoothed
counterpart.  We note that Proposition~\ref{thm:main_bound} addresses the
resolution of the smoothed problem~\eqref{prob:main_unconstrained}.  Next, we
leverage the notion of approximate Clarke stationarity introduced earlier to
relate any stationary point to the smoothed problem with an approximate
stationary point to the original nonsmooth constrained
problem~\eqref{prob:main}. {Furthermore, any solution to $\eta$-smoothed
problem is \us{provably} $\mathcal{O}(\eta)$-feasible with respect to $X$.}
\begin{proposition}\em \label{proposition:1} Consider problem~\eqref{prob:main_unconstrained} and let Assumption~\ref{assum:Heterogeneous} hold. 

\noindent (i) Assume that $X_i=\mathbb{R}^n$ for all $i \in [m]$. Then, for any $\eta>0$, we have $\nabla f^{\eta}(x) \in \partial_{\yqro{\eta}} f(x)$. Consequently,  if $\nabla f^{\eta}(x)=0$, then $0 \in \partial_{\yqro{\eta}} f(x)$.

\noindent  (ii) Assume that the sets $X_i$ are identical for all $i \in [m]$. Let $\delta>0$ be an arbitrary scalar. If $\nabla \mathbf{f}^{\eta}(x) =0$ and $\eta \leq \frac{\delta}{\max\{\yqro{1,4\sqrt[4]{2\pi}L_0\sqrt{n}}\}}$, then $\|x - {\cal P}_X(x)\| \le \yqro{4\sqrt[4]{2\pi}\eta L_0\sqrt{n}}$ and $0_n \in \partial_{\delta} \left(  f + \mathbb{I}_X \right)(x).$

\end{proposition}

\section{Constrained bilevel federated optimization}\label{sec:bl}
Here, we design and analyze an FL method for {resolving} the bilevel optimization problem ~\eqref{prob:main_BL}.  Consider the following smoothed implicit problem.

\begin{definition}[\underline{Unconstrained smoothed implicit problem}] 
Given {\eqref{prob:main_BL} and} {$\eta > 0$}, consider an unconstrained smoothed {implicit} problem given as
\begin{align}
&\min_{x \in \mathbb{R}^n} \ {\bf f}^{\eta}(x) \left\{\triangleq \tfrac{1}{m}\textstyle\sum_{i=1}^m \left(\mathbb{E}_{\xi_i,u\in\mathbb{B}}[{\tilde{f}_i}\left(x+\eta u,y(x+\eta u),\xi_i\right)] +\tfrac{1}{2\eta}\mbox{dist}^2(x,X_i) \right)\right\} .\label{prob:main_unconstrained_bilevel}
\end{align}
\end{definition}

{The main assumptions utilized in this section are formally stated below.}
\begin{assumption}\label{assum:bilevel}\em 
Consider problem \eqref{prob:main_BL}. Let the following assumptions hold.
 
\noindent (i) For all $i \in [m]$, ${\tilde{f}_i}(\bullet,y(\bullet),\xi_i)$ is $L_0^{\tiny \mbox{imp}}(\xi_i)$-Lipschitz and ${\tilde{f}_i}(x,\bullet,\xi_i)$ is $L_{0,y}^f(\xi_i)$-Lipschitz for any $x$, where $L_0^{\tiny \mbox{imp}}\triangleq \displaystyle{\max_{i=1,\ldots,m}}\sqrt{\mathbb{E}_{\yqro{\xi_i\in \mathcal D_i}}[( L_0^{\tiny \mbox{imp}}(\xi_i))^2]}<\infty$ and $L_{0,y}^f \triangleq \displaystyle{\max_{i=1,\ldots,m}}\sqrt{\mathbb{E}_{\yqro{\xi_i\in \mathcal D_i}}[({L_{0,y}^f}(\xi_i))^2]}<\infty$.
\yqro{Furthermore, ${\tilde{f}_i}(x,y(x),\xi_i)$ \fyy{satisfies} the following for every $x$ in an {almost-sure} sense: $\mathbb{E}_{\xi_i\in \mathcal D_i}[{{\tilde{f}_i}(x,y(x),\xi_{i})} \mid x] = f_i(x,y(x))$.}

\noindent (ii) For all $i\in [m]$, for any $x$, $h_i(x,\bullet)\triangleq \mathbb{E}_{\zeta_i \in\tilde{\mathcal{D}}_i}[\tilde{h}_i(x,\bullet,\zeta_i) ]$ is $L_{1,y}^h$-smooth and $\mu_h$-strongly convex. Further, for any $y$, the map $\nabla_y h_i(\bullet,y)$ is Lipschitz continuous with parameter $L_{0,x}^{\nabla h}$. 

\noindent (iii) The sets $X_i$, for $i \in [m]$, satisfy Assumption~\ref{assum:Heterogeneous} (iii).

\noindent (iv) For any $x \in  \mathbb{R}^{n}$, the set $Y(x)\subseteq \mathbb{R}^{\tilde n}$ is nonempty,  closed, and convex.
\end{assumption}
\begin{remark}
As stated in Assumption~\ref{assum:bilevel} (i), we assume that the implicit function is Lipschitz continuous.  We note that this assumption is weaker than the differentiability of the implicit function which has been recently utilized frequently in the analysis of bilevel optimization problems in the literature.  Indeed,  unlike the differentiability property, the Lipschitz continuity  of the implicit function has been shown to hold under mild conditions in prior work, e.g., see~\cite[Lemma 2.2]{dafermos1988sensitivity}.  To elaborate, in the following result we consider a special case where the lower-level problem in \eqref{prob:main_BL} is unconstrained and proceed to not only show that the implicit function admits the Lipschitz continuity property, but also we obtain the Lipschitz parameter $L_0^{\tiny \mbox{imp}}(\xi_i) $ explicitly in terms of problem parameters. 

\yqro{We also note that the assumption that ${\tilde{f}_i}(x,y,\xi_i)$ is Lipschitz with respect to $y$ at any given \us{$x$} is realistic. The hyperparameter optimization problem \fyy{considered} in \fyy{S}ection \ref{sec:hyper} is shown \fyy{to meet this condition}. The details are presented in \fyy{S}ection \ref{sec:hyper}.}
$\hfill$ $\Box$  
\end{remark}
\begin{lemma}[Properties of the implicit function]\label{lem:implicit_props} \em
Consider the implicit function $ f(x)$ and mapping $y(x)$ given by \eqref{prob:main_BL}.  Assume that for all $i \in [m]$, {${\tilde{f}_i}(\bullet,y,\xi_i)$ is $L_{0,x}^f(\xi_i)$-Lipschitz for any $y$ and ${\tilde{f}_i}(x,\bullet,\xi_i)$ is $L_{0,y}^f(\xi_i)$-Lipschitz for any $x$}, where $L_{0,x}^f \triangleq \displaystyle{\max_{i=1,\ldots,m}}\sqrt{\mathbb{E}_{\yqro{\xi_i\in \mathcal D_i}}[({L_{0,x}^f}(\xi_i))^2]}<\infty$ and $L_{0,y}^f \triangleq \displaystyle{\max_{i=1,\ldots,m}}\sqrt{\mathbb{E}_{\yqro{\xi_i\in \mathcal D_i}}[({L_{0,y}^f}(\xi_i))^2]}<\infty$.  Further, let Assumption~\ref{assum:bilevel} (ii) and (iii) hold. Then, the following  hold.

\noindent (i) $y(\bullet)$ is $\left(\tfrac{L_{0,x}^{\nabla h}}{\mu_h}\right)$-Lipschitz continuous,  i.e., for any $x_1,x_2 \in \mathbb{R}^n$, $\|y(x_1) -y(x_2)\| \leq \left( \tfrac{L_{0,x}^{\nabla h}}{\mu_h}\right)\|x_1-x_2\|.$

\noindent (ii) The random implicit function is Lipschitz with constant $L_0^{\tiny \mbox{imp}}(\xi_i):=\tfrac{L_{0,y}^f(\xi_i) L_{0,x}^{\nabla h}}{\mu_h}+ L_{0,x}^f(\xi_i)$, i.e., 
\begin{align*} 
 | {\tilde{f}_i}(x_1,y(x_1),\xi_i) -{\tilde{f}_i}(x_2,y(x_2),\xi_i) |& \leq \ \left( \tfrac{L_{0,y}^f(\xi_i) L_{0,x}^{\nabla h}}{\mu_h}+ L_{0,x}^f(\xi_i)\right)\|x_1-x_2\| , \qquad \hbox{for all } x_1,x_2 \in \mathbb{R}^n .
\end{align*}

\noindent (iii) The implicit function is Lipschitz with parameter ${l_0^{\tiny \mbox{imp}}}:=\tfrac{L_{0,y}^f L_{0,x}^{\nabla h}}{\mu_h}+ L_{0,x}^f$, i.e.,  
\begin{align*} 
 | f (x_1)-  f (x_2) |& \leq \ \left( \tfrac{L_{0,y}^f L_{0,x}^{\nabla h}}{\mu_h}+ L_{0,x}^f\right)\|x_1-x_2\| , \qquad \hbox{for all } x_1,x_2 \in \mathbb{R}^n.
\end{align*}
\end{lemma}

\begin{algorithm}[htb]                                                                                                                                 
	\caption{Randomized Implicit \fy{Delay-Enabled} Zeroth-Order Federated Averaging (FedRZO$_{\texttt {bl}}$)}\label{alg:FedRiZO_upper}
{\begin{algorithmic}[1] 
	 \State \textbf{input:} Server chooses a random $\hat x_0 \in X$, stepsize $\gamma$, smoothing parameter $\eta$, synchronization indices $T_0:=0$ and $T_r \geq 1$, where $r\geq 1$ is the upper-level communication round index
	 \For {$r = 0,1, \dots$}
	 \State Server generates a random \fyy{sample} $v_{T_r} \in \eta\mathbb{S}$
	 \State Server calls FedAvg to receive $y_{\varepsilon_r}(\hat x_r+v_{T_r})$ and $y_{\varepsilon_r}(\hat x_r\yqro{-v_{T_r}})$, denoting the inexact evaluations of $y(\hat x_r+v_{T_r})$ and $y(\hat x_r\yqro{-v_{T_r}})$, respectively, {where $\mathbb{E}[\|y_{\varepsilon_r}(\bullet)-y(\bullet)\|^2] \leq \varepsilon_r$}
	  \State Server broadcasts $\hat x_{r}\yqro{-v_{T_r}}$,  $\hat x_{r}+v_{T_r}$, $y_{\varepsilon_r}(\hat x_r\yqro{-v_{T_r}})$, and $y_{\varepsilon_r}(\hat x_r+v_{T_r})$ to all {clients}; $x_{i,T_r}:= \hat x_r, \ \forall i$
	 \For {$k = T_r,  \ldots,T_{r+1}-1$} {\bf in parallel by {clients}}
	 \State Client $i$ generates \fyy{a} random \fyy{sample} $\xi_{i,k} \in \mathcal{D}_i$ 
	 \State  $g^{\eta,\varepsilon_r}_{i,k}: =\frac{n}{\yqro{2}\eta^2} \left( \tilde{f}_i(x_{i,k}+v_{T_r},y_{\varepsilon_r}(\hat x_r+v_{T_r}),\xi_{i,k})- \tilde{f}_i(x_{i,k}\yqro{-v_{T_r}},y_{\varepsilon_r}(\hat x_r\yqro{-v_{T_r}}),\xi_{i,k})\right)v_{T_r}$
	 \State  Client $i$ does a local update as $x_{i,k+1} :=  x_{i,k} - \gamma \left(  g^{\eta,\varepsilon_r}_{i,k} + \frac{1}{\eta}\left(x_{i,k}-\mathcal{P}_{X_i}(x_{i,k}) \right) \right)    $

	 \EndFor	
	\State Server receives $x_{i,T_{r+1}}$ from all clients and aggregates, i.e., $\hat x_{r+1} := \frac{1}{m}\sum_{i=1}^m x_{i,T_{r+1}}$

	 \EndFor
	\end{algorithmic}}
\end{algorithm}

{To address the smoothed problem~\eqref{prob:main_unconstrained_bilevel}, we
propose FedRZO$_{\texttt {bl}}$}, {outlined in
Algorithm~\ref{alg:FedRiZO_upper}}. We make the following remarks: (i) At each
global step, the server makes two calls to a lower-level FL method to inexactly
compute $y(\hat x_r+v_{T_r})$ and $y(\hat x_r\yqro{-v_{T_r}})$. These lower-level FL calls are
performed by the same clients, on the lower-level FL problem. (ii) {As
discussed later,} the inexactness error can be carefully controlled by
terminating the lower-level FL oracle after a prescribed number of iterations
characterized by $r$, where $r$ denotes the upper-level communication round
index. (iii) During the local steps in solving the upper-level FL problem,
FedRZO$_{\texttt {bl}}$ skips making inexact calls to the lower-level FL
oracle.  This helps with reducing the overall communication overhead
significantly. To accommodate this, unlike in FedRZO$_{\texttt {nn}}$, here we
employ a global randomized smoothing denoted by $v_{T_r}$ during the
communication round $r$ in the upper level.  
\begin{remark} 
    A technical challenge in designing Algorithm \ref{alg:FedRiZO_upper} is that inexact evaluations of $y(x)$ should be avoided during the local steps. This is because we consider bilevel problems of the form \eqref{prob:main_BL} where both levels are distributed. Because of this, the inexact evaluation of $y(x)$ by each client in the local step in the upper level would require significant communication,  {an undesirable outcome in an FL setting.} 
    We carefully address this challenge by introducing {\bf delayed inexact computation} of $y(x)$. As observed in step 8 of Algorithm \ref{alg:FedRiZO_upper}, we note how $y_{\varepsilon}$ \yqro{are} evaluated at $\hat x_r +v_{T_r}$ \yqro{and $\hat x_r -v_{T_r}$} which are different than the \yqro{vectors} used by the client, i.e., $x_{i,k} +v_{T_r}$ \yqro{and $x_{i,k} -v_{T_r}$}. At each communication round in the upper level, we only compute $y(x)$ inexactly twice in the global step \yqro{by calling FedAvg (Algorithm \ref{alg:FedZO_lower}) with $\tilde H$ local steps, where $\tilde H$ is defined in Algorithm \ref{alg:FedZO_lower}.} \yqro{Then,} this delayed information \yqro{is used} in the \yqro{upper-level} local steps. \yqro{The length of each delay is the same as the upper-level local iteration number $H$, as \fyy{the values of} $y_{\varepsilon_r}(\bullet)$ are obtained only at global steps.} This delayed inexact computation of $y$ {introduces new}  challenges  in the convergence analysis, {making} the the design and analysis of Algorithm~\ref{alg:FedRiZO_upper} a {\it non-trivial extension} of Algorithm~\ref{alg:PZO_LSGD}. $\hfill$ $\Box$
\end{remark}
{In the following theorem, we present the main convergence result and complexity bounds for FedRZO$_{\texttt {bl}}$. }
%
\begin{theorem}[Complexity bounds for FedRZO$_{\texttt {bl}}$]\label{thm:bilevel}\em
    Consider Algorithm~\ref{alg:FedRiZO_upper}. Let Assumption~\ref{assum:bilevel} hold {and let} $k^*$ be chosen uniformly at random from $0,\ldots,K:=T_R-1$.  Let $\gamma \leq \min\left\{\tfrac{\max\left\{{2, \sqrt{0.1\Theta_3},4B_2\sqrt{3\Theta_2},4B_2\sqrt{3\Theta_3}} \right\} ^{-1}}{4H},\yqro{\tfrac{\eta}{24({L_0^{\text{imp}}}\sqrt{n}+1)}}\right\}$. Let $\varepsilon_r$ denote the inexactness in obtaining the lower-level solution, i.e., $\mathbb{E}\left[\|{y}_{\varepsilon_r}(x) -y(x) \|^2 {\mid x }\right] \leq \varepsilon_r$ for $x \in \cup_{r=0}^{R}\{\hat x_r,\hat x_r+v_{T_r} \}$. 

\noindent (i) {\bf [Error bound]} \fyy{The following bound holds.}
\begin{align*} 
 & \mathbb{E}\left[\|\nabla {\bf f}^{\eta}(\bar x_{k^*})\|^2\right]    \leq   8(\gamma K)^{-1}(\mathbb{E}\left[{\bf f}^{\eta}(x_{0}) \right] -  {\bf f}^{\eta,*}  )   +\tfrac{8\gamma \Theta_1}{m}  + 8H^2\gamma^2\max\{\Theta_2,\Theta_3\}\Theta_5 \\
 & \qquad\qquad\qquad\qquad    + 8\left(H^2\gamma^2\max\{\Theta_2,\Theta_3\}\Theta_4+\Theta_3\right)H\tfrac{\textstyle\sum_{r=0}^{R-1}\varepsilon_r}{K}     , \\
&\mbox{where }\Theta_1 :=\yqro{\tfrac{48\sqrt{2\pi}(L_0^{\text{imp}}\sqrt{n}+1)(L_0^{\text{imp}})^2n}{2m\eta}}, \Theta_2:=\yqro{\tfrac{5(L_0^{\text{imp}}\sqrt{n}+1)^2}{8\eta^2}}, \text{ and } {\Theta_3:= \left( \tfrac{L_{0,x}^{\nabla h}}{\mu_h}\right)^2\tfrac{\yqro{5}  n^2}{\eta^2}(L_{0,y}^f)^2}. \\
& \Theta_4:=\tfrac{{\yqro{24}}n^2}{\eta^2}(L_{0,y}^f)^2, \text{ and } \Theta_5:= \tfrac{48 B_1^2}{\eta^2}+ \yqro{384\sqrt{2\pi}(4B^2_2+1)(L_0^{\text{imp}})^2n}.   
\end{align*}

\noindent (ii) {\bf  [Iteration complexity]} Let $\gamma:=\yqro{\sqrt{\tfrac{\eta m}{Kn^{1.5}(L_0^{\text{imp}})^3}}}$
    and $H:=   \left\lceil\sqrt[4]{\tfrac{K}{m^3}}\right\rceil$ where $\eta>0$.
    Let $\epsilon>0$ be an arbitrary scalar and $K_\epsilon$ denote the number
    of iterations such that $\mathbb{E}\left[\|\nabla {\bf f}^\eta(\bar
    x_{k^*})\|^2\right] \leq \epsilon$. Also, suppose we employ an FL method in
    the lower level  {with} a sublinear convergence rate {characterized by}  a
    linear speedup in terms of the number of clients,  i.e.,
    $\varepsilon_r:=\tilde{\mathcal{O}}(\tfrac{1}{m \tilde{T}_{\tilde{R}_r}})$
    where $\tilde{R}_r$ denotes the number of communication rounds performed in
    the lower-level FL method when called in round $r$ of FedRZO$_{\texttt
    {bl}}$ and $\tilde{T}_{\tilde{R}_r}$ denotes the number of iterations
    performed in the lower-level FL scheme to do $\tilde{R}_r$ rounds of
    upper-level communication. Further, suppose
    $\tilde{T}_{\tilde{R}_r}:=\tilde{\mathcal{O}}\left(m^{-1}(r+1)^{\frac{2}{3}}\yqro{n}\right)$.
    Then, the iteration complexity of FedRZO$_{\texttt {bl}}$ (upper level) is \yqro{$$K_{\epsilon}=\mathcal{O}\left(\left(\tfrac{n^{0.75}(L_0^{\text{imp}})^{1.5}}{\eta^{0.5}}+\tfrac{n^{1.5}(L_0^{\text{imp}})^{-3}}{\eta^3}+\tfrac{n}{\eta^2}+\tfrac{\sqrt n(L_0^{\text{imp}})^{-3}}{\eta^3}+\tfrac{n^{1.5}(L_0^{\text{imp}})^{-1}}{\eta} \right)^2\tfrac{1}{m\epsilon^2}\right).$$}

\noindent (iii) {\bf  [Communication complexity]}  Suppose $K_\epsilon \geq m^3$. Then, the number of communication rounds in FedRZO$_{\texttt {bl}}$ (upper level only) to achieve the accuracy level in (ii) is $R:=\mathcal{O}\left((mK_\epsilon)^{3/4}\right)$. 

\end{theorem}
\begin{remark}\label{rem:thm2} 
\noindent (i) Matching the complexity bounds of FedRZO$_{\texttt {nn}}$:
    Importantly, Theorem~\ref{thm:bilevel} is equipped with {an} explicit
    iteration complexity
    $K_\epsilon:=\tilde{\mathcal{O}}\left(\frac{1}{m\epsilon^2} \right)$ and
    communication complexity $R:=\mathcal{O}\left((mK_\epsilon)^{3/4}\right)$,
    matching those of single-level nonsmooth nonconvex problems, obtained in
    Proposition~\ref{thm:main_bound}. This is promising, implying that as long
    as the lower-level FL oracle has a rate of
    $\varepsilon_r:=\tilde{\mathcal{O}}(\tfrac{1}{m \tilde{T}_{\tilde{R}_r}})$,
    the inexactness in approximating the lower-level solutions does not
    {adversely} affect the complexity bounds of the method in the
    upper level.  This is {achieved}  by carefully terminating the
    lower-level oracle after
    $\tilde{T}_{\tilde{R}_r}:=\tilde{\mathcal{O}}\left(m^{-1}(r+1)^{\frac{2}{3}}\yqro{n^a}\right)$
    iterations, \yqro{where $a>0$ is a constant. \fyy{Invoking} this choice of $\tilde{T}_{\tilde{R}_r}$ \fyy{in} (i) leads to the iteration complexity of $\mathcal{\tilde O}(n^{\max\{5-2a,3\}})$. We note that $a=1$ leads to the optimal iteration complexity and overall communication complexity in terms of $n$ in all the three cases discussed below.}


\noindent  (ii) \yqro{Non-iid data}: As noted in Assumption~\ref{assum:bilevel}, we do not impose the condition $\mathbb{E}_{\yqro{\xi_i\in\mathcal D_i}}[{{\tilde{f}_i}(x,y(x),\xi_{i})} \mid x] = f(x)$, allowing \yqro{clients to have non-iid data \fyy{in the upper level}.}

    \noindent  (iii)  Overall communication complexity of FedRZO$_{\texttt
    {bl}}$: The overall communication complexity of FedRZO$_{\texttt {bl}}$
    depends on the standard FL method employed for the inexact approximation of
    the lower-level solutions. To elaborate, we provide specific complexity
    bounds in Table~\ref{Table:FLbl} for three instances where we employ Local
    SGD \cite{khaled2020tighter}, FedAC~\cite{FedAc}, and LFD
    \cite{haddadpour2019convergence} for inexactly solving the lower-level FL
    problem.  Notably,  all these schemes meet the rate condition
    $\varepsilon_r:=\tilde{\mathcal{O}}(\tfrac{1}{m \tilde{T}_{\tilde{R}_r}})$,
    stated in Theorem~\ref{thm:bilevel}. Amongst these, the last one
    allows for the presence of \yqro{non-iid datasets \fyy{in the lower level}}.  The detailed analysis {of}
    these cases is provided {next.}

    \noindent {\bf Case 1:} FedRZO$_{\texttt {bl}}$ employs Local SGD~\cite{khaled2020tighter} with iid datasets in round $r$ to obtain an $\varepsilon_r$-inexact solution to the lower-level problem~\eqref{prob:main_BL}.  {An} outline of Local SGD is presented in Algorithm~\ref{alg:FedZO_lower}. In view of Corollary 1 in~\cite{khaled2020tighter} under suitable settings and terminating Local SGD after $ \tilde{T}_{\tilde{R}_r}$ iterations, $\mathbb{E}\left[\|{y}_{\varepsilon_r}(x) -y(x) \|^2\right] \leq \varepsilon_r := \mathcal{\tilde O}\left(\tfrac{1}{\mu_h^2 m \tilde{T}_{\tilde{R}_r}}\right)$
    in $\tilde{R}_r=m$ {$^{\rm th}$} round of communications (in the lower level). Invoking Theorem \ref{thm:bilevel} (iii), the overall communication complexity is $m\times \mathcal{O}\left((mK_\epsilon)^{3/4}\right)$ that is $ \mathcal{O}\left(m^{7/4}K_\epsilon^{3/4}\right)$.

\noindent {\bf Case 2:} FedRZO$_{\texttt {bl}}$ employs FedAC~\cite{FedAc} with iid datasets in round $r$ to obtain $\varepsilon_r$-inexact solution to the lower-level problem~\eqref{prob:main_BL}. Similarly, invoking Theorem 3.1 in~\cite{FedAc}, the communication complexity of the lower-level calls is of the order $m^{1/3} $. Invoking Theorem \ref{thm:bilevel} (iii) again, the overall communication complexity is $m^{1/3}\times \mathcal{O}\left((mK_\epsilon)^{3/4}\right)$ that is $ \mathcal{O}\left(m^{{13}/12}K_\epsilon^{3/4}\right)$. \yqro{Note that in both case 1 and case 2, the lower-level communication complexity does not depend on $n$, and the overall communication complexity depend\us{s} only on the iteration $K_{\epsilon}$. Therefore, by minimizing the dependence \fyy{on} $n$ in $K_{\epsilon}$, the overall communication complexity is minimized in terms of \fyy{the choice of $a$}  \fyy{for} $a=1$.}

\noindent {\bf Case 3:} FedRZO$_{\texttt {bl}}$ employs
LFD~\cite{haddadpour2019convergence} with \yqro{non-iid} datasets in round
$r$ to obtain $\varepsilon_r$-inexact solution to the lower-level
problem~\eqref{prob:main_BL}. In view of Theorem
4.2~\cite{haddadpour2019convergence}, the number of local steps in the
lower-level calls is of the order $
m^{-\frac{1}{3}}\tilde{T}_{\tilde{R}_r}^{\frac{2}{3}}$, which implies the
lower-level communication complexity of the order $
m^{\frac{1}{3}}\tilde{T}_{\tilde{R}_r}^{\frac{1}{3}}$.  Invoking Theorem
\ref{thm:bilevel} (iii) and
$\tilde{T}_{\tilde{R}_r}:=\tilde{\mathcal{O}}\left(m^{-1}(r+1)^{\frac{2}{3}}\yqro{n^a}\right)$,
the overall communication complexity is
$\sum_{r=0}^{(mK_\epsilon)^{3/4}}\left(m^{\frac{1}{3}}\tilde{T}_{\tilde{R}_r}^{\frac{1}{3}}\right)
= \sum_{r=0}^{(mK_\epsilon)^{3/4}}\tilde{\mathcal{O}}\left(
r^{2/9}\yqro{n^{a/3}}\right)=\tilde{\mathcal{O}}\left(\left(mK_\epsilon\right)^{11/12}\yqro{n^{a/3}}
\right)$, \yqro{this can be written as
$\tilde{\mathcal{O}}(n^{\max\{5-2a,3\}11/12+a/3})$. We can observe that $a=1$
leads to the best overall communication complexity. Note that the \us{multiplication by} $n$ in the lower-level iteration number $\tilde{T}_{\tilde{R}_r}$ will not \us{worsen} \fyy{the} dependence on $n$ in the overall communication complexity. This is because without having dependence of $n^a$ in $\tilde{T}_{\tilde{R}_r}$,
the overall iteration complexity would be $K_{\epsilon}=\tilde{\mathcal{O}}(n^5)$, leading to the overall communication complexity $\tilde{\mathcal{O}}\left(\left(mK_\epsilon\right)^{11/12}\right)=\tilde{\mathcal{O}}\left(n^{55/12}\right)$,
which is worse than the \fyy{bound of} $\tilde{\mathcal{O}}(n^{37/12})$ \fyy{obtained when the dependence on $n^a$ is captured in $\tilde{T}_{\tilde{R}_r}$}.} $\hfill$
$\Box$
\end{remark}

\begin{algorithm}[htb]                                                                                                                                 
	\caption{Local SGD $(x,r,y_{0,r},m,\tilde{\gamma},\tilde{H},\tilde{T}_{\tilde{R}})$ for lower-level}\label{alg:FedZO_lower}
\begin{algorithmic}[1] 
	 \State \textbf{input:} $x$, $r$, server chooses a random initial point $\hat y_{0}:= y_{0,r} \in \mathbb{R}^{\tilde{n}}$, $a_r:=\max\{m,4\kappa_h,r\}+1$ where $\kappa_h:=\tfrac{L_{1,y}^h}{\mu_h}$, $\tilde{\gamma}:=\tfrac{1}{\mu_h a_r}$, $\tilde{T}_{\tilde{R}_r}:=2a_r\ln(a_r)$, and $\tilde{H}:=\lceil\tfrac{\tilde{T}_{\tilde{R}_r}}{m}\rceil$ 
	 \For {$\tilde{r} = 0,1,\ldots, \tilde{R}_r-1$}
	  \State Server broadcasts $\hat y_{\tilde{r}}$ to all agents: $y_{i,\tilde{T}_{\tilde{r}}}:= \hat y_{\tilde{r}}, \ \forall i$
	 \For {$t = \tilde{T}_{\tilde{r}},  \ldots,\tilde{T}_{\tilde{r}+1}-1$} {\bf in parallel by agents}
	  
	 \State  {Agent $i$ does a local update as $y_{i,t+1} :=  y_{i,t} -  \tilde{\gamma} \nabla_y \yqro{\tilde h_i}(x,y_{i,t},\yqro{\zeta_{i,t}}) $}
	 \State Agent $i$ sends $\yqro{y_{i,\tilde T_{\tilde r+1}}}$ to the server  
	 \EndFor
	\State Server aggregates, i.e., $\hat y_{\tilde{r}+1} := \frac{1}{m}\sum_{i=1}^m y_{i,\yqro{\tilde T_{\tilde{r}+1}}}$

	 \EndFor
	\end{algorithmic}
\end{algorithm}


\subsection{Constrained minimax federated optimization}\label{sec:mm}
Next, we consider the decentralized federated minimax problem of the form
\eqref{prob:main_minimax} introduced earlier. This problem is indeed a zero-sum
game and under the nonconvex-strongly concave setting, this is 
an instance of the non-zero sum game \eqref{prob:main_BL} where $\tilde{h}_i: =
-\tilde{f}_i $.  We elaborate on this equivalence by {considering}
\eqref{prob:main_minimax}. Recall that $y(x) = \mbox{arg}\max_{y\in Y(x)
} \frac{1}{m}\sum_{i=1}^m  \mathbb{E}_{\zeta_i
\in\tilde{\mathcal{D}}_i}[\tilde{f}_i(x,y,\xi_i) ]$ is assumed to be unique.
Via this definition,  for any $x$, $$\left[\max_{y\in
Y(x) \subseteq \mathbb{R}^{\tilde{n}}} \ \frac{1}{m}\sum_{i=1}^m
\mathbb{E}_{\xi_i \in\mathcal{D}_i} [\, \tilde{f}_i(x,y,\xi_i)\,  ] \right]
\quad = \quad \frac{1}{m}\sum_{i=1}^m  \mathbb{E}_{\xi_i \in\mathcal{D}_i} [\,
\tilde{f}_i(x,y(x),\xi_i)\,  ]. $$
This implies that 
\begin{align*} 
    &\left[\min_{x \in X \triangleq \cap_{i=1}^m X_i} \max_{y\in \fy{Y(x) \subseteq }\mathbb{R}^{\tilde{n}}} \  \frac{1}{m}\sum_{i=1}^m  \mathbb{E}_{\xi_i \in\mathcal{D}_i} [\, \tilde{f}_i(x,y,\xi_i)\,  ] \right] \quad =  \quad \min_{x \in X \triangleq \cap_{i=1}^m X_i}\frac{1}{m}\sum_{i=1}^m  \mathbb{E}_{\xi_i \in\mathcal{D}_i} [\, \tilde{f}_i(x,y(x),\xi_i)\,  ]
\end{align*}
where we equivalently write $y(x)= \mbox{arg}\min_{y\in \fy{Y(x)  }} \frac{1}{m}\sum_{i=1}^m  \mathbb{E}_{\zeta_i \in\tilde{\mathcal{D}}_i}[-\tilde{f}_i(x,y,\xi_i) ]$. Thus,  \eqref{prob:main_minimax} is an instance of \eqref{prob:main_BL}  in that we set $\tilde{h}_i: = -\tilde{f}_i $ and $\mathcal{D}_i:=\tilde{\mathcal{D}}_i$. Accordingly, we have the following result. 
\begin{corollary}\label{thm:minimax}\em 
Consider Algorithm~\ref{alg:FedRiZO_upper} for solving \eqref{prob:main_minimax}. Let Assumption~\ref{assum:bilevel} hold {for $\tilde{h}_i: = -\tilde{f}_i $ and $\mathcal{D}_i:=\tilde{\mathcal{D}}_i$}. Then, all the results in Theorem~\ref{thm:bilevel} hold true.
\end{corollary}
 
\section{Federated two-stage SMPECs}\label{sec:2s}
In this section, we address the two-stage stochastic mathematical program
with equilibrium constraints, defined as ~\eqref{prob:fed2smpecs}.  In
contrast with standard models, in this setting, the data is distributed among
$m$ clients. The local loss function of client $i$ is endogenized in the
random solution map  $y(\bullet,\xi_i)$; in fact, $y(x,\fyy{\xi_i})$ is the
unique solution to the parametric variational inequality (VI) problem, given as
$ \text{VI}(Y(x, \xi_i), G(x, \bullet , \xi_i))$.  Recalling the definition of a
VI,  $y(x,\xi_i)$ is a vector in $Y(x, \xi_i)$ such that $G(x, y(x,\xi_i) ,
\xi_i)^{\fyy{\top}}(y-y(x,\xi_i)) \geq 0$ for all $y \in Y(x, \xi_i)$.  To address
2s-SMPECs in a federated setting, we consider the following
implicit problem. 
\begin{definition}[\underline{Unconstrained smoothed implicit problem}] 
Given \eqref{prob:fed2smpecs} and $\eta > 0$, consider an unconstrained smoothed {implicit} problem given as
\begin{align}
&\min_{x \in \mathbb{R}^n} \ {\bf f}^{\eta}(x) \left\{\triangleq \tfrac{1}{m}\textstyle\sum_{i=1}^m \left(\mathbb{E}_{\xi_i,u\in\mathbb{B}}[{\tilde{f}_i}\left(x+\eta u,y(x+\eta u,\xi_i),\xi_i\right)] +\tfrac{1}{2\eta}\mbox{dist}^2(x,X_i) \right)\right\} .\label{prob:main_unconstrained_2s}
\end{align}
\end{definition}
\begin{remark}[Distinctions with bilevel FL]\label{rem:bl_vs_2s}
We note that the implicit problem of~\eqref{prob:fed2smpecs} given by~\eqref{prob:main_unconstrained_2s} is distinct from that of ~\eqref{prob:main_BL} given by~\eqref{prob:main_unconstrained_bilevel}.  A major difference is that in the former,  the solution to the lower-level problem is independent of $\xi_i$ while in the latter,  $y(x,\xi_i)$ is characterized by $\xi_i$.  Consequently,  FedRZO$_{\texttt {bl}}$ is not suitable for addressing 2s-SMPECs.  In fact,  unlike in the bilevel problem~\eqref{prob:main_BL}  where the lower-level problem is itself a federated optimization problem,  the lower-level problem in~\eqref{prob:fed2smpecs} is a  non-federated (deterministic) VI problem parametrized by $\xi_i$.  In view of these observations, to develop an inexact zeroth-order scheme for solving 2s-SMPECs. \fyy{One} may consider employing a standard iterative solver available, \us{that leverages} gradient or extragradient methods for VIs with strongly monotone maps. $\hfill$ $\Box$  
\end{remark} 

\begin{algorithm}[htb]                                                                                                                                 
	\caption{Two-Stage Randomized Implicit Zeroth-Order Federated Averaging (FedRZO$_{\texttt {2s}}$)}\label{alg:fed2:upper}
\begin{algorithmic}[1] 
	 \State \textbf{input:} Server chooses a random $\hat x_0 \in X$, stepsize $\gamma$, smoothing parameter $\eta$, synchronization indices $T_0:=0$ and $T_r \geq 1$, where $r\geq 1$ is the upper-level communication round index
	 \For {$r = 0,1, \dots$}
	 \State Server broadcasts $\hat x_{r}$ to all clients; $x_{i,T_r}:= \hat x_r, \ \forall i$
	 \For {$k = T_r,  \ldots,T_{r+1}-1$} {\bf in parallel by {clients}}
	 \State Client $i$ generates \fyy{a} random \fyy{sample} $\xi_{i,k} \in \mathcal{D}_i$ and $v_{i,k}\in \eta\mathbb{S}$, and calls Algorithm \ref{alg:fed2:lower} twice to receive $y_{\varepsilon_k}(x_{i,k}+v_{i,k}, \xi_{i,k})$ and $y_{\varepsilon_k}(x_{i,k}\yqro{-v_{i,k}}, \xi_{i,k})$, {where $\|y_{\varepsilon_k}(\bullet,\xi)-y(\bullet,\xi)\|^2 \leq \varepsilon_k$}
	 \State  $g^{\eta,\varepsilon_k}_{i,k}: =\frac{n}{\yqro{2}\eta^2} ( \tilde{f}_i(x_{i,k}+v_{i,k},y_{\varepsilon_k}(x_{i,k}+v_{i,k}, \xi_{i,k}),\xi_{i,k})- \tilde{f}_i(x_{i,k}\yqro{-v_{i,k}},y_{\varepsilon_k}(x_{i,k}\yqro{-v_{i,k}}, \xi_{i,k}),\xi_{i,k}))v_{i,k}$
	 \State  Client $i$ does a local update as $x_{i,k+1} :=  x_{i,k} - \gamma \left(  g^{\eta,\varepsilon_k}_{i,k} + \frac{1}{\eta}\left(x_{i,k}-\mathcal{P}_{X_i}(x_{i,k}) \right) \right)    $

	 \EndFor	
	\State Server receives $x_{i,T_{r+1}}$ from all clients and aggregates, i.e., $\hat x_{r+1} := \frac{1}{m}\sum_{i=1}^m x_{i,T_{r+1}}$

	 \EndFor
	\end{algorithmic}
\end{algorithm}
The main assumptions in this section are stated in the following. 
\begin{assumption}\label{assum:fed2smpecs} \em
Consider problem \eqref{prob:fed2smpecs}. Let the following hold.
 

\noindent (i) For all $i \in [m]$, $\tilde
    f_i(\bullet,y(\bullet,\xi_i),\xi_i)$ is $L_0^{\tiny
    \mbox{imp}}(\xi_i)$-Lipschitz and ${\tilde{f}_i}(x,\bullet,\xi_i)$ is
    $L_{0,y}^f(\xi_i)$-Lipschitz for any $x$, where $L_0^{\tiny
    \mbox{imp}}\triangleq \displaystyle{\max_{i=1,\ldots,m}}\sqrt{\mathbb{E}_{\yqro{\xi_i\in \mathcal D_i}}[(
    L_0^{\tiny \mbox{imp}}(\xi_i))^2]}<\infty$ and $L_{0,y}^f \triangleq
    \displaystyle{\max_{i=1,\ldots,m}}\sqrt{\mathbb{E}_{\yqro{\xi_i\in \mathcal D_i}}[({L_{0,y}^f}(\xi_i))^2]}<\infty$.


\noindent (ii) For all $i \in [m]$, $G(x, \bullet, \xi_i)$ is a $\mu_F(\xi_i)$-strongly monotone and $L_F(\xi_i)$-Lipschitz continues mapping on $\mathbb{R}^{\tilde n}$ uniformly in $x\in X$, where there exist some positive scalars $\mu_F\leq \displaystyle{\inf _{\xi_i \in \mathcal{D}_i}}\mu_F(\xi_i)$ and $L_F\geq \displaystyle{\sup _{\xi_i \in \mathcal{D}_i}}L_F(\xi_i)$. 

\noindent (iii) The sets $X_i$, for all $i \in [m]$, satisfy Assumption~\ref{assum:Heterogeneous} (iii).

\noindent (iv) For any $x \in  \mathbb{R}^{n}$ and all $\xi_i \in \mathcal{D}_i$ where $i \in [m]$,  the set $Y(x,\xi_i)\subseteq \mathbb{R}^{\tilde n}$ is nonempty,  compact, and convex.

\end{assumption}

\begin{algorithm}[htb]
	\caption{Projection method for the lower-level VI in \eqref{prob:fed2smpecs}}\label{alg:fed2:lower}                                                                                                                               
\begin{algorithmic}[1]
    \State \textbf{input:} An arbitrary $y_0$, integer $k {> 0}$, vector $\hat x_k$, random variable $\xi$, stepsize $\alpha > 0$, scalar $\tau >0$
       \State Set $t_k = \lceil\tau \text{ln}(k+1)\rceil$
	 \For {$t = 0,1, \dots, t_k -1$}
	 \State Evaluate $G(\hat x_k, y_t , \xi)$ an update $y_t$ as $y_{t+1} = \mathcal{P}_{Y(x,\xi)}\left(y_t - \alpha G(\hat x_k, y_t , \xi)\right)$
	 \EndFor	
       \State Return $y_{t_k}$
	\end{algorithmic}
\end{algorithm}
In the following theorem, we present the main convergence result and complexity bounds for FedRZO$_{\texttt {2s}}$.
\begin{theorem}[Rate and complexity statements for \eqref{prob:fed2smpecs}]\label{thm:fed2smpec}\em  Consider Algorithms \ref{alg:fed2:upper} and \ref{alg:fed2:lower} for solving \eqref{prob:fed2smpecs}. Let Assumption \ref{assum:fed2smpecs} hold.

\noindent (i) Consider Algorithm~\ref{alg:fed2:lower} being called by client $i$,  given $\hat x_{i,k}$ and $\xi_{i,k}$. Suppose the stepsize in Algorithm~\ref{alg:fed2:lower} be given as $\alpha := \tfrac{ \mu_F}{L_F^2}$.  Let $y_{\varepsilon_{k}}(\hat x_k, \xi_{i,k})$ denote the output of Algorithm~\ref{alg:fed2:lower} and $y(\hat x_k, \xi_{i,k})$ denote the solution to the lower-level VI in~\eqref{prob:fed2smpecs}. Then, $\|y_{\varepsilon_{k}}(\hat x_k, \xi_{i,k})-y(\hat x_k, \xi_{i,k})\|^2\leq  \varepsilon_k$ where $\varepsilon_k=B (1-1/\kappa_F^2)^{t_k}$, $\kappa_F\triangleq \tfrac{L_F}{\mu_F}$,  $t_k$ denotes the number of iterations used in Algorithm~\ref{alg:fed2:lower}, and $B$ is an upper bound on $\|y(\hat x_k, \xi_{i,k})-y_0\|^2$. 

\noindent (ii) {\bf [Error bound]} Let $\gamma \leq \min\{\tfrac{1}{12\sqrt{6}HB_2\sqrt{\Theta_2}}, \yqro{\tfrac{\eta}{16(L_0^{\tiny \mbox{imp}}\sqrt n+1)}}\}$ and $k^*$ be chosen uniformly at random from $0,\ldots,K:=T_R-1$, we have
\begin{align*} 
\mathbb{E}\left[\|\nabla {\bf f}^{\eta}(\bar x_{k^*})\|^2\right]   & \leq   8(\gamma T_R)^{-1}(\mathbb{E}\left[{\bf f}^{\eta}(x_{0}) \right] -  {\bf f}^{\eta,*}  )   +\tfrac{8\gamma \Theta_1}{m} +8H^2\gamma^2\Theta_2\Theta_5\\
 & +\left(8(T_R)^{-1}H^2\gamma^2\Theta_2\Theta_4  + 8(T_R)^{-1}\Theta_3\right)\sum_{k=0}^{T_R-1}\varepsilon_k.
\end{align*}
where $\Theta_1 :=\yqro{\tfrac{16\sqrt{2\pi}(L_0^{\text{imp}}\sqrt n+1)(L_0^{\text{imp}})^2n}{m\eta}}$, $\Theta_2:=\tfrac{5\yqro{(L_0^{\text{imp}}\sqrt n+1)^2}}{8\eta^2}$, ${\Theta_3:= \tfrac{\yqro{5} n^2}{\eta^2}(L_{0,y}^f)^2}$, ${\Theta_4:=\tfrac{\yqro{9}n^2}{\eta^2}(L_{0,y}^f)^2}$ and $\Theta_5:={\yqro{18}\left(\yqro{\tfrac{B_1^2}{\eta^2}} + (4B_2^2+1)\yqro{8\sqrt{2\pi}(L_0^{\text{imp}})^2n} \right)}$.

\noindent(iii) {\bf  [Iteration complexity]} Let $t_k= \lceil\tau \text{ln}(k+1)\rceil$ where $\tau\geq\tfrac{-1}{\text{ln}(1-1/\kappa^2_F)}$, $\gamma:=\sqrt{\tfrac{m}{K}}$, and $H:=   \left\lceil\sqrt[4]{\tfrac{K}{m^3}}\right\rceil$. Let $\epsilon>0$ be an arbitrary scalar and $K_\epsilon$ denote the number of iterations such that $\mathbb{E}\left[\|\nabla {\bf f}^\eta(\bar x_{k^*})\|^2\right] \leq \epsilon$. Then, the iteration complexity of FedRZO$_{\texttt {2s}}$ (upper level) is  \yqro{$$K_\epsilon:=\mathcal{\tilde{O}}\left( \tfrac{\yqro{\tfrac{n^{1.5}(L_0^{\text{imp}})^{3}}{\eta}+\tfrac{n(L_0^{\text{imp}})^{-2}}{\eta^{2}}}}{m\epsilon^2} + \tfrac{\yqro{\tfrac{n(L_0^{\text{imp}})^{-2}}{\eta^{-2/3}}}}{m^{1/3}\epsilon^{2/3}} + \tfrac{\yqro{n^2}}{\epsilon}\right).$$}  

\noindent (iv)  {\bf  [Overall communication complexity and projections]}  Suppose $K_\epsilon \geq m^3$.  The overall communication complexity is (upper and lower levels) is $R = \mathcal{O}\left( (mK)^{\tfrac{3}{4}}\right)$.  The total number of lower-level projection steps is $\tilde{\mathcal{O}}\left(mK /\text{ln}\left(\tfrac{1}{1- \kappa^{-2}_F}\right) \right).$
\end{theorem}

\begin{remark}
 (i) The  iteration complexity of FedRZO$_{\texttt {2s}}$ (upper level) is given as $K_\epsilon:=\mathcal{\tilde{O}}\left( \frac{1}{m\epsilon^2} + \frac{1}{m^{1/3}\epsilon^{2/3}} + \frac{1}{\epsilon}\right).$ This is indeed comparable with that of the convergence rate of FedRZO$_{\texttt {nn}}$ given as $\mathcal{O}\left(\tfrac{1}{ {m\epsilon^2}}\right)$ (see  Remark~\ref{rem:fedrzonn}). Indeed, if $\epsilon \leq \tfrac{1}{m}$, the iteration complexity of FedRZO$_{\texttt {2s}}$ is
$K_\epsilon:=\mathcal{\tilde{O}}\left( \frac{1}{m\epsilon^2} + \frac{1}{m\epsilon^2} + \frac{1}{m\epsilon^2}\right)$, matching with the iteration complexity of FedRZO$_{\texttt {nn}}$.   (ii) The upper-level communication complexity, given as $\mathcal{O}\left( (mK)^{\tfrac{3}{4}}\right)$, is also the total communication complexity of FedRZO$_{\texttt {2s}}$. This is because at any iteration $k$ of Algorithm~\ref{alg:fed2:upper},  implementing Algorithm~\ref{alg:fed2:lower} by each client does not require any communications with other clients. Further, unlike FedRZO$_{\texttt {bl}}$, FedRZO$_{\texttt {2s}}$ does not require any additional communications during the global step. Consequently, the overall communication overhead of FedRZO$_{\texttt {2s}}$ is significantly smaller than that of FedRZO$_{\texttt {bl}}$ (see Remark~\ref{rem:thm2} (iii) for further details about FedRZO$_{\texttt {bl}}$). $\hfill$ $\Box$ 
\end{remark}
\section{\us{Numerical} experiments}\label{sec:num} 
We present four sets of experiments to validate the performance of the
proposed algorithms. In Section~\ref{sec:NNs}, 
we implement
FedRZO$_{\texttt {nn}}$ on ReLU neural networks (NNs) and compare it with some
recent FL methods. In Sections~\ref{sec:hyper} and~\ref{sec:fair},
we implement FedRZO$_{\texttt {bl}}$ on federated hyperparameter
learning and a minimax
formulation in FL.  Lastly,  in  Section~\ref{sec:num2s} we implement FedRZO$_{\texttt {2s}}$ for solving an instance of a two-stage SMPEC. Throughout, we use the MNIST dataset.  In Section~\ref{sec:NNs}, we also present additional experiments on a higher dimensional dataset (i.e., CIFAR-10). \yqro{\us{We} simulate the non-iid setting with Dirichlet distribution $Dir(\alpha)$, which is commonly used in FL experiments (e.g. \cite{li2022federated,wu2023personalized,xiong2023feddm}) and $\alpha$ \us{denotes} the concentration parameter used to determine the non-iid level. We set $\alpha=0.5$.}
\begin{table}[H]
\setlength{\tabcolsep}{0pt}
\centering{
\begin{tabular}{c || c  c  c}
{\footnotesize {Setting}\ \ }& {\footnotesize  (S1) with $\eta =0.1$} & {\footnotesize (S2) with $\eta =0.01$} & {\footnotesize (S3) with $\eta =0.001$}\\
\hline\\
\rotatebox[origin=c]{90}{{\footnotesize {Log of average loss}}}
&
\begin{minipage}{.3\textwidth}
\includegraphics[scale=.21, angle=0]{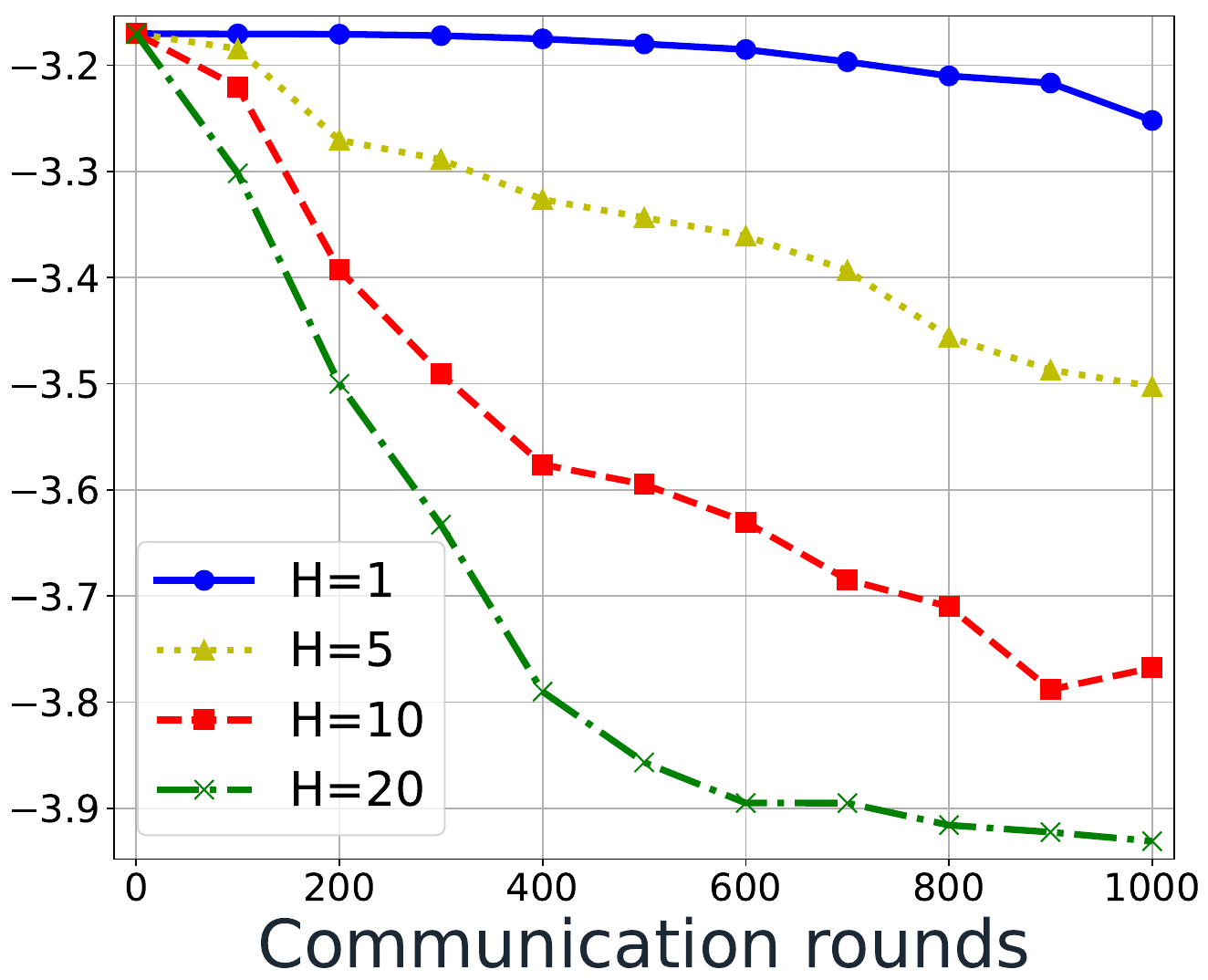}
\end{minipage}
&
\begin{minipage}{.3\textwidth}
\includegraphics[scale=.21, angle=0]{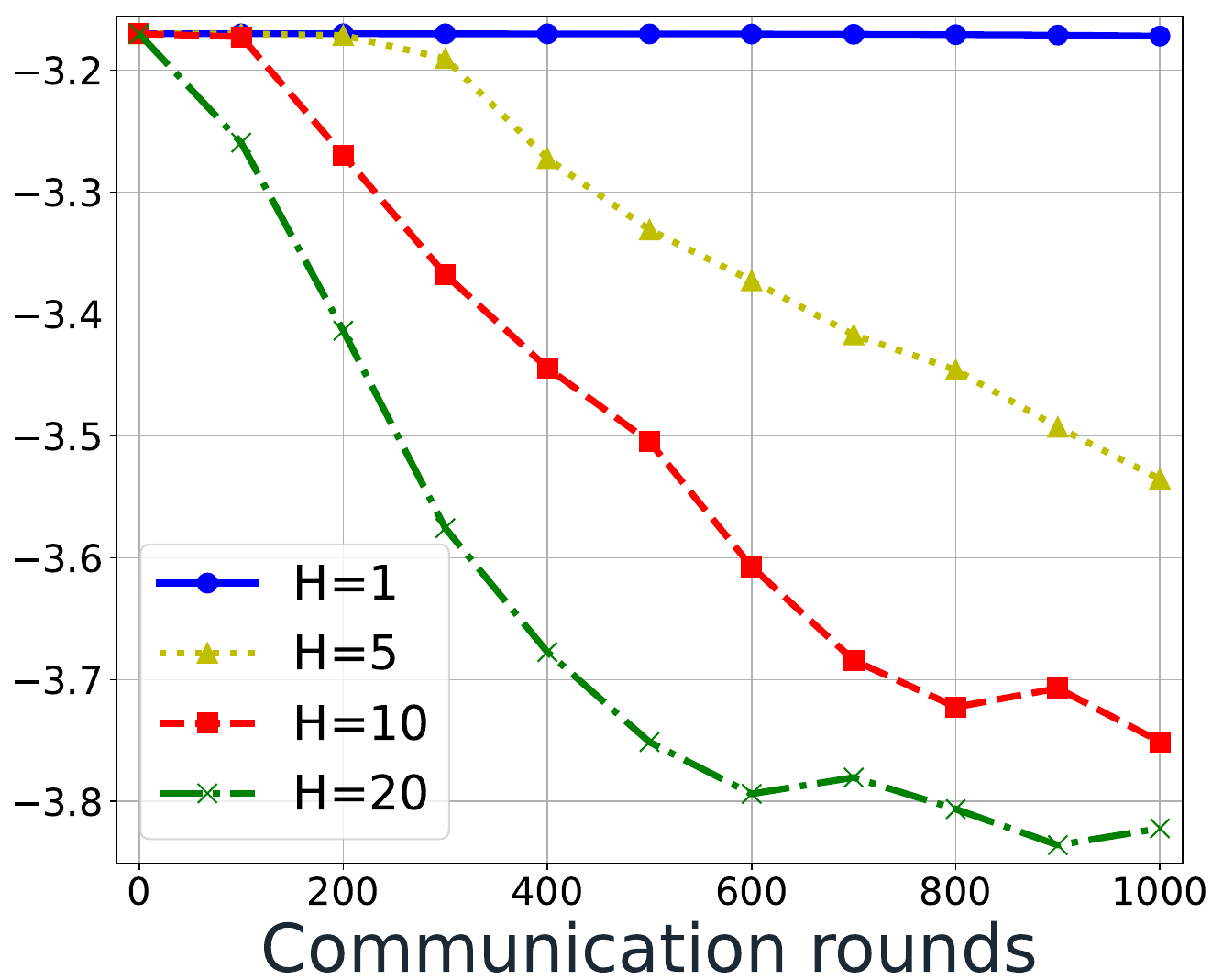}
\end{minipage}
&
\begin{minipage}{.3\textwidth}
\includegraphics[scale=.21, angle=0]{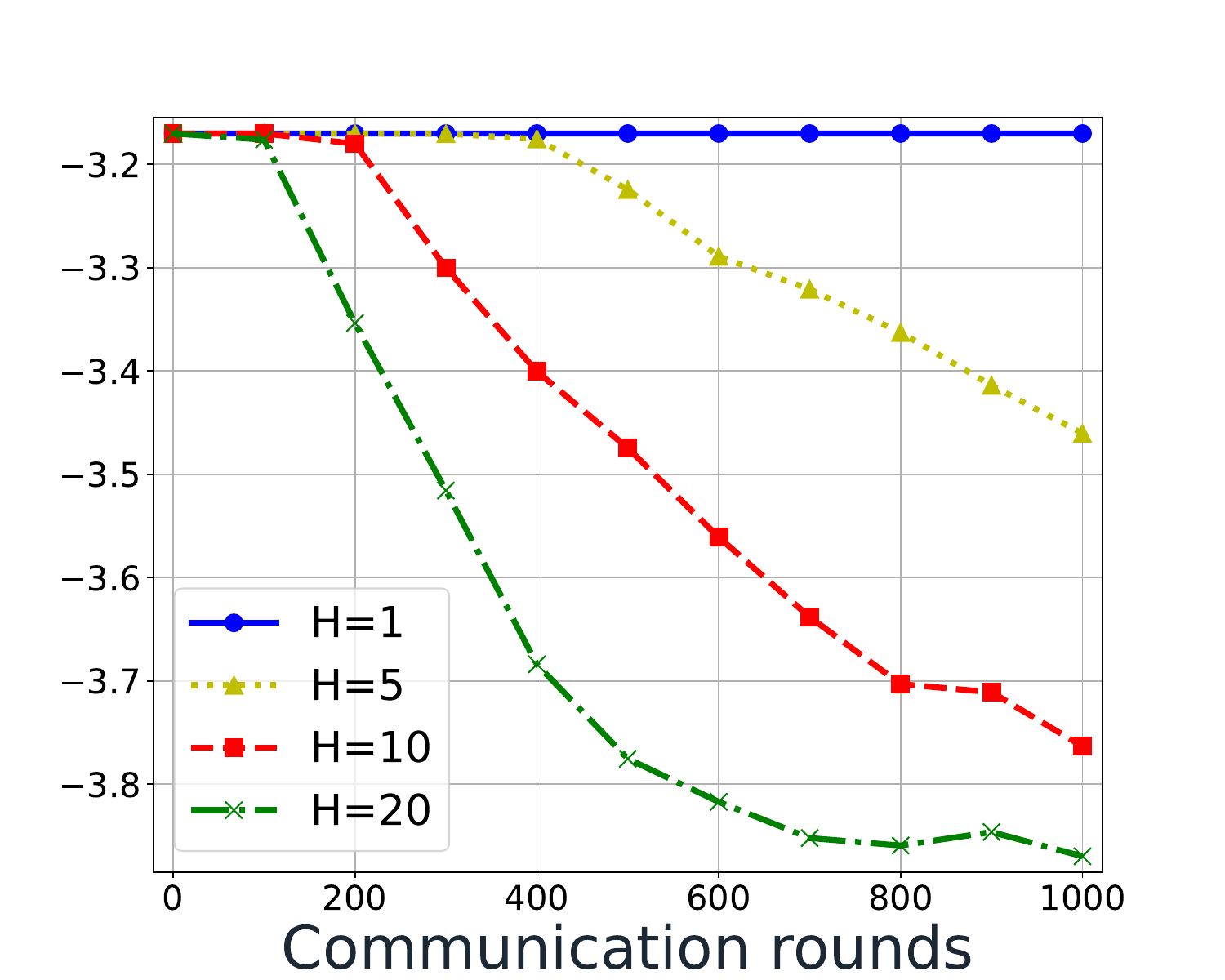}
\end{minipage}
\\
\end{tabular}}
\vspace{0.1in}
\captionof{figure}{\yqro{Performance of FedRZO$_{\texttt {nn}}$ on a two-layer ReLU NN regarding communication rounds for different no. of local steps and different values of the smoothing parameter $\eta$. FedRZO$_{\texttt {nn}}$ benefits from a larger number of local steps and shows robustness with respect to the choice of \fyy{$\eta$.}}}
\label{fig1:NN:cr:iter}
\vspace{-.2in}
\end{table}
\subsection{Federated training of ReLU {NNs}. }\label{sec:NNs}
We implement FedRZO$_{\texttt {nn}}$ for federated training in a \yqro{two-layer (one hidden layer)} ReLU {NN} with $N_1$ neurons. This is a nondifferentiable nonconvex optimization problem, aligning with \eqref{prob:main} and taking the form 
$\min_{x:=(Z,w) \in \mathcal{X}} \quad \frac{1}{2m}\sum_{i=1}^m\sum_{\ell\in \yqro{\mathcal{S}}_i}  (v_{i,\ell}-\sum_{q=1}^{N_1}w_q \sigma(Z_{\bullet,q}U_{i,\ell}))^2+\frac{\lambda}{2}\left( \|Z\|_F^2+\|w\|^2 \right),$
where $m$ denotes the number of clients, $Z\in\mathbb{R}^{N_1\times N_0}$, $w\in\mathbb{R}^{N_1}$, $N_0$ is the feature dimension, \yqro{$\mathcal D_i \triangleq \{(U_{i,\ell},v_{i,\ell})\in \mathbb R^{N_0}\times \{-1,1\} \mid \ell\in \mathcal S_i\}$ denotes the local dataset of client $i$ and $\mathcal S_i=\{1,\dots,|\mathcal D_i|\}$ denotes the local index set,} respectively,  $\sigma(x):=\max \{0,x\}$, and $\lambda$ is the regularization parameter.

{\bf Setup.} We distribute the training dataset among $m:=5$ clients and implement FedRZO$_{\texttt {nn}}$ for the FL training with $N_1:=4$ neurons under three different settings for the smoothing parameter, $\eta \in \{0.1, 0.01, 0.001\}$, $\gamma:=10^{-5}$, and $\lambda:=0.01$. We study the performance of the method under different number of local steps {with} $H \in \{1,5,10,20\}$.

\textbf{Results and insights.} Figure~\ref{fig1:NN:cr:iter} presents the first
set of {numerics} for FedRZO$_{\texttt {nn}}$ under the aforementioned
settings. In terms of communication rounds, {we observe that} the
performance of the method improves by using a larger number of local steps. In
fact, {in} the case where $H:=1$, FedRZO$_{\texttt {nn}}$ is equivalent to a
parallel zeroth-order SGD that employs communication among clients at each
iteration, {resulting} in a poor performance, motivating the need for the FL
framework. In terms of $\eta$, while we observe robustness of the scheme in
terms of the original loss function, we also {note a} slight improvement in
the {empirical} speed of convergence in early steps, as $\eta$ increases.
This is indeed aligned with the dependence of convergence bound in
Proposition~\ref{thm:main_bound} on $\eta$. 

\textbf{Comparison with other FL methods.} While \fy{we are} unaware of other FL
methods for addressing nondifferentiable nonconvex problems, we compare
FedRZO$_{\texttt {nn}}$ with other FL methods including FedAvg~\cite{FedAvg17},
FedProx~\cite{li2020federated}, FedMSPP~\cite{yuan2022convergence}, and
Scaffnew~\cite{pmlr-v162-mishchenko22b} when {applied} on a {NN} with a smooth rectifier.   
Throughout this experiment, 
$\eta:=0.01$ for FedRZO$_{\texttt {nn}}$. For the other methods, we use the
softplus rectifier defined as $\sigma^{\beta}(x):=\beta^{-1}{\ln(1+\exp(\beta
x))} $ as a smoothed approximation of the ReLU function. Our goal lies in
observing the sensitivity of the methods to $\beta$ when we compare them in terms of the original ReLU loss function.
Figure~\ref{fig4:FedZo:gradeval:cr} presents the results. We observe that the
performance of the standard FL methods is sensitive to the choice of $\beta$.
It is also worth noting that the smooth approximation of activation functions,
such as softplus, hypertangent, or sigmoid, involves the computation of an
exponential term which might render some scaling issues in large-scale
datasets. 
\begin{table}[htb]
\setlength{\tabcolsep}{0pt}
\centering{
\begin{tabular}{c || c  c  c}
{\footnotesize {Setting}\ \ }& {\footnotesize $\beta=1$} & {\footnotesize $\beta=10$} & {\footnotesize $\beta=50$ }\\
\hline\\
\rotatebox[origin=c]{90}{{\footnotesize {Log of average loss}}}
&
\begin{minipage}{.3\textwidth}
\includegraphics[scale=.21, angle=0]{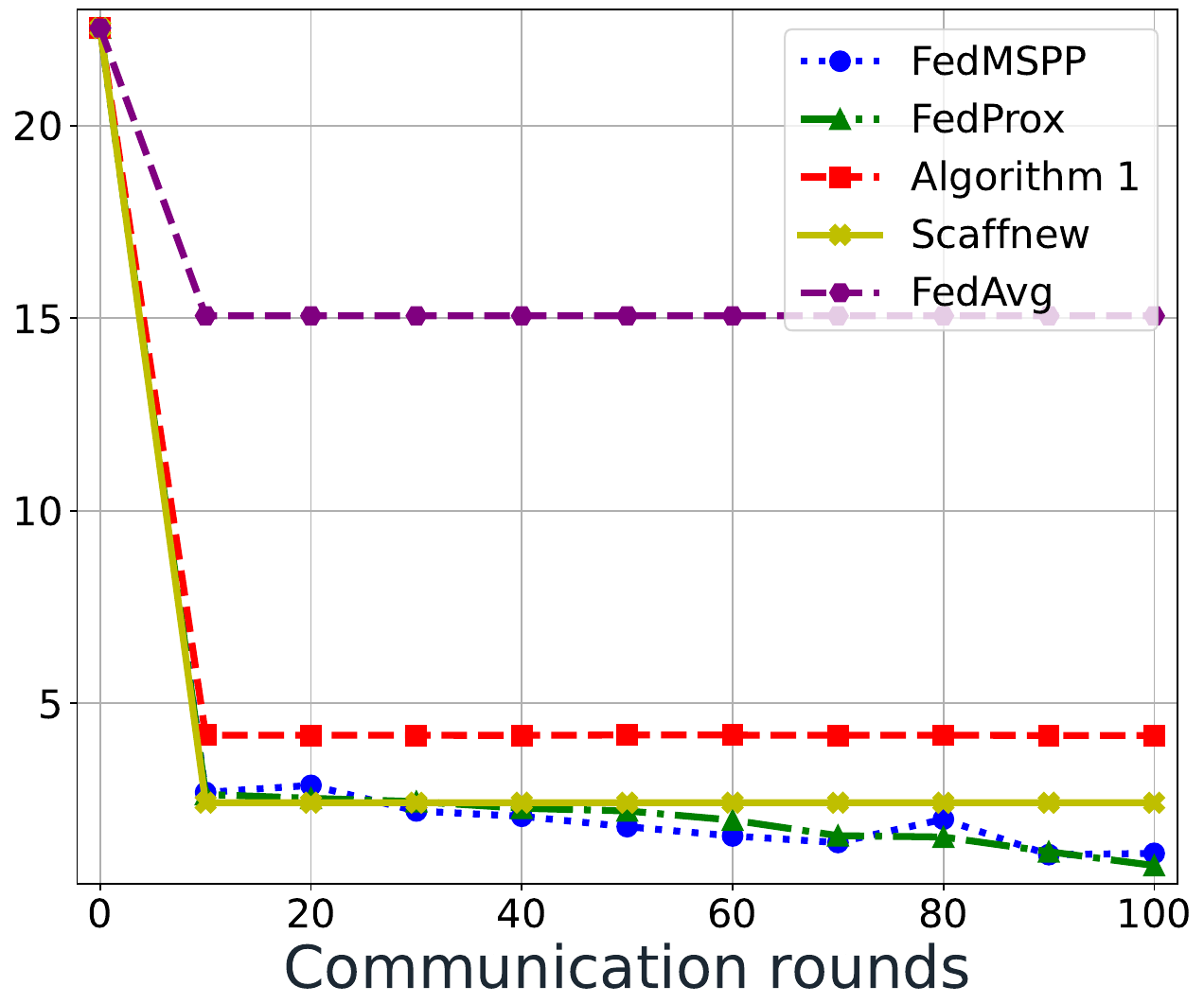}
\end{minipage}
&
\begin{minipage}{.3\textwidth}
\includegraphics[scale=.21, angle=0]{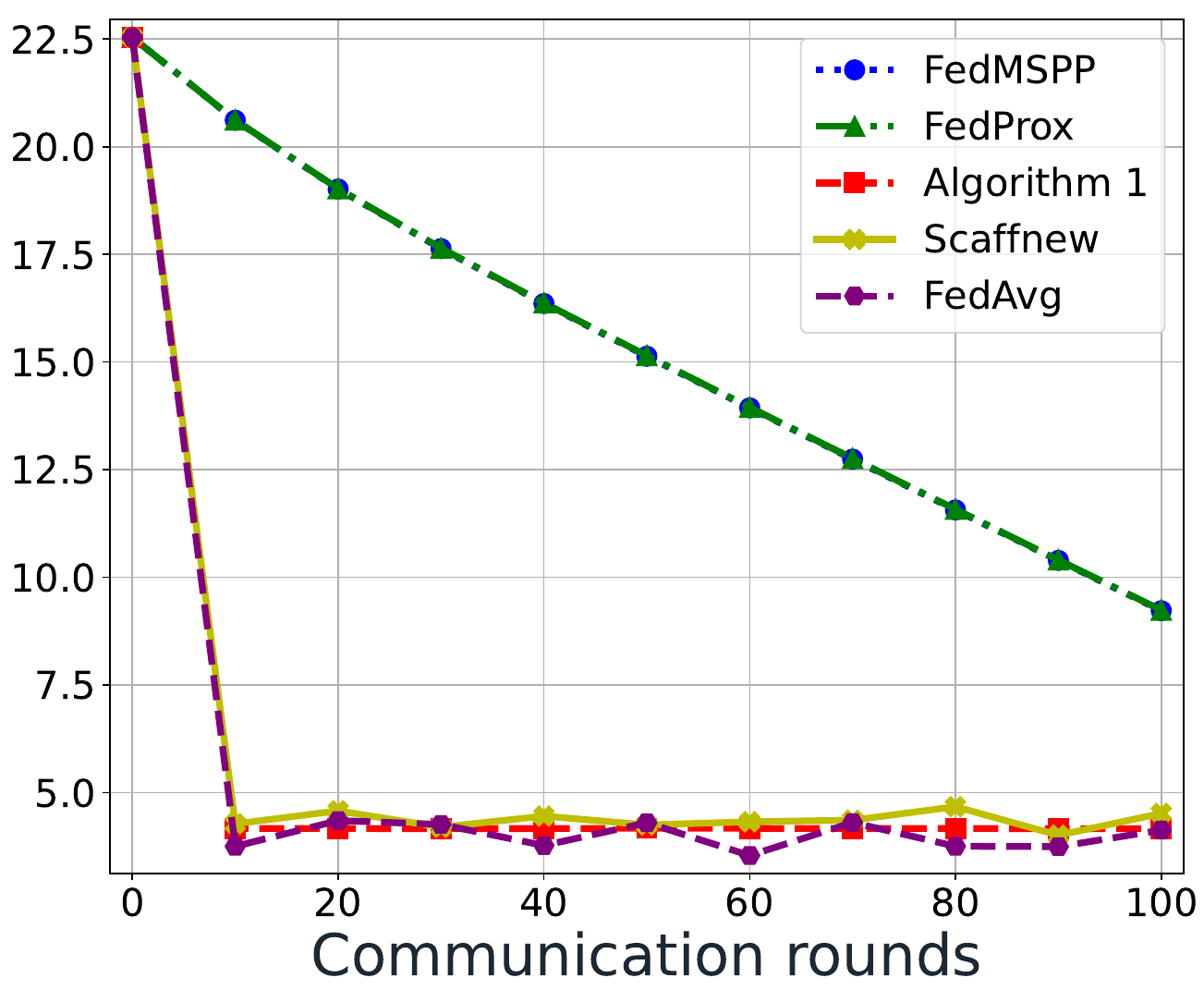}
\end{minipage}
&
\begin{minipage}{.3\textwidth}
\includegraphics[scale=.21, angle=0]{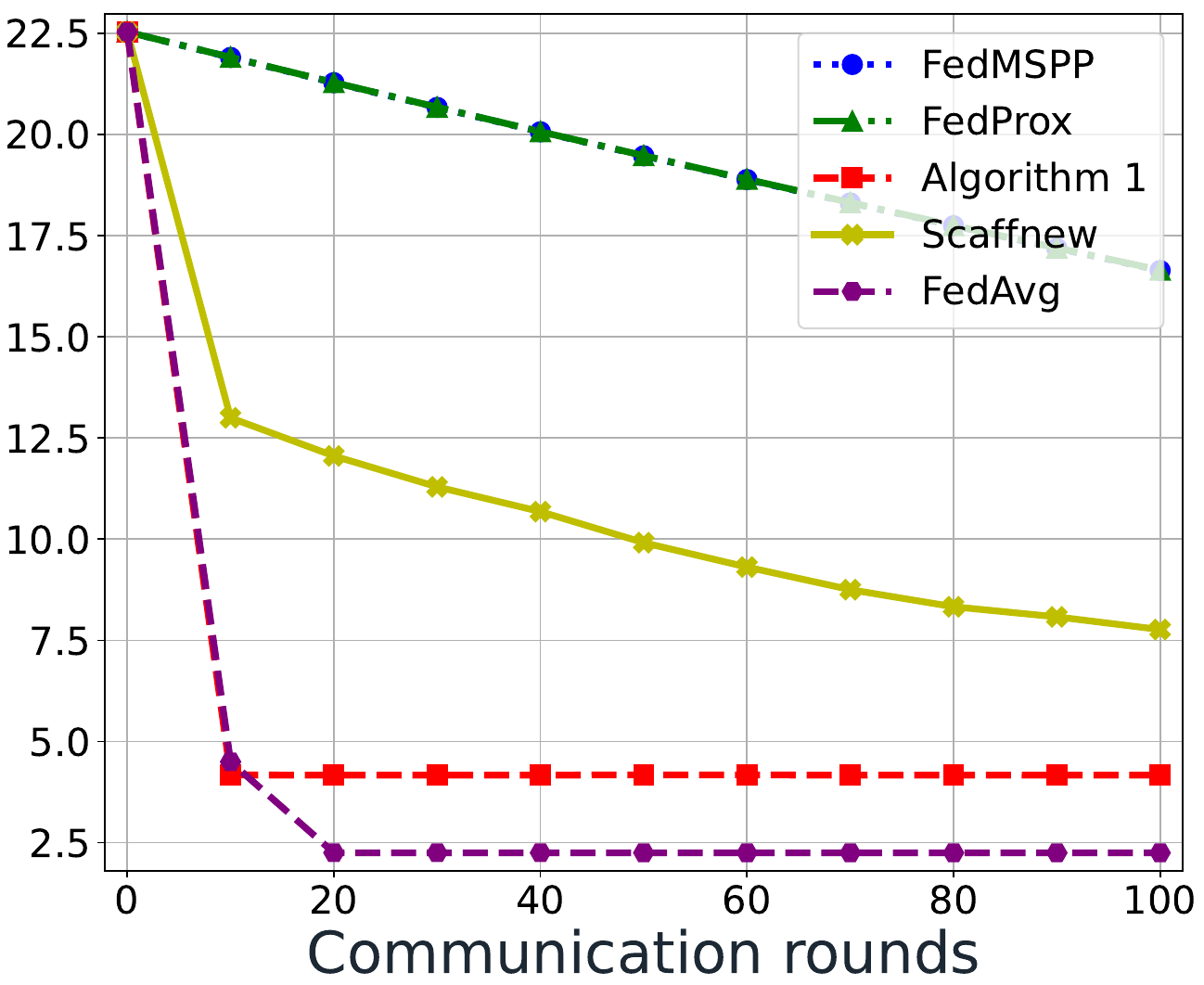}
\end{minipage}

\\
\hbox{}& & &\\
\hline\\
\rotatebox[origin=c]{90}{{\footnotesize {Log of average loss}}}
&
\begin{minipage}{.3\textwidth}
\includegraphics[scale=.21, angle=0]{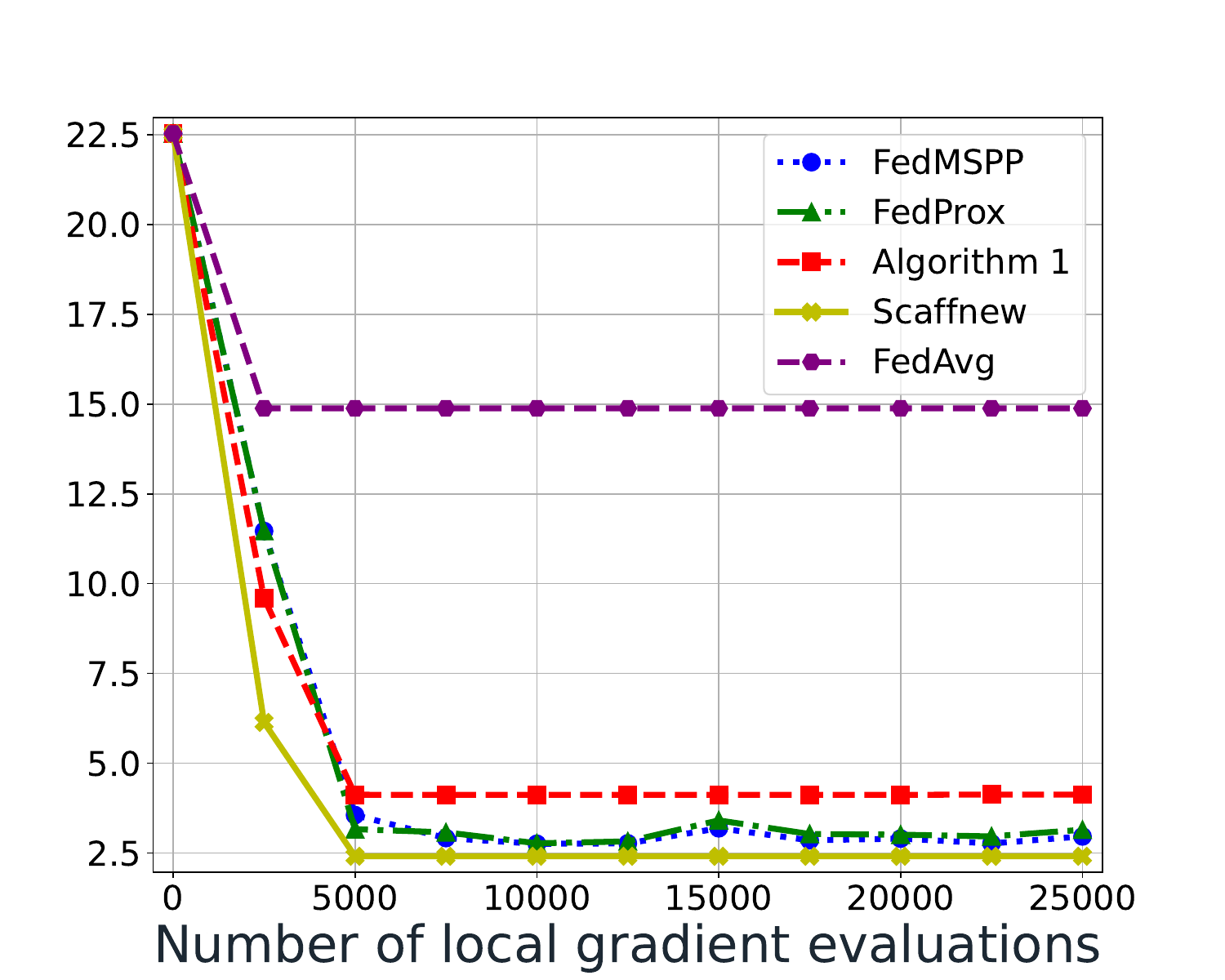}
\end{minipage}
&
\begin{minipage}{.3\textwidth}
\includegraphics[scale=.21, angle=0]{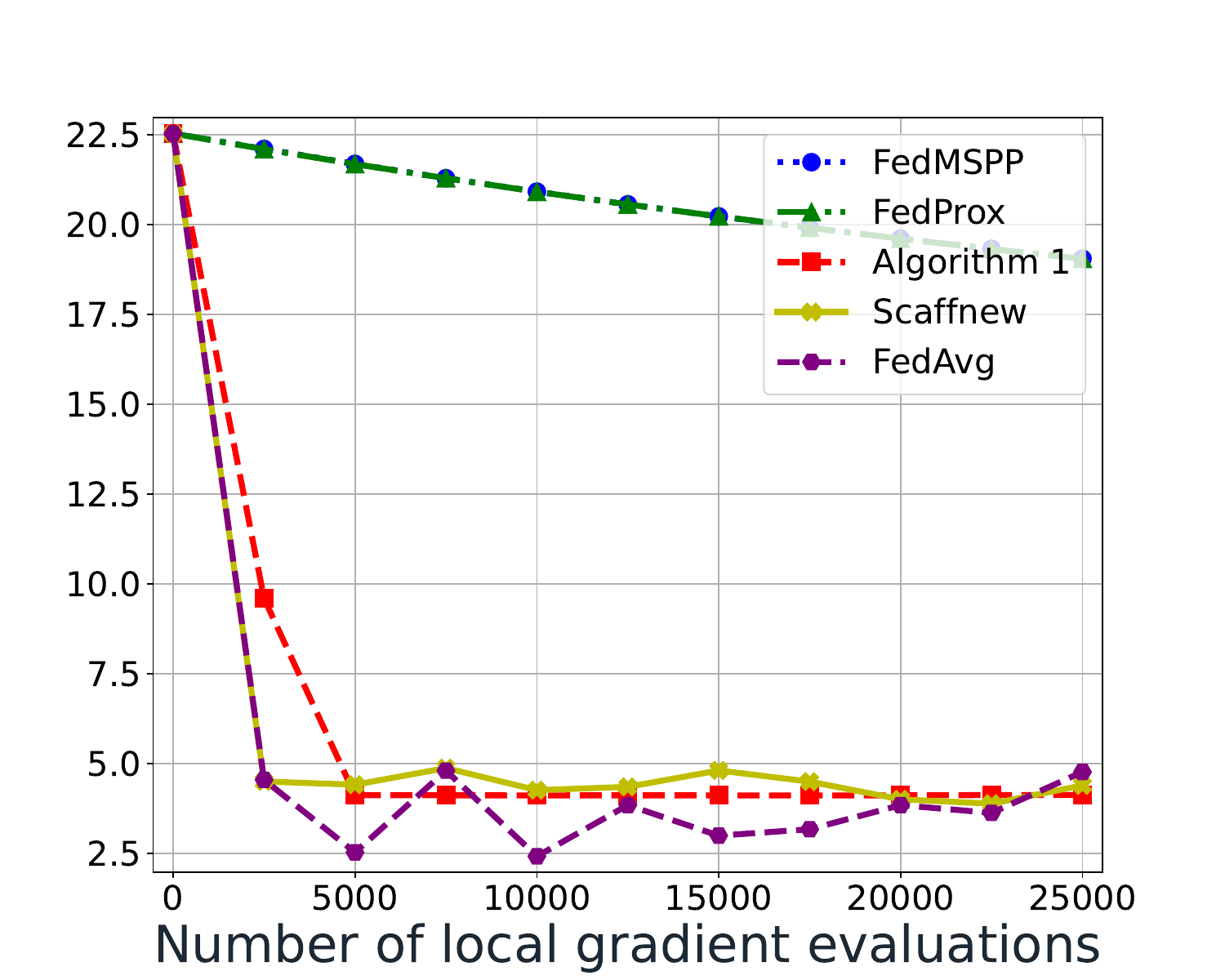}
\end{minipage}
&
\begin{minipage}{.3\textwidth}
\includegraphics[scale=.21, angle=0]{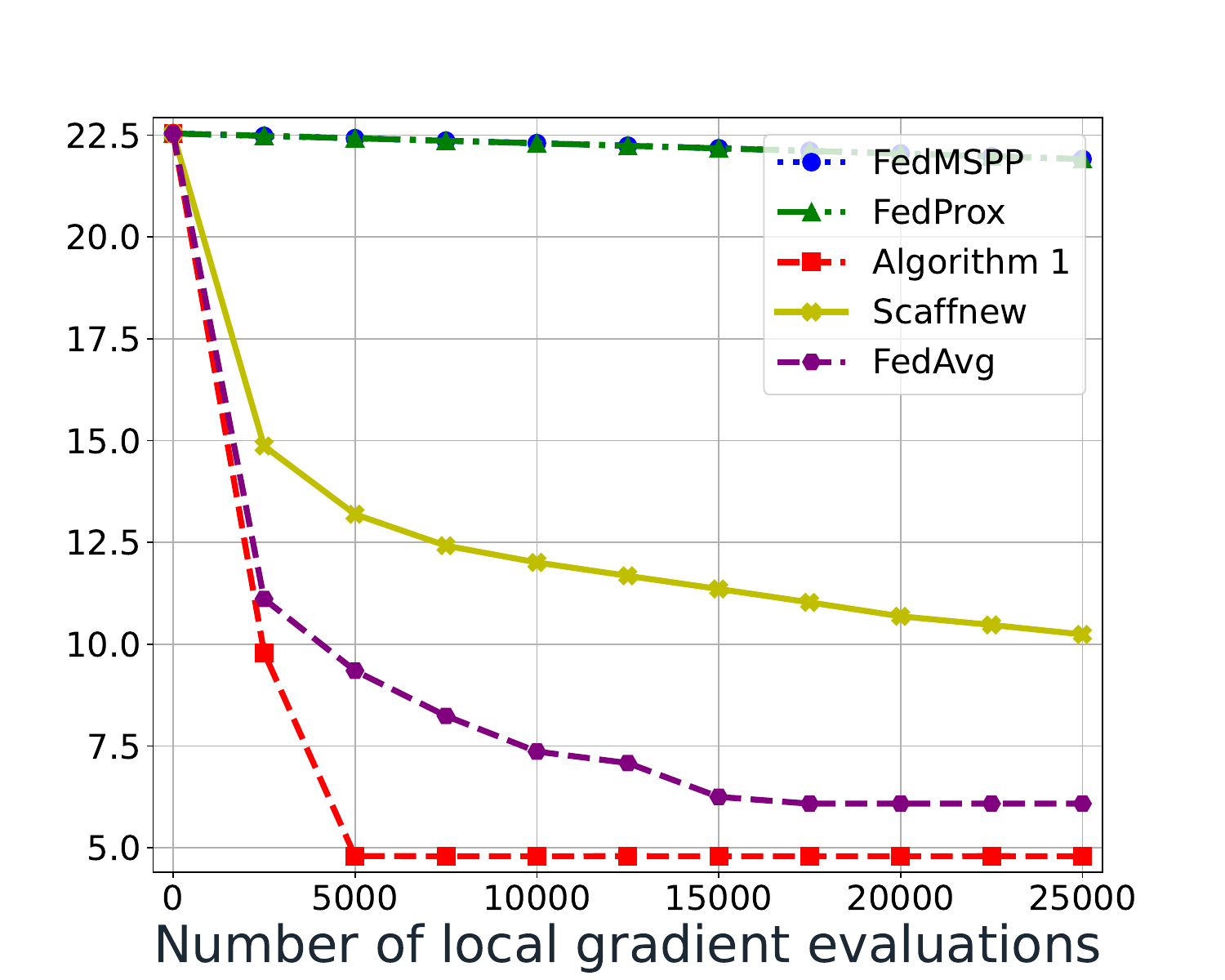}
\end{minipage}
\end{tabular}}
\vspace{.1in}
\captionof{figure}{Comparison between FedRZO$_{\texttt {nn}}$ on ReLU NN and standard FL methods when they are implemented on an NN with a smoothed variant of ReLU,  characterized by a parameter $\beta$. The performance of the standard FL methods appears to be sensitive to the choice of $\beta$.}
\label{fig4:FedZo:gradeval:cr}
\vspace{-.3in}
\end{table}

We also implement Algorithm~\ref{alg:PZO_LSGD} on the CIFAR-10 dataset to test its performance on problems with higher dimensions. We use the same objectives as defined in Section~\ref{sec:NNs}. We set the number of neurons in the layer to be 4, 20, and 100 in three sets of experiments, \us{respectively}. The results are presented in Figure \ref{fig:cifar10:sup}. 
\begin{table}[htb]
\setlength{\tabcolsep}{0pt}
\centering{
\begin{tabular}{c || c  c  c}
{\footnotesize {Setting}\ \ }& {\footnotesize $\eta=0.1$} & {\footnotesize $\eta=0.01$} & {\footnotesize $\eta=0.001$ }\\
\hline\\
\rotatebox[origin=c]{90}{{\footnotesize {4 neurons}}}
&
\begin{minipage}{.3\textwidth}
\includegraphics[scale=.21, angle=0]{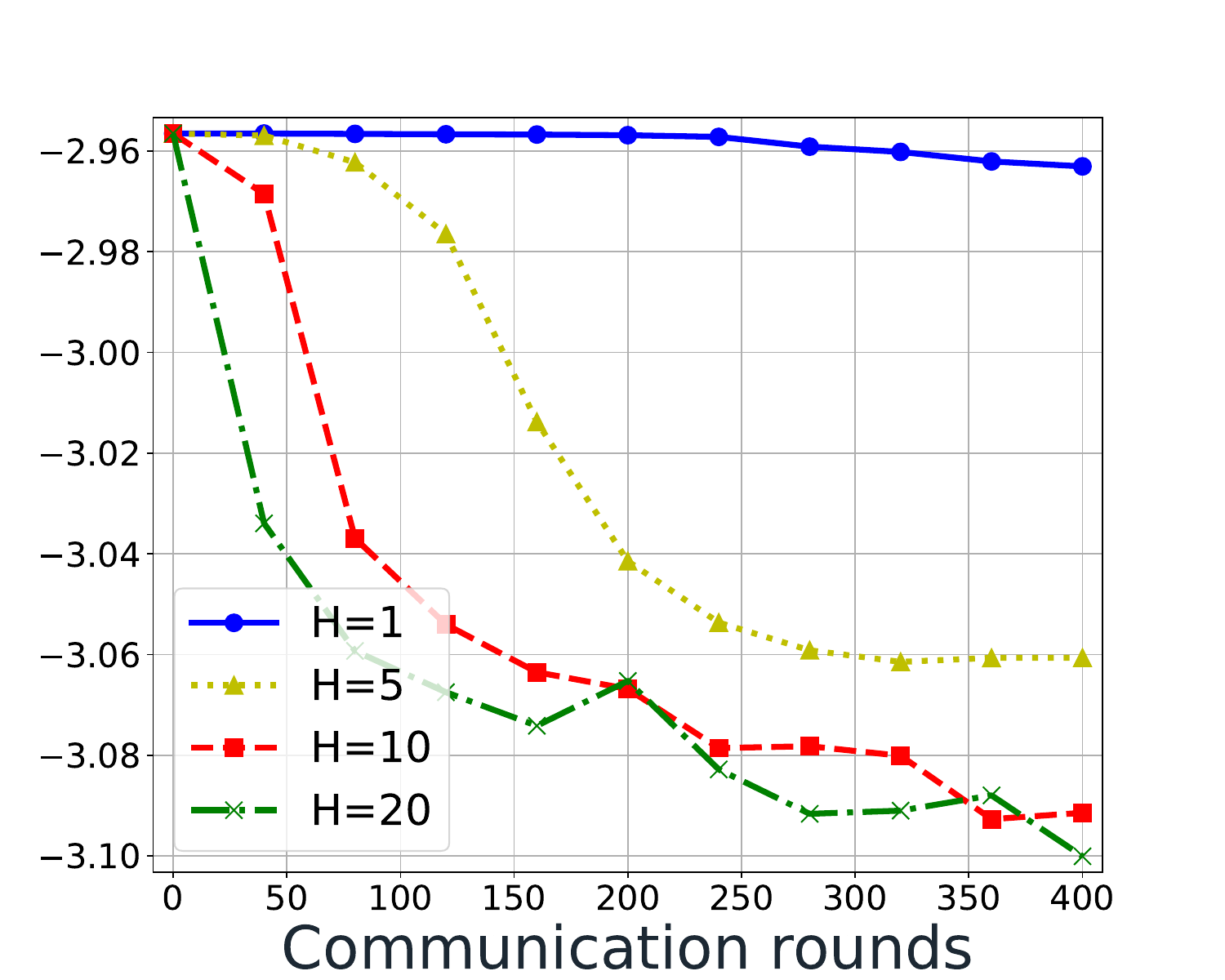}
\end{minipage}
&
\begin{minipage}{.3\textwidth}
\includegraphics[scale=.21, angle=0]{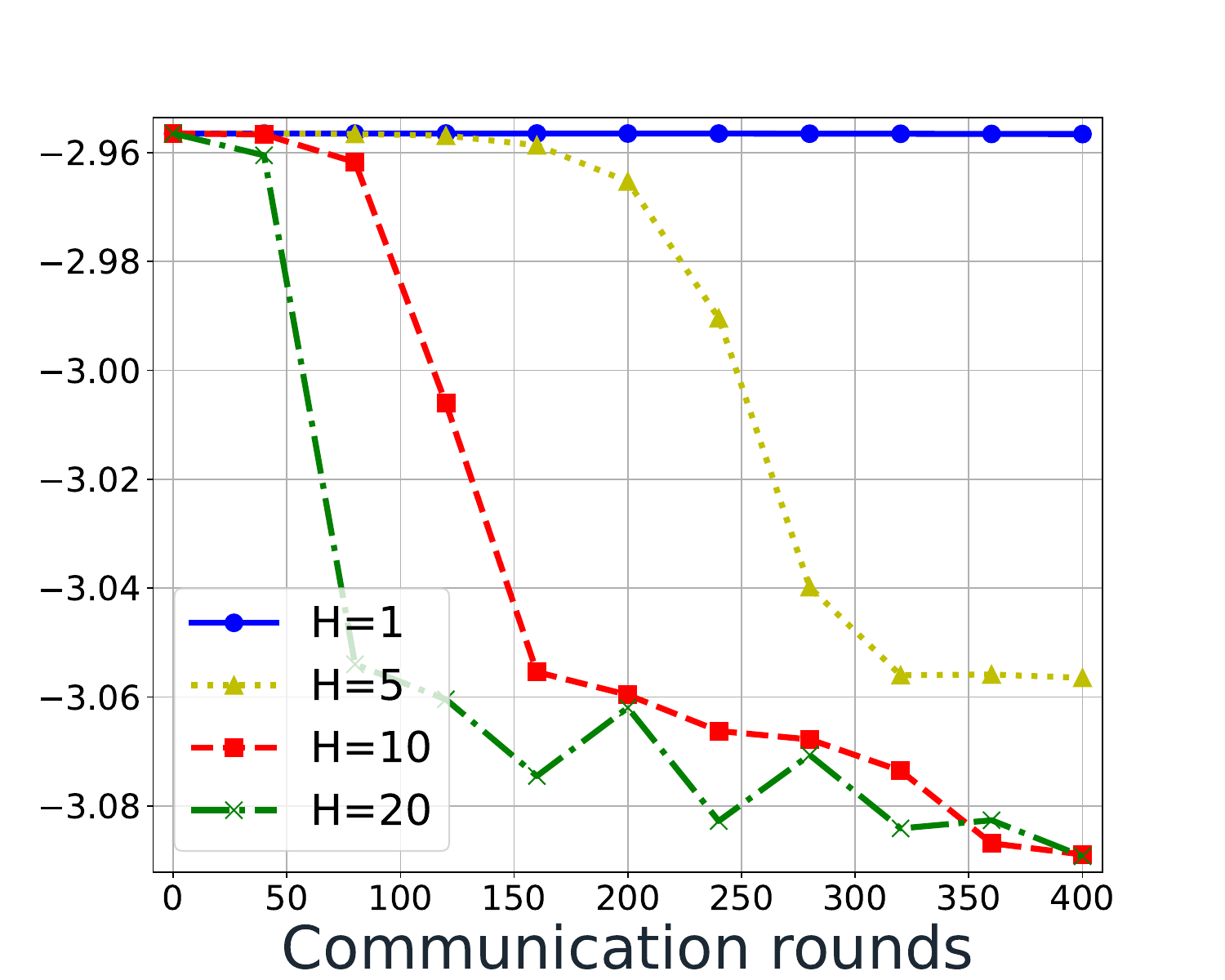}
\end{minipage}
&
\begin{minipage}{.3\textwidth}
\includegraphics[scale=.21, angle=0]{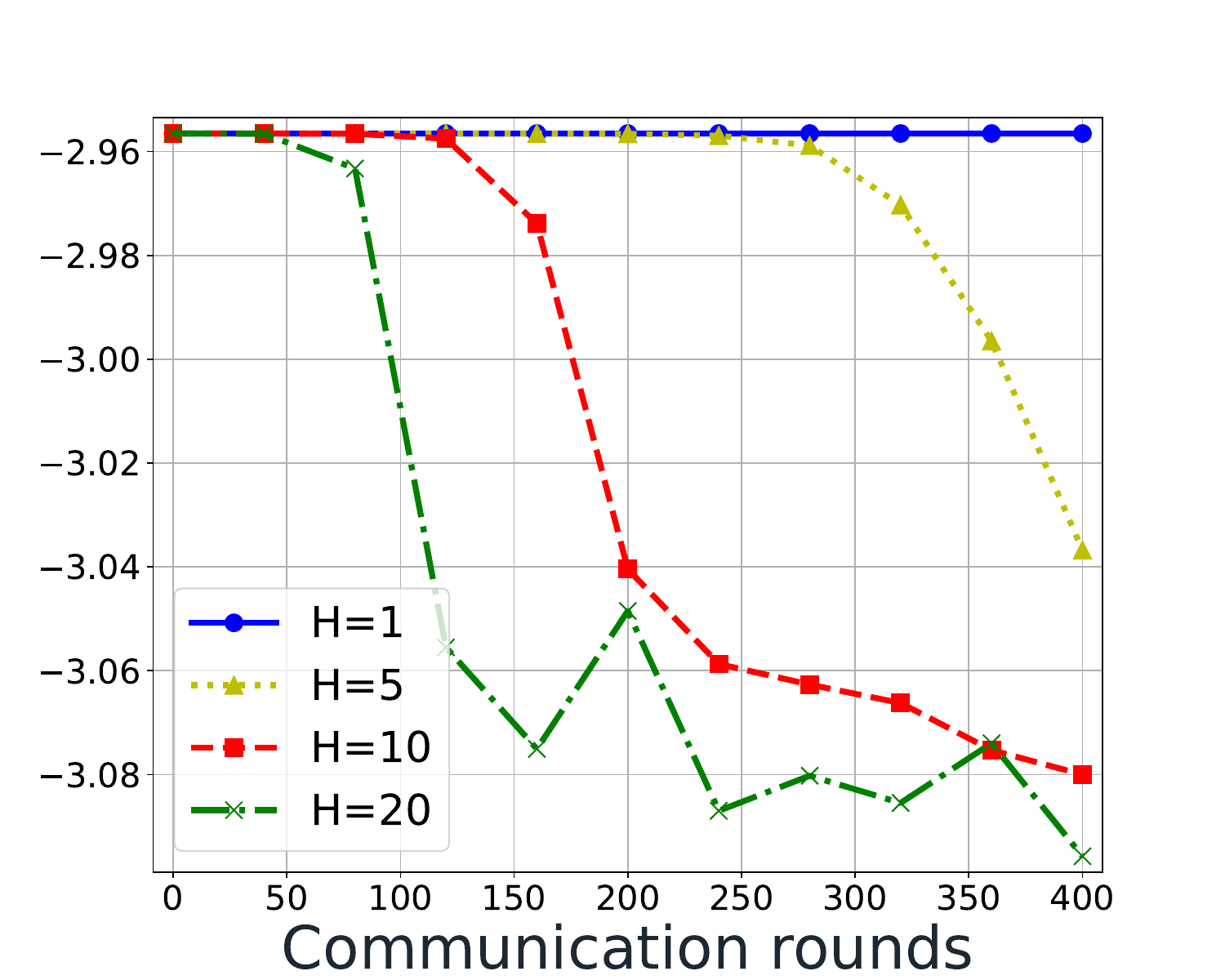}
\end{minipage}

\\
\hbox{}& & &\\
\hline\\
\rotatebox[origin=c]{90}{{\footnotesize {20 neurons}}}
&
\begin{minipage}{.3\textwidth}
\includegraphics[scale=.21, angle=0]{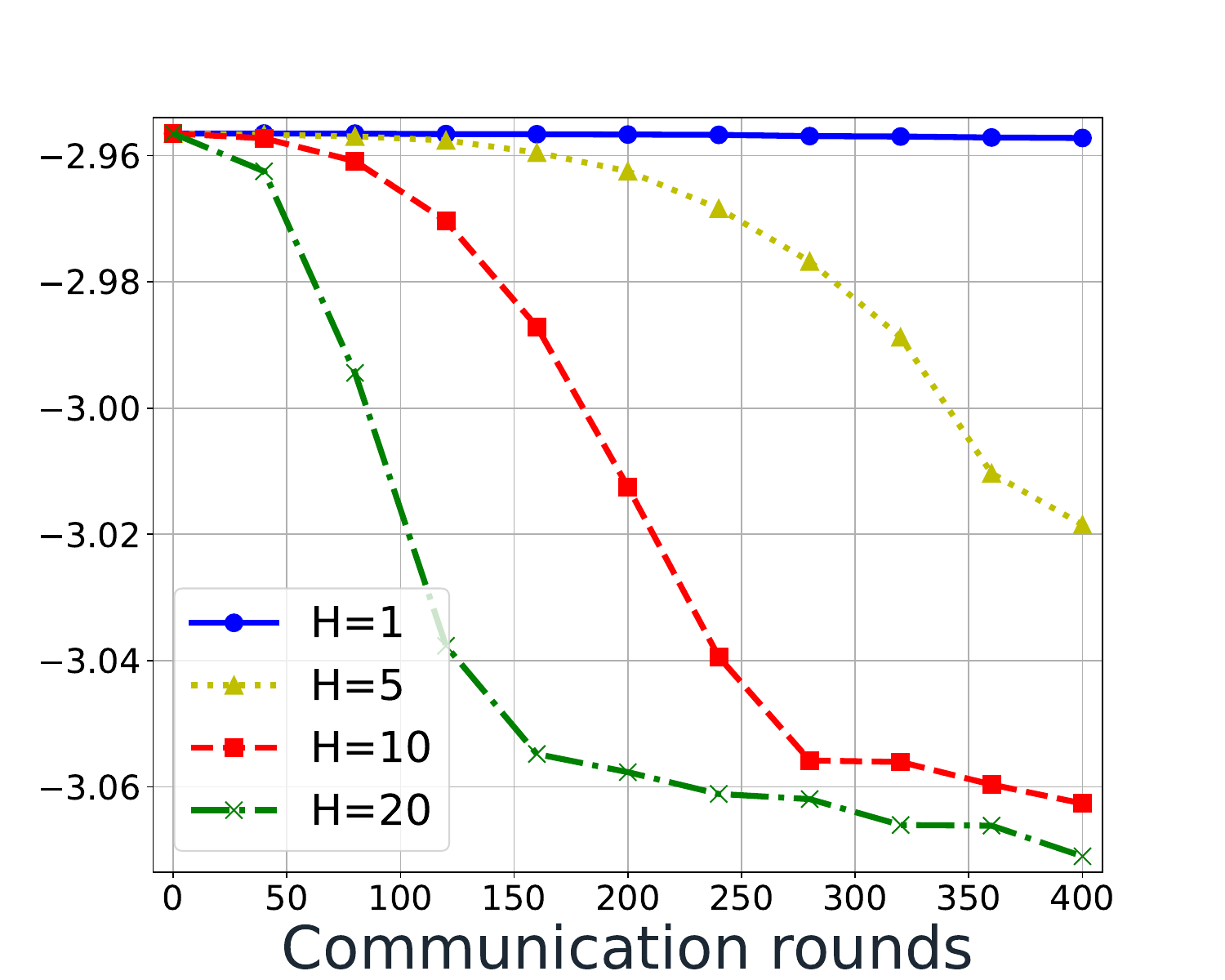}
\end{minipage}
&
\begin{minipage}{.3\textwidth}
\includegraphics[scale=.21, angle=0]{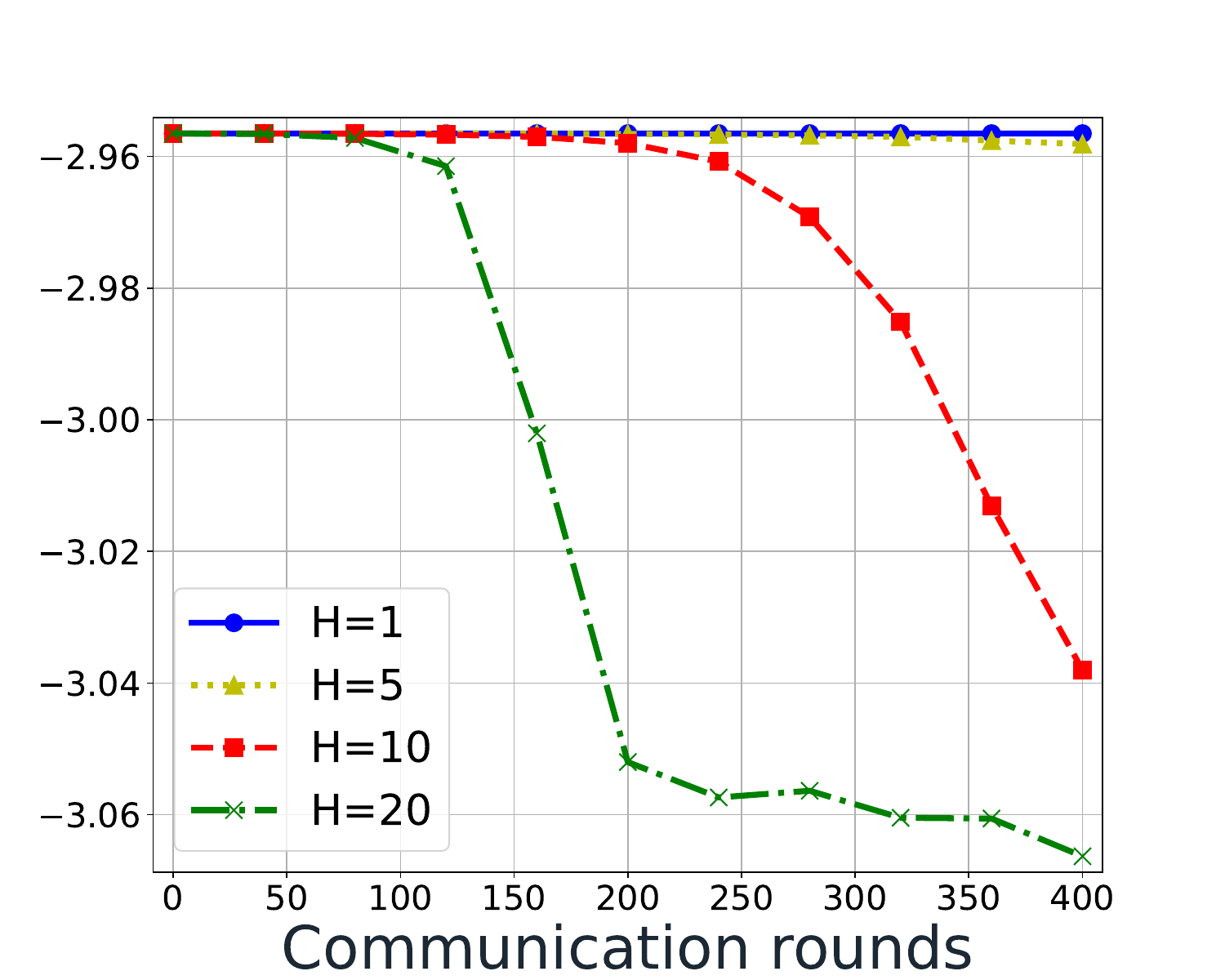}
\end{minipage}
&
\begin{minipage}{.3\textwidth}
\includegraphics[scale=.21, angle=0]{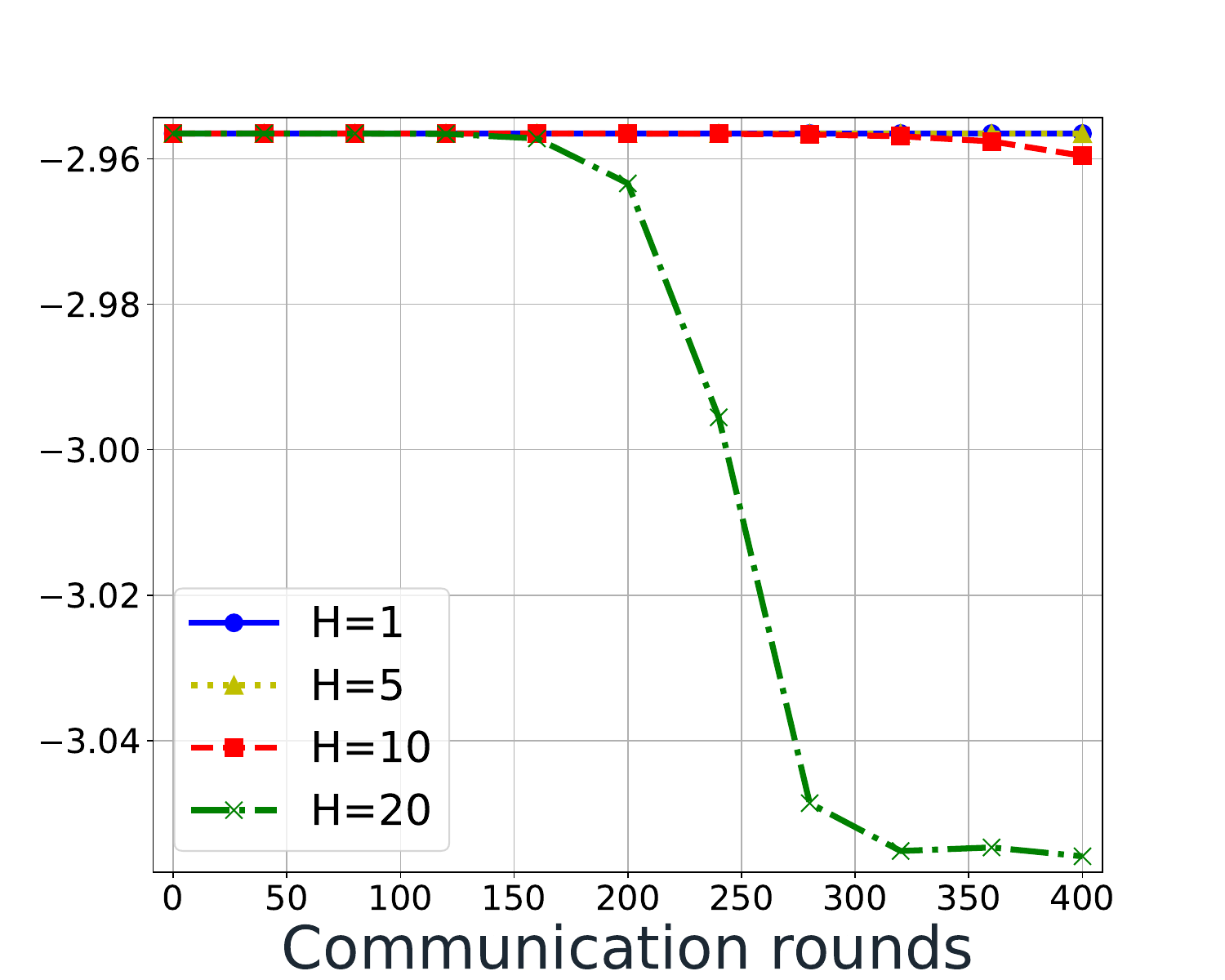}
\end{minipage}

\\
\hbox{}& & &\\
\hline\\
\rotatebox[origin=c]{90}{{\footnotesize {100 neurons}}}
&
\begin{minipage}{.3\textwidth}
\includegraphics[scale=.21, angle=0]{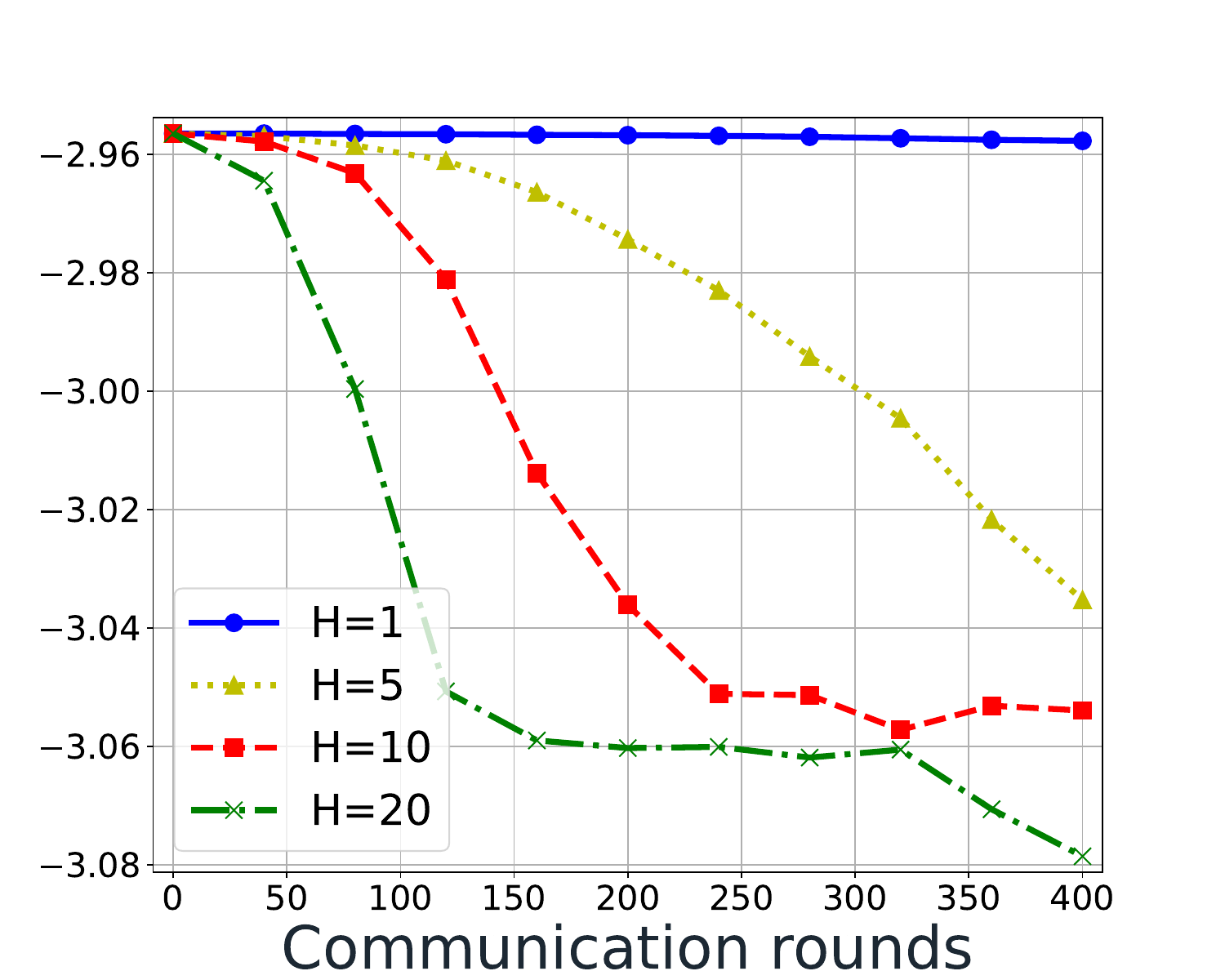}
\end{minipage}
&
\begin{minipage}{.3\textwidth}
\includegraphics[scale=.21, angle=0]{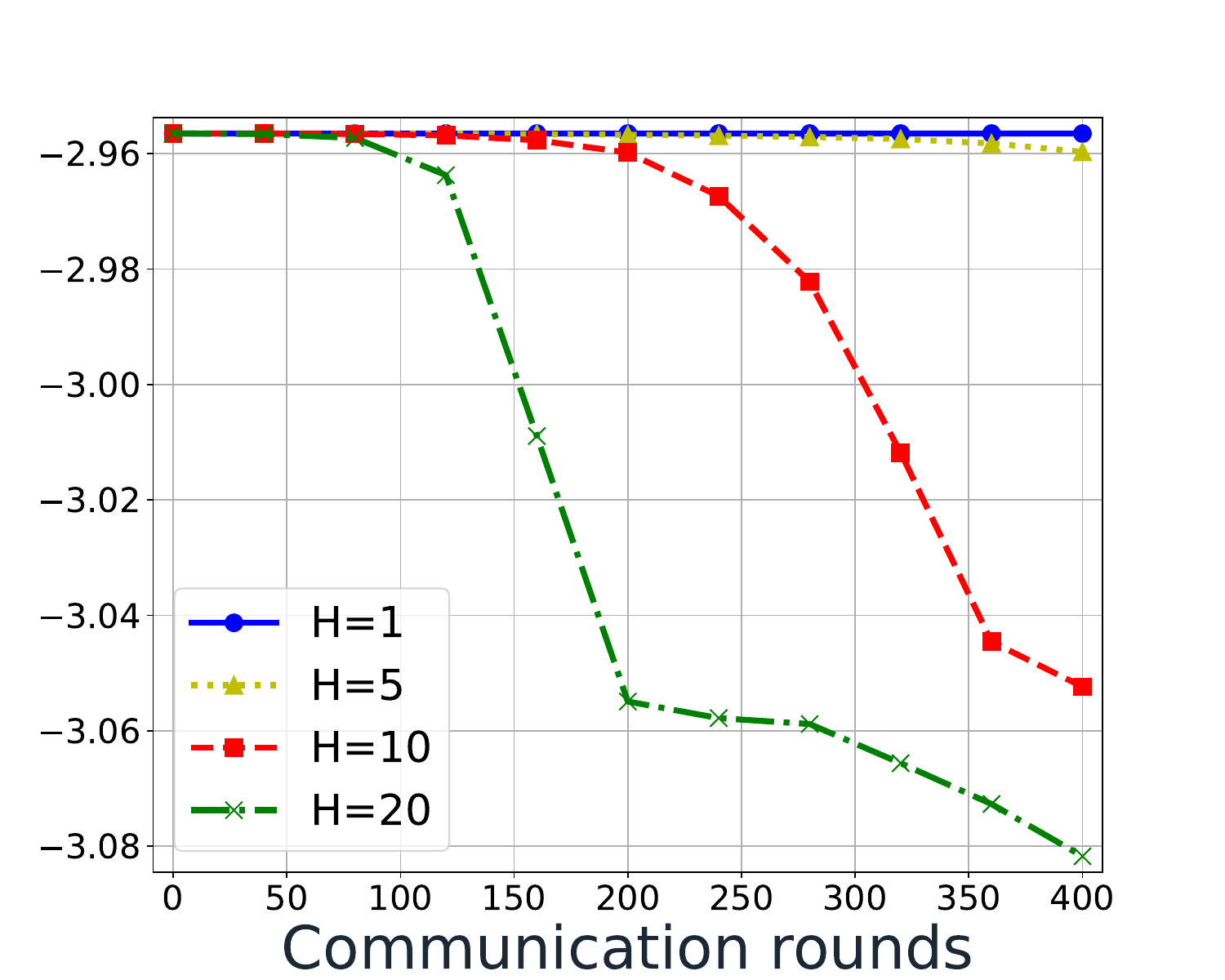}
\end{minipage}
&
\begin{minipage}{.3\textwidth}
\includegraphics[scale=.21, angle=0]{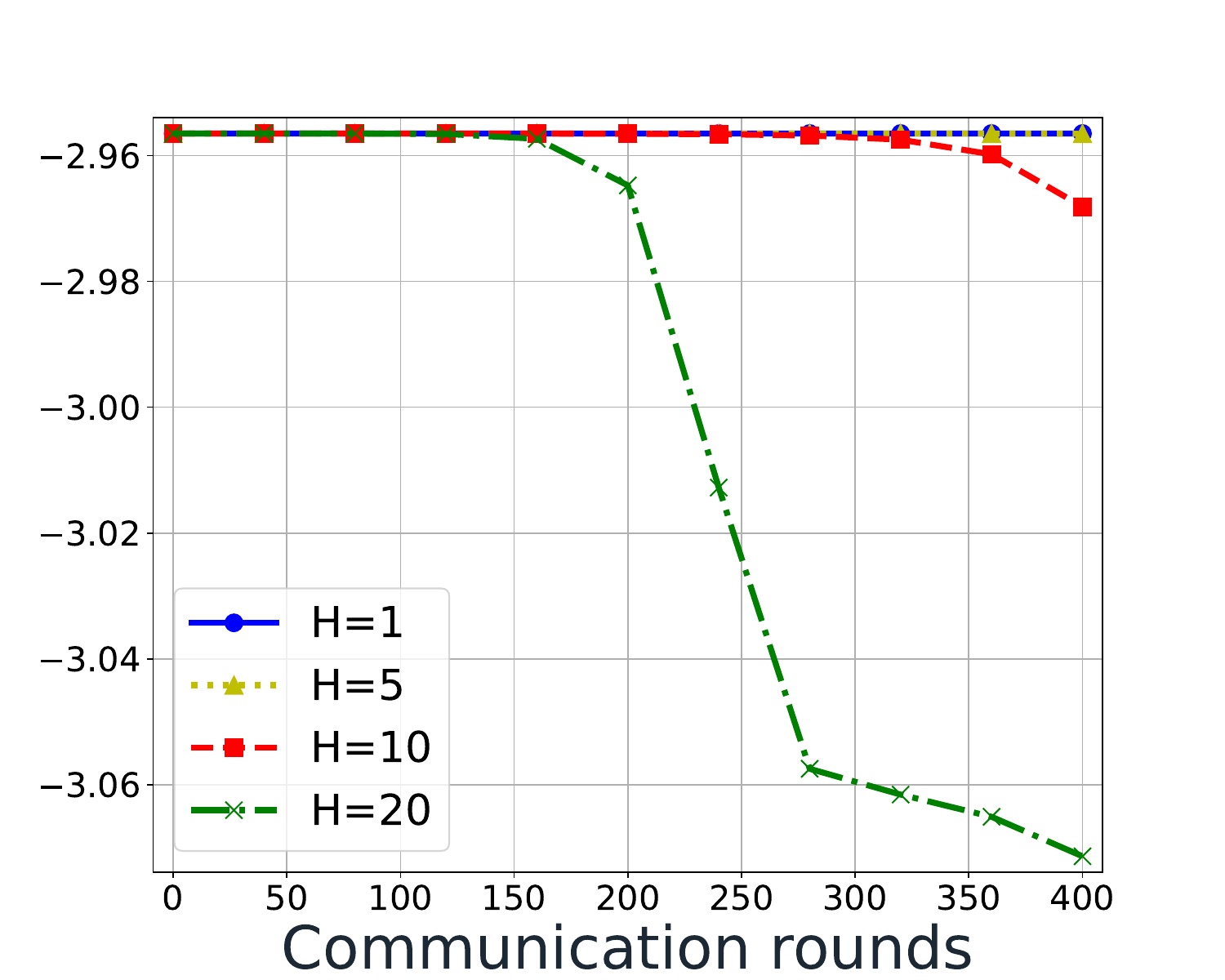}
\end{minipage}
\end{tabular}}
\vspace{.1in}
\captionof{figure}{{Performance of FedRZO$_{\texttt {nn}}$ on ReLU NN with different number of neurons in the layer and different value of the smoothing parameter $\eta$, using the CIFAR-10 dataset. The method performs better with smaller number of neurons and less communication frequency given a certain communication rounds. We also observe the smaller the smoothing parameter $\eta$, the longer it takes for the method to converge. This is aligned with our theory and analysis.}}
\label{fig:cifar10:sup}
\vspace{-.1in}
\end{table}

\subsection{Federated hyperparameter learning}\label{sec:hyper}
To validate FedRZO$_{\texttt {bl}}$, we consider {the following} FL hyperparameter learning problem for binary classification using logistic loss.
\begin{align*}
&\mbox{min}_{x \in X,\ y \in \mathbb{R}^{n} }\quad f(x,y) \triangleq  \frac{1}{m}\sum_{i=1}^m  {\sum_{\ell \in \yqro{\mathcal{S}}_i}\log\left( 1+\exp(-v_{i,\ell}U_{i,\ell}^{\fyy{\top}}y)\right)} \\
&\hbox{subject to }  \quad y \in \mbox{arg}\min_{y\in \mathbb{R}^{n}}\quad   h(x,y) \triangleq \frac{1}{m}\sum_{i=1}^m  \left({\sum_{\tilde\ell \in \yqro{\tilde S}_i}\log\left( 1+\exp(-v_{i,\tilde\ell}U_{i,\tilde\ell}^{\fyy{\top}}y)\right)}+x_i\frac{\|y\|^2}{2}\right) ,
\end{align*}
where $m$ is number of {clients}, $x$ denotes the regularization parameter for client $i$, \yqro{$\mathcal D_i \triangleq \{(U_{i,\ell},v_{i,\ell})\in \mathbb R^{n}\times \{-1,1\} \mid \ell\in \mathcal S_i\}$ denotes the local validation dataset of client $i$ and $\mathcal S_i=\{1,\dots,|\mathcal D_i|\}$ denotes the local index set, respectively. $\tilde{\mathcal D}_i \triangleq \{(U_{i,\tilde\ell},v_{i,\tilde\ell})\in \mathbb R^{n}\times \{-1,1\} \mid \tilde\ell\in \tilde{\mathcal S}_i\}$ denotes the local training dataset of client $i$ and $\tilde{\mathcal S}_i=\{1,\dots,|\tilde{\mathcal D}_i|\}$ denotes the local index set, respectively.}  The constraint set $X$ is considered as {$X := \{x \in \mathbb{R}^m \mid x \geq \munderbar{\mu}\mathbf{1}_m\}$, where $\munderbar{\mu}>0$. }
This problem is an instance of \eqref{prob:main_BL}, where the lower-level problem is $\ell_2$-regularized and the regularization parameter is a decision variable of the upper-level FL problem. The convergence results are presented in Fig.~\ref{fig5:FedRizo:cr:bilevel}.  \yqro{Let us define ${\tilde{f}_i}(x,y,\xi_i)=\log\left( 1+\exp(-v_{i,\ell}U_{i,\ell}^{\fyy{\top}}y)\right)$, where $\xi_i\in \mathcal D_i$ denotes a random data sample of client $i$. We note that ${\tilde{f}_i}(x,y,\xi_i)$ is $L_{0,y}^f(\xi_i)$-Lipschitz with respect to $y$ at any given $x$, where $L_{0,y}^f(\xi_i)=\|U_{i,\ell}\|$, which matches with our Assumption \ref{assum:bilevel} (i). The proof is in Section \ref{sec:appendix}.}

\subsection{Federated fair classification learning}\label{sec:fair} 
Here, we study the convergence of FedRZO$_{\texttt {bl}}$ in minimax FL. We
consider solving {the following} FL minimax formulation of the fair classification
problem~\cite{nouiehed2019solving} using ReLU neural networks. 
\begin{align*}
&\min_{x \in \mathbb{R}^n} \max_{y \in \mathbb{R}^c} \quad  \frac{1}{m}\sum_{i=1}^m    { \sum_{c=1}^C \sum_{\ell\in \yqro{\mathcal{S}}_{i,c}}  \left(v_{i,\ell}-\sum_{q=1}^{N_1}w_q \sigma(Z_{\bullet,q}U_{i,\ell})\right)^2- \frac{\lambda}{2}\|y\|^2} , 
\end{align*}
where $c$ denotes the class index and \yqro{$\mathcal D_{i,c} \triangleq \{(U_{i,\ell},v_{i,\ell})\in \mathbb R^{N_0}\times \{-1,1\} \mid \ell\in \mathcal S_{i,c}\}$ denotes the portion of client $i$'s local datset that is associated with class $c$ and $\mathcal S_{i,c}=\{1,\dots,|\mathcal D_{i,c}|\}$ denotes its index set, respectively.} The loss function follows the same formulation in Section~\ref{sec:NNs}, where an ReLU neural network is employed. 
This problem is nondifferentiable {and} nonconvex-strongly concave,  fitting well with the assumptions
in our work in addressing minimax FL problems. The performance of our algorithm
is presented in Figure~\ref{fig5:FedRizo:cr:minimax}.
\begin{table}[htbp]
\setlength{\tabcolsep}{0pt}
\centering{
\begin{tabular}{c  c  c }
\rotatebox[origin=c]{90}{{\footnotesize  }}
&
\begin{minipage}{.5\textwidth}
\centering
\includegraphics[scale=.3, angle=0]{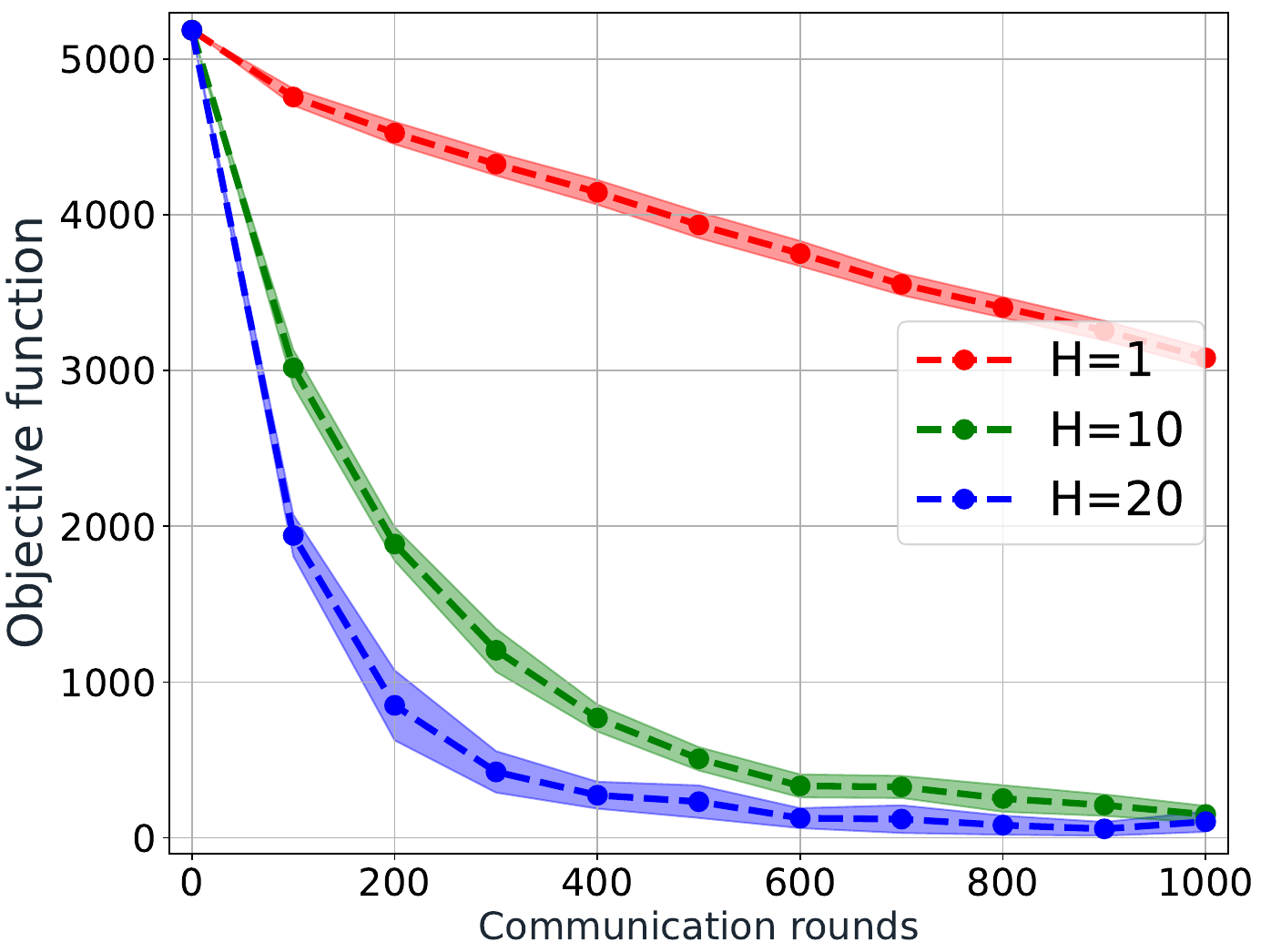} 
\end{minipage}
\end{tabular}}
\vspace{.1in}
\captionof{figure}{\yqro{Convergence of FedRZO$_{\texttt {bl}}$ in hyperparameter FL for $\ell_2$ regularized logistic loss, where we plot the loss function on test data for different values of local steps with $95\%$ CIs.}}
\label{fig5:FedRizo:cr:bilevel}
\end{table}

\begin{table}[htbp]
\setlength{\tabcolsep}{0pt}
\centering{
\begin{tabular}{c  c  c }
\rotatebox[origin=c]{90}{{\footnotesize  }}
&
\begin{minipage}{.5\textwidth}
\centering
\includegraphics[scale=.3, angle=0]{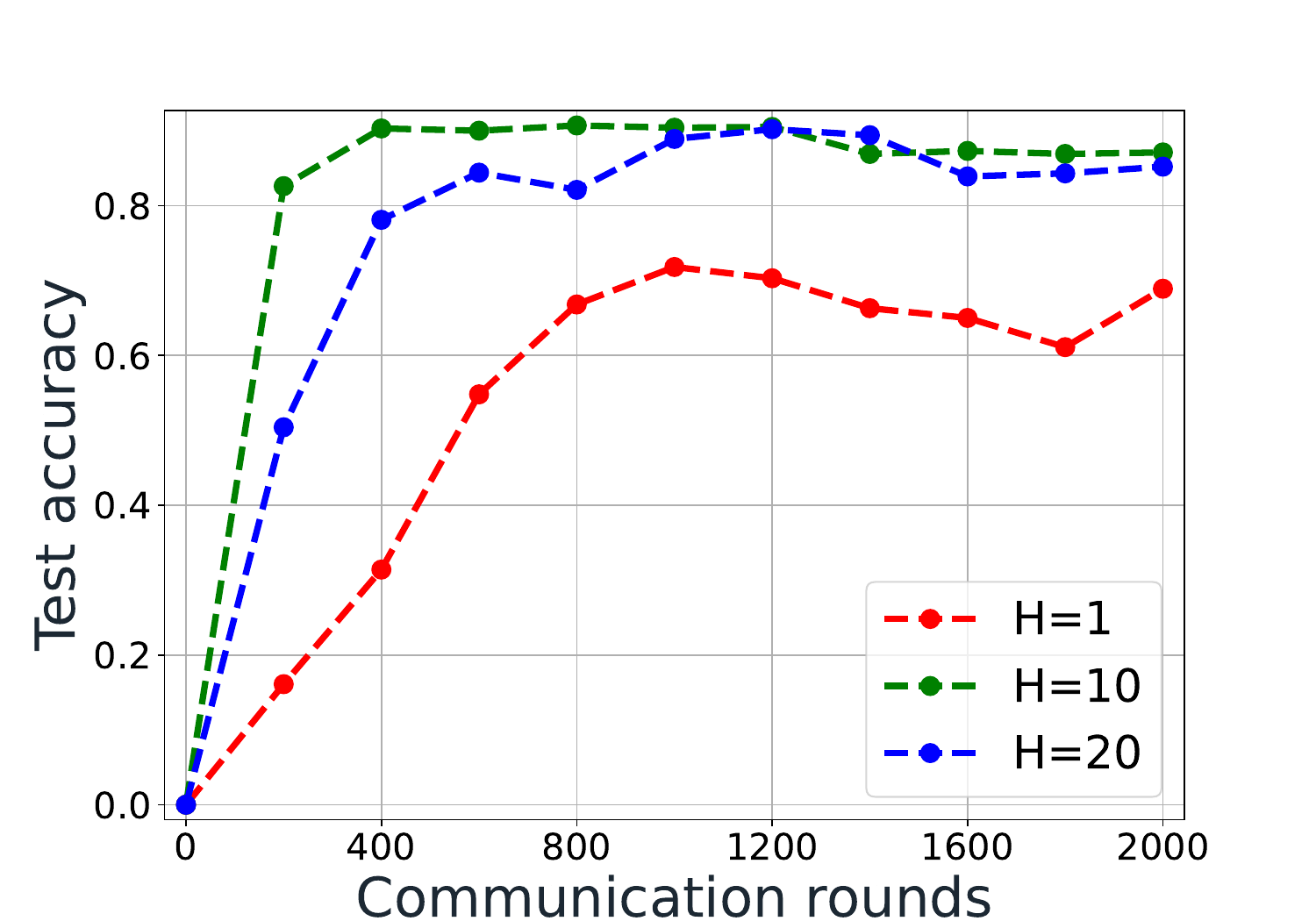} 
\end{minipage}
\end{tabular}}
\vspace{.1in}
\captionof{figure}{\yqro{Convergence of FedRZO$_{\texttt {bl}}$ in minimax FL, where we present test results in solving a nondifferentiable nonconvex-strongly concave FL minimax formulation of the fair classification problem~\cite{nouiehed2019solving}.}}
\label{fig5:FedRizo:cr:minimax}
\end{table}

\subsection{Federated two-stage stochastic single-leader multi-follower game}\label{sec:num2s}
Here, we consider a Stackelberg-Nash-Cournot equilibrium (SNCE) problem studied
in~\cite{sherali1983stackelberg}. The stochastic extension of this hierarchical
problem was considered in~\cite{cui2022complexity}. Our goal here lies in
implementing FedRZO$_{\texttt {2s}}$ for approximating the solution to the
smoothed implicit problem associated with the game, in a
communication-efficient and decentralized fashion.  The mathematical model  is {presented next}.  Consider a single-commodity market with a
firm acting as a single leader and $\tilde n$ number of firms acting as
followers.  Given the number of units generated by the leader (denoted
by $x \in \mathbb{R}$) and the realization of random state of the market
(characterized by $\xi$),  the follower firms compete through a Cournot game
mathematically captured by the {collection of the}  following optimization problems for
$j=1,\ldots,\tilde n$. 
\begin{align}\label{eqn:cournot_model}
&\max_{0\leq y_j  \leq \hat{y}_j}\quad  y_j p\left(y_j + x +\textstyle\sum_{\ell=1,  \ \ell \neq j}^{\tilde n} y_\ell(x,\xi), \xi\right) - c_j(y_j),
\end{align}
where $y_j$ denotes the firm $j$'s sales  assuming that $y_\ell(x,\xi)$ is fixed for all $\ell \neq j$, where $y_\ell(x,\xi)$ denotes the sales by other firms characterized by $x$ and $\xi$, $p(\bullet, \xi)$ denotes the random inverse demand curve at the aggregated sales by all firms (leader and followers), $c_j(\bullet)$ denotes the cost function of the $j$th firm,  and $\hat{y}_j$ is the capacity limit of the firm.  The leader's expectation-valued stochastic optimization problem is given as
\begin{align*}
&\max_{0\leq x \leq \hat{x}}\quad \mathbb{E}_{\xi}\left[xp\left( x +\textstyle\sum_{\ell=1}^{\tilde n} y_\ell(x,\xi), \xi\right)\right] - c_0(x),
\end{align*}
where $c_0(\bullet)$ and $\hat{x}$ are the production cost function and  capacity limit associated with the leading firm, respectively, and the tuple $(y_\ell(x,\xi))_{\ell=1}^{\tilde{n}}$ denotes the equilibrium to the lower-level Cournot competition.  

\noindent {\bf Model and algorithm setup.} In the experiments, we consider the
case that the demand curve is linear,  as $p(q ,\xi) := a(\xi) - bq$ where
$a(\xi)$ is uniformly distributed within the range $[7.5,12.5]$ and $b$ is a
constant taking values in $\{0.5, 1\}$. We let the cost functions be given as
$c_j(q) : = \frac{1}{2}{c_j} q^2$ for $j=0,\ldots,\tilde{n}$ where $c_j \in
[0.09,0.11]$ is a fixed constant.  Also, we set the capacity limits as
{$\hat{x} := 10$} and {$\hat{y}_j \in[2, 4]$}, all being constants.  We
consider $\eta \in \{1, 0.1, 0.01\}$ and choose the upper-level stepsize as
$\gamma:=10^{-3}$. In the Cournot setting, we consider different values
for the number of followers, as $\tilde n = \{10, 20 ,100, 1000\}$. As $\tilde
n$ increases, the condition number of the mapping, denoted by $\kappa_F$, would
also increase. In view of the condition $\tau\geq
\tfrac{-1}{\log(1-1/\kappa^2_F)}$, we tune $\tau$ as $\tau\in\{20, 35, 55,
80\}$, respectively, which results in the choices $\alpha\in
\{10^{-4},10^{-6},10^{-7},10^{-8}\}$ for the lower-level stepsize,
respectively. To highlight the benefit of the federated setting, we consider
$H=\{1, 10, 20\}$ and use a fixed total number of communication rounds as
$R=100$.

\noindent {\bf VI characterization of Cournot game.} Let us define the parametric mapping $G(x,\bullet,\xi):\mathbb{R}^{\tilde n} \to \mathbb{R}^{\tilde n} $ as $$G(x,y,\xi)\triangleq (c+b\mathbf{1}_{\tilde n})\odot y -a(\xi)\mathbf{1}_{\tilde n} +  \left(x+\textstyle\sum_{\ell=1}^{\tilde n} y_\ell\right)b \mathbf{1}_{\tilde n}, $$
where $\odot $ denotes the Hadamard product operator, $c \triangleq [c_1,\ldots,c_{\tilde n}]^{\fyy{\top}}$ and $\mathbf{1}_{\tilde n} \in \mathbb{R}^{\tilde n}$ denotes the vector of ones.   Under the prior setting, it can be shown that the mapping $G$ is a strongly monotone and Lipschitz continuous linear mapping with respect to $y$~\cite{demiguel2009stochastic}.  Consequently, a vector $y(x,\xi)\in \mathbb{R}^{\tilde n}$ is the unique Nash equilibrium to the Cournot game~\eqref{eqn:cournot_model} if and only if it is the unique solution to $\mbox{VI}\left(\prod_{j=1}^{\tilde n} [0,\hat{y}_j],G(x,\bullet, \xi)\right)$.  

\noindent {\bf Federated setting.} To compute an approximate solution to the SNCE problem, we consider a federated computational setting with $m$ clients, associated with iid random variables $\xi_i$ that admit the same probability distribution as that of $\xi$.  Let us define the stochastic local implicit loss functions as 
$${\tilde{f}_i}(x,y(x,\xi_i),\xi_i)\triangleq 0.5 \,c_0x^2 - x   \cdot (a(\xi_i) - b \cdot (x + \mathbf{1}_{\tilde n}^{\fyy{\top}} y(x,\xi_i)).$$ 
Then,  the SNCE problem can be compactly cast a federated two-stage SMPEC of the form
\begin{align*}
  \min_{x  \, \in \,  [0,\hat{x}]}\left\{  \frac{1}{m}\sum_{i=1}^m  \mathbb{E}_{\xi_i \in\mathcal{D}_i}  [\, {\tilde{f}_i}(x,y(x,\xi_i),\xi_i) \, ] \, \mid \, 
    y(x,\xi_i) \in \text{SOL}\left(\prod_{j=1}^{\tilde n} [0,\hat{y}_j], G(x, \bullet , \xi_i)\right), \ \hbox{for a.e. $\xi_i$}\right\}.
\end{align*}

\noindent {\bf Results and insights.} 
Table~\ref{table:smpec_numerics} provides the implementation results for FedRZO$_{\texttt {2s}}$ in addressing the SNCE problem. Notably, in the case when $H=1$, the communication occurs at every iteration, resulting in a standard distributed and parallel variant of the ZSOL scheme in~\cite{cui2022complexity}. Importantly, for each setting of $\tilde{n}$, as the number of local steps, i.e., $H$, increases, the global objective value appears to be improving. This is indeed interesting, as it supports the premise of using a federated framework  and highlights the advantage of employing a communication-efficient decentralized scheme in addressing two-stage SMPECs. We further observe that  FedRZO$_{\texttt {2s}}$ displays robustness with respect to the choice of the smoothing parameter $\eta$, as the global loss function only moderately changes as $\eta$ varies. To elaborate {further}, in Figure~\ref{fig:Fed2s:smpec} we display the progress of the algorithms with respect to the choices of $H$ in two of the settings in the table, where on the left we have $\tilde n = 100$, $\eta = 1$ and $b=0.5$, and on the right we have $\tilde n = 1000$, $\eta = 0.1$, and $b=1$. 

\begin{table}[htb]
\caption{Implementation results of FedRZO$_{\texttt {2s}}$ addressing the SNCE problem. FedRZO$_{\texttt {2s}}$ displays improvement as $H$ increases.}
\centering{\small{\yqro{
\begin{tabular}{clc|cclc|clcc|clcc}
\hline
\multicolumn{3}{c|}{\multirow{2}{*}{Setting}}             & \multicolumn{4}{c|}{$f(\bar x_K)$, $H=1$}                                       & \multicolumn{4}{c|}{$f(\bar x_K)$, $H=10$}                                      & \multicolumn{4}{c}{$f(\bar x_K)$, $H=20$}                                      \\ \cline{4-15} 
\multicolumn{3}{c|}{}                                     & \multicolumn{1}{c|}{$\eta =1$} & \multicolumn{2}{c|}{$\eta=0.1$} & $\eta =0.01$ & \multicolumn{2}{c|}{$\eta =1$} & \multicolumn{1}{c|}{$\eta=0.1$} & $\eta =0.01$ & \multicolumn{2}{c|}{$\eta=1$}  & \multicolumn{1}{c|}{$\eta=0.1$} & $\eta=0.01$ \\ \hline
\multicolumn{2}{c|}{$\tilde n=10,$}   & $b=0.5$           & \multicolumn{1}{c|}{$-7.810$}  & \multicolumn{2}{c|}{$-7.812$}   & $-7.813$     & \multicolumn{2}{c|}{$-35.211$} & \multicolumn{1}{c|}{$-35.211$}  & $-35.212$    & \multicolumn{2}{c|}{$-39.450$} & \multicolumn{1}{c|}{$-39.450$}  & $-39.450$   \\ \cline{3-3}
\multicolumn{2}{c|}{$\tau = 20$}      & $b=1$             & \multicolumn{1}{c|}{$-6.444$}  & \multicolumn{2}{c|}{$-6.444$}   & $-6.446$     & \multicolumn{2}{c|}{$-18.995$} & \multicolumn{1}{c|}{$-18.995$}  & $-18.999$    & \multicolumn{2}{c|}{$-19.408$} & \multicolumn{1}{c|}{$-19.409$}  & $-19.409$   \\ \hline
\multicolumn{2}{c|}{$\tilde n=20,$}   & $b=0.5$           & \multicolumn{1}{c|}{$-8.617$}  & \multicolumn{2}{c|}{$-8.616$}   & $-8.616$     & \multicolumn{2}{c|}{$-39.767$} & \multicolumn{1}{c|}{$-39.767$}  & $-39.767$    & \multicolumn{2}{c|}{$-44.158$} & \multicolumn{1}{c|}{$-44.159$}  & $-44.158$   \\ \cline{3-3}
\multicolumn{2}{c|}{$\tau=35$}        & $b=1$             & \multicolumn{1}{c|}{$-7.912$}  & \multicolumn{2}{c|}{$-7.913$}   & $-7.913$     & \multicolumn{2}{c|}{$-23.110$} & \multicolumn{1}{c|}{$-23.110$}  & $-23.110$    & \multicolumn{2}{c|}{$-23.349$} & \multicolumn{1}{c|}{$-23.349$}  & $-23.349$   \\ \hline
\multicolumn{2}{c|}{$\tilde n=100,$}  & $b=0.5$           & \multicolumn{1}{c|}{$-9.025$}  & \multicolumn{2}{c|}{$-9.023$}   & $-9.022$     & \multicolumn{2}{c|}{$-40.739$} & \multicolumn{1}{c|}{$-40.739$}  & $-40.739$    & \multicolumn{2}{c|}{$-45.174$} & \multicolumn{1}{c|}{$-45.174$}  & $-45.174$   \\ \cline{3-3}
\multicolumn{2}{c|}{$\tau=55$}        & $b=1$             & \multicolumn{1}{c|}{$-8.222$}  & \multicolumn{2}{c|}{$-8.222$}   & $-8.222$     & \multicolumn{2}{c|}{$-23.560$} & \multicolumn{1}{c|}{$-23.560$}  & $-23.560$    & \multicolumn{2}{c|}{$-23.901$} & \multicolumn{1}{c|}{$-23.901$}  & $-23.901$   \\ \hline
\multicolumn{2}{c|}{$\tilde n=1000,$} & $b=0.5$           & \multicolumn{1}{c|}{$-8.881$}  & \multicolumn{2}{c|}{$-8.887$}   & $-8.889$     & \multicolumn{2}{c|}{$-40.322$} & \multicolumn{1}{c|}{$-40.321$}  & $-40.321$    & \multicolumn{2}{c|}{$-44.800$} & \multicolumn{1}{c|}{$-44.800$}  & $-44.799$   \\ \cline{3-3}
\multicolumn{2}{c|}{$\tau=80$}        & $b=1$             & \multicolumn{1}{c|}{$-7.983$}  & \multicolumn{2}{c|}{$-7.981$}   & $-7.981$     & \multicolumn{2}{c|}{$-23.332$} & \multicolumn{1}{c|}{$-23.332$}  & $-23.332$    & \multicolumn{2}{c|}{$-23.681$} & \multicolumn{1}{c|}{$-23.681$}  & $-23.680$   \\ \hline
\end{tabular}}}}
\label{table:smpec_numerics}
\end{table}

\begin{table}[htbp]
\setlength{\tabcolsep}{0pt}
\centering{
\begin{tabular}{c  c  c }
\rotatebox[origin=c]{90}{{\footnotesize  }}
&
\begin{minipage}{.5\textwidth}
\centering
\includegraphics[scale=.28, angle=0]{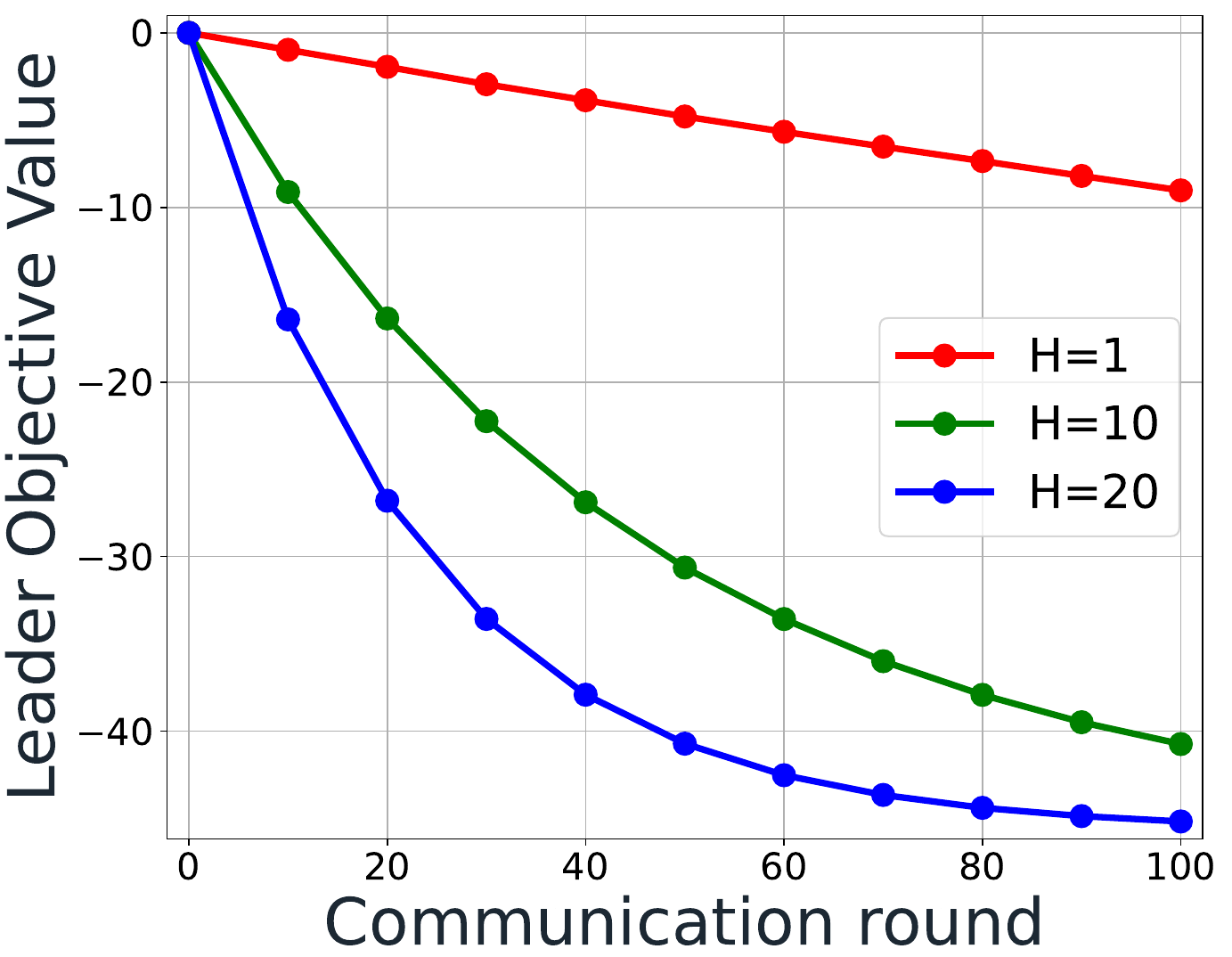} 
\end{minipage}
&
\begin{minipage}{.5\textwidth}
\centering
\includegraphics[scale=.28, angle=0]{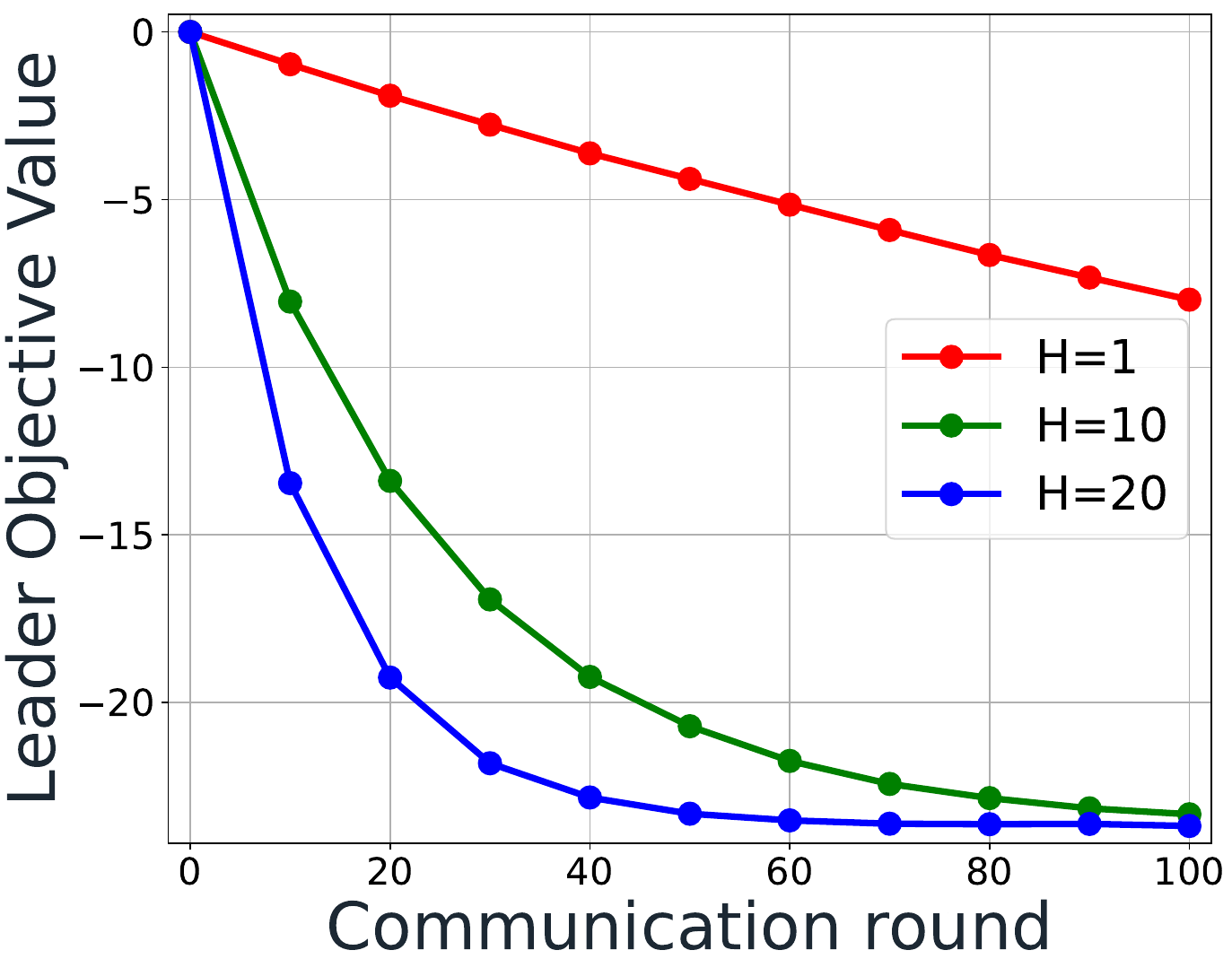} 
\end{minipage}
\end{tabular}}
\vspace{.1in}
\captionof{figure}{ (Left) Convergence of FedRZO$_{\texttt {2s}}$ for the SNCE problem with $\tilde n =100$ followers.  (Right) Convergence of FedRZO$_{\texttt {2s}}$ for the SNCE problem with $\tilde n =1000$ followers. In both settings, the performance of FedRZO$_{\texttt {2s}}$ improves as the number of local steps increases. This supports the need for communication-efficient decentralized optimization in addressing two-stage SMPECs.}
\label{fig:Fed2s:smpec}
\end{table}

\section{Concluding \fyy{r}emarks}\label{sec:conc} Federated learning has assumed growing relevance in machine learning. However, most practical problems are characterized by the presence of local objectives, jointly afflicted by
nonconvexity and nondifferentiability, precluding resolution by most FL schemes, which can cope with nonconvexity in only smooth settings. We resolve this gap via a zeroth-order communication-efficient FL framework that can
contend with both nondifferentiability and nonsmoothness with rate and complexity guarantees for computing approximate Clarke-stationary points. Extensions to nonconvex bilevel problems with a strongly convex minimization in
the lower level, nonconvex-strongly concave minimax  problems, and two-stage mathematical programs with equilibrium constraints are developed via inexact generalizations. In each case, \fyy{assuming that clients are associated with non-iid data,} we provide explicit iteration and communication
complexity bounds. Our results are novel in that this appears to be the first time that these classes of nonsmooth nonconvex and hierarchical problems are provably addressed in federated learning. An interesting, yet challenging
research question lies in investigating whether we can address these hierarchical problems in the absence of uniqueness of the lower-level solution. This is left as a direction to our future research. 


\bibliographystyle{plain}
\bibliography{references}

\section{Proofs of the main results}\label{sec:proofs}

\subsection{Preliminaries for the proof of Proposition~\ref{thm:main_bound}}
 \fy{In this subsection, we provide preliminary results to prove Theorem~\ref{thm:fed2smpec}.  We begin with introducing some terms in the following.}
\begin{definition}\label{def:basic_terms}                                                                                                         
Let {the} function $d_i(x)$ be defined as $d^{\eta}_i(x)\triangleq \frac{1}{2\eta}\mbox{dist}^2(x,X_i)$. Let $x_{i,k}$ be given by Algorithm~\ref{alg:PZO_LSGD} for all $i \in [m]$ and $k\geq 0$. Let us define {the following terms:}  
\begin{align*}
& g^{\eta}_{i,k} \triangleq \frac{n}{\eta^2} \left({\tilde{f}_i}(x_{i,k}+v_{i,k},\xi_{i,k}) - {\tilde{f}_i}(x_{i,k}\yqro{-v_{i,k}},\xi_{i,k})\right)v_{i,k},\quad \nabla_{i,k}^{\eta,d} \triangleq \frac{1}{\eta}(x_{i,k}-\mathcal{P}_{X_i}(x_{i,k})),\\
& \bar{g}^{\eta}_k \triangleq  \frac{\sum _{i=1}^m g^{\eta}_{i,k}}{m}, \quad \bar{\nabla}_{k}^{\eta,d} \triangleq \frac{\sum _{i=1}^m {\nabla}_{i,k}^{\eta,d}}{m}, \quad \bar{x}_k   \triangleq \frac{\sum_{i=1}^m x_{i,k}}{m}, \quad  \bar{e}_k \triangleq  \frac{\sum_{i=1}^m \|x_{i,k}-\bar{x}_k\|^2}{m}.
\end{align*}
Here, $\bar{x}_k$ is an auxiliary sequence {that} denotes the average iterates of the {clients} at any iteration $k$ and $\bar{e}_k$ denotes an average consensus error at that iteration.
\end{definition}
We will make use of the following notation for the history of FedRZO$_{\texttt {nn}}$. Let $ \mathcal{F}_0=\{\fy{\hat{x}_0}\}$ and
\begin{align*}
\mathcal{F}_k \triangleq \cup_{i=1}^{m} \cup_{t=0}^{k-1} \{ \xi_{i, t}, v_{i, t} \}\cup\{\fy{\hat{x}_0}\}, \quad \hbox{for all } k\geq 1.
\end{align*}
Throughout the proofs in this subsection, we assume that {the $i$th client} {generates} iid random samples $\xi_i\in \mathcal{D}_i$ and $v_i\in \eta\mathbb{S}$. These samples are assumed to be independent across {clients}.

\begin{remark}[Compact local representation of Algorithm~\ref{alg:PZO_LSGD}] 
Let us define $\mathcal{I}\triangleq \{K_1,K_2,\ldots\}$ where $K_r\triangleq T_r-1$ for $r\geq 1$. The following equation, for $k\geq 0$, compactly represents the update rules of Algorithm~\ref{alg:PZO_LSGD}.
\begin{align}\label{eqn:compact_update_rules}
x_{i,k+1} := 
\begin{cases}
\frac{1}{m}\sum_{j=1}^m\left(x_{j,k}-\gamma  \left(g^{\eta}_{j,k}+\nabla^{\eta,d}_{j,k}\right)\right), & k \in \mathcal{I}\\
 x_{i,k} -\gamma \left(  g^{\eta}_{i,k}+\nabla^{\eta,d}_{i,k}\right) ,  & k \notin \mathcal{I}.							
\end{cases} 
\end{align}
\end{remark}

\yqro{{\bf Proof of Lemma \ref{SphericalSmooth}.}

\proof{Proof.} (i) From Lemma 1 in \cite{cui2022complexity}, we have $\nabla h^\eta(x) = \tfrac{n}{\eta}\mathbb{E}_{v \in \mathbb{\eta S}}  \left[  h(x+ v)  \tfrac{v}{\|v\|}  \right]$ for any $x\in \mathbb{R}^n$. By the symmetric property of the distribution of $v$, we have $\mathbb{E}_{v \in \mathbb{\eta S}}  \left[  h(x+ v)  \tfrac{v}{\|v\|}  \right]=\mathbb{E}_{v \in \mathbb{\eta S}}  \left[  h(x- v)  \tfrac{-v}{\|v\|}  \right]$. Therefore, we have $\tfrac{n}{2\eta}\mathbb{E}_{v}  \left[  h(x+ v)-h(x-v) \tfrac{v}{\|v\|}  \right]=\tfrac{n}{2\eta}\mathbb{E}_{v}  \left[  h(x+ v) \tfrac{v}{\|v\|}  \right]+\tfrac{n}{2\eta}\mathbb{E}_{v}  \left[h(x-v) \tfrac{-v}{\|v\|}  \right]=\tfrac{n}{\eta}\mathbb{E}_{v}  \left[  h(x+ v)  \tfrac{v}{\|v\|}  \right] = \nabla h^\eta(x)$.

(ii, iii) See Lemma 1 in \cite{cui2022complexity}.

(iv) From  \cite[Lemma 8]{YNS_Automatica12}, we have $\| \nabla h^{\eta}(x) - \nabla h^{\eta} (y)\| \leq \frac{C_nL_0}{\eta}  \|x-y \|$, where $C_n=\tfrac{2}{\pi}\tfrac{n!!}{(n-1)!!}$ if $n$ is even, and $C_n=\tfrac{n!!}{(n-1)!!}$ if $n$ is odd. It remains to show that $C_n\leq \sqrt n$, \fyy{for any $n\geq 1$}.

Let $w_n=\int_0^{\pi/2}\text{sin}^n(x)dx$ \fyy{denote} the Wallis's integral. \fyy{We} have two properties: (i) $w_n=\tfrac{n-1}{n}w_{n-2}$ for $n\geq 2$, where we define $w_0=\tfrac{\pi}{2}$ and $w_1=1$; (ii) $w_{n+1}<w_n$ \cite{spivak2006calculus}. Then, by \fyy{using} induction, we have $w_n = \tfrac{\pi}{2}\tfrac{(n-1)!!}{n!!}$ for all even $n$, and $w_n = \tfrac{(n-1)!!}{n!!}$ for all odd $n$. Therefore, we have $C_n=\tfrac{1}{w_n}$.

Next, by noting that $(n+1)w_nw_{n+1}=(n+1)\tfrac{\pi}{2}\tfrac{(n-1)!!}{n!!}\tfrac{n!!}{(n+1)!!}=\tfrac{\pi}{2}$, and by property (ii), we have $w_n > \sqrt{\tfrac{\pi}{2(n+1)}} = \sqrt{\tfrac{\pi n}{2(n+1)}}\tfrac{1}{\sqrt n} \geq \sqrt{\tfrac{\pi }{3}}\tfrac{1}{\sqrt n}$, for \fyy{any} $n\geq 2$. Therefore, we obtain $C_n=\tfrac{1}{w_n}<\sqrt{\tfrac{3}{\pi}}\sqrt n < \sqrt n$.

}

\begin{lemma}\label{lem:bar_x_recursion}                                                                                                          
Consider Algorithm~\ref{alg:PZO_LSGD} and Definition~\ref{def:basic_terms}. For all $k\geq 0$, we have $\bar x_{k+1} = \bar x_k - \gamma \left(\bar{g}^{\eta}_k+\bar{\nabla}_k^{\eta,d}\right).$
\end{lemma}
 \proof{Proof.}
\noindent {\it Case 1:} When $k \in \mathcal{I}$, from equation~\eqref{eqn:compact_update_rules} we can write
\begin{align*}
 x_{i,k+1}& =\tfrac{1}{m}\textstyle\sum_{j=1}^m\left(x_{j,k}-\gamma \left(g^{\eta}_{j,k}+\nabla^{\eta,d}_{j,k}\right)\right) \\ 
& =  \tfrac{1}{m}\textstyle\sum_{j=1}^m x_{j,k}-\gamma\frac{1}{m}\sum_{j=1}^m(g^{\eta}_{j,k}+\nabla^{\eta,d}_{j,k})  = \bar x_k - \gamma \left(\bar{g}^{\eta}_k+\bar{\nabla}_k^{\eta,d}\right),
\end{align*}
where the last equation is implied by the definition of $\bar x_k$, $\bar{g}_k^{\eta}$ and $\bar{\nabla}_k^{\eta,d}$. {Averaging} on both sides over $i=1,\ldots,m$, we obtain $\bar x_{k+1} = \bar x_k - \gamma \left(\bar{g}^{\eta}_k+\bar{\nabla}_k^{\eta,d}\right)$.

\noindent {\it Case 2:} When $k \notin \mathcal{I}$, from equation~\eqref{eqn:compact_update_rules} we can write $x_{i,k+1}=x_{i,k} - \gamma g^{\eta}_{i,k}+\nabla^{\eta,d}_{i,k}.$ Summing over $i=1,\ldots,m$ on both sides and then dividing both sides by $m$, we obtain
$$\tfrac{1}{m}\textstyle\sum_{i=1}^m x_{i,k+1}=\tfrac{1}{m}\textstyle\sum_{i=1}^m x_{i,k} - \gamma \tfrac{1}{m}\textstyle\sum_{i=1}^m (g^{\eta}_{i,k}+\nabla^{\eta,d}_{i,k}).$$
Invoking Definition~\ref{def:basic_terms}, we obtain the desired result.

\yqro{
\begin{lemma}[L\'evy concentration on sphere $\mathbb S$ {\cite[Proposition 3.11 and Example 3.12]{wainwright2019high}}]\label{lem:Levy}\em
Let $h:\mathbb{R}^n\to \mathbb{R}$ be a given $L_0$-Lipschitz continuous function \fyy{and $v$ be} uniformly distributed on sphere $\mathbb S$. Then,
\begin{align}
\mathbb P(|h(v)-\mathbb E[h(v)]|\geq \epsilon) \leq 2\sqrt{2\pi} e^{-\tfrac{n\epsilon^2}{8L_0^2}}, \text{ for all } \epsilon>0. 
\end{align}
\end{lemma}}

\begin{lemma}[Properties of zeroth-order stochastic local gradient]\label{hetero:local:gradient}\em                              
Let Assumption~\ref{assum:Heterogeneous} hold. Consider Algorithm~\ref{alg:PZO_LSGD}. Then, the following relations hold for all $k\geq 0$ and all $i \in [m]$ {in an \fy{almost-sure} sense}.

\noindent (i) $\mathbb{E}\left[g^{\eta}_{i,k}+ \nabla^{\eta,d}_{i,k}\mid \mathcal{F}_k\right] = \nabla {\bf f}_i^{\eta}(x_{i, k})$.

\noindent (ii) $\mathbb{E}\left[\|g^{\eta}_{i,k}+ \nabla^{\eta,d}_{i,k}\|^2 \mid \mathcal{F}_k\right] \leq \yqro{32\sqrt{2\pi}L_0^2n} +2\|\nabla^{\eta,d}_{i,k}\|^2$.

\noindent (iii) $\mathbb{E}\left[\|g^{\eta}_{i,k}+ \nabla^{\eta,d}_{i,k} -  \nabla {\bf f}_i^{\eta}(x_{i, k})\|^2 \mid \mathcal{F}_k\right] \leq \yqro{16\sqrt{2\pi}L_0^2n}$. 
\end{lemma}
\proof{Proof.}
\noindent (i) From the definition of the zeroth-order stochastic gradient, we can write
\begin{align*}
\mathbb{E}\left[g^{\eta}_{i,k} \mid \mathcal{F}_k\right]  &= \mathbb{E}_{v_{i, k}}\left[\mathbb{E}_{\xi_{i, k}}\left[\tfrac{n}{\yqro{2}\eta^2} \left({\tilde{f}_i}(x_{i,k}+v_{i,k},\xi_{i,k}) - {\tilde{f}_i}(x_{i,k}\yqro{-v_{i,k}},\xi_{i,k})\right)v_{i,k} \mid \mathcal{F}_k\cup\{v_{i,k}\}\right]\right]\\
& =\mathbb{E}_{v_{i, k}}\left[\tfrac{n}{\yqro{2}\eta^2}  f_i(x_{i,k}+v_{i,k})   v_{i,k} \mid \mathcal{F}_k\right] \yqro{+\mathbb{E}_{v_{i, k}}\left[\tfrac{n}{2\eta^2}  f_i(x_{i,k}-v_{i,k})(-v_{i,k}) \mid \mathcal{F}_k\right]}\\
&\yqro{=\fyy{\left(\tfrac{n}{\eta^2}\right)}\mathbb{E}_{v_{i, k}}\left[  f_i(x_{i,k}+v_{i,k})   v_{i,k} \mid \mathcal{F}_k\right]}\\
&\stackrel{\tiny \mbox{Lemma~\ref{SphericalSmooth} (i)}}{=} \nabla {{ f}_{i}^\eta(x_{i, k})}.
\end{align*}
Adding the preceding equality {to} $\nabla^{\eta,d}_{i,k}$, noting that $\nabla^{\eta,d}_{i,k}=\nabla d_i^{\eta}(x_{i,k})$, and using equation~\eqref{prob:main_unconstrained}, we obtain the relation in (i).

\noindent (ii) From the definition of $g^{\eta}_{i,k}$ and that $\|v_{i,k}\| =\eta$, we have
\begin{align*}
\mathbb{E}\left[\left\|g^{\eta}_{i,k}\right\|^2 \mid  \mathcal{F}_k\cup\{\xi_{i,k}\} \right]  &=\left(\tfrac{n^2}{\yqro{4}\eta^4}\right) \mathbb{E}\left[ \large|  {\tilde{f}_i}(x_{i,k}+v_{i,k},\xi_{i,k}) - {\tilde{f}_i}(x_{i,k}\yqro{-v_{i,k}},\xi_{i,k})  \large|^2 \mid \mathcal{F}_k\cup\{\xi_{i,k}\} \right]\|v_{i,k}\|^2 \\
& =\left(\tfrac{n^2}{\yqro{4}\eta^2}\right) \mathbb{E}\left[ \large|{\tilde{f}_i}(x_{i,k}+v_{i,k},\xi_{i,k}) - {\tilde{f}_i}(x_{i,k}\yqro{-v_{i,k}},\xi_{i,k}) \large|^2 \mid \mathcal{F}_k\cup\{\xi_{i,k}\} \right].
\end{align*}
\yqro{Adding and subtracting $\mathbb E_{\hat v}[{\tilde{f}_i}(x_{i,k}+\hat v,\xi_{i,k})]$, where $\hat v$ is uniformly distributed on $\eta \mathbb S$, and using the fact that $(a-b)^2\leq 2a^2+2b^2$ \fyy{for any $a,b \in \mathbb{R}$}, we obtain
\begin{align*}
\mathbb{E}\left[\left\|g^{\eta}_{i,k}\right\|^2 \mid  \mathcal{F}_k\cup\{\xi_{i,k}\} \right] &\leq \tfrac{n^2}{2\eta^2} \mathbb{E}\left[ \large|{\tilde{f}_i}(x_{i,k}+v_{i,k},\xi_{i,k}) - \mathbb E_{\hat v}[{\tilde{f}_i}(x_{i,k}+\hat v,\xi_{i,k})] \large|^2 \mid \mathcal{F}_k\cup\{\xi_{i,k}\} \right]\\
&+\tfrac{n^2}{2\eta^2} \mathbb{E}\left[ \large|{\tilde{f}_i}(x_{i,k}-v_{i,k},\xi_{i,k}) - \mathbb E_{\hat v}[{\tilde{f}_i}(x_{i,k}+\hat v,\xi_{i,k})] \large|^2 \mid \mathcal{F}_k\cup\{\xi_{i,k}\} \right].
\end{align*}
Next, by the symmetric property of the uniform distribution of $v_{i,k}$, we have
\begin{align*}
&\mathbb{E}\left[ \large|{\tilde{f}_i}(x_{i,k}+v_{i,k},\xi_{i,k}) - \mathbb E_{\hat v}[{\tilde{f}_i}(x_{i,k}+\hat v,\xi_{i,k})] \large|^2 \mid \mathcal{F}_k\cup\{\xi_{i,k}\} \right]\\ 
&= \mathbb{E}\left[ \large|{\tilde{f}_i}(x_{i,k}-v_{i,k},\xi_{i,k}) - \mathbb E_{\hat v}[{\tilde{f}_i}(x_{i,k}+\hat v,\xi_{i,k})] \large|^2 \mid \mathcal{F}_k\cup\{\xi_{i,k}\} \right].
\end{align*}
Then we obtain
\begin{align*}
\mathbb{E}\left[\left\|g^{\eta}_{i,k}\right\|^2 \mid  \mathcal{F}_k\cup\{\xi_{i,k}\} \right] 
\leq \tfrac{n^2}{\eta^2} \mathbb{E}\left[ \large|{\tilde{f}_i}(x_{i,k}+v_{i,k},\xi_{i,k}) - \mathbb E_{\hat v}[{\tilde{f}_i}(x_{i,k}+\hat v,\xi_{i,k})] \large|^2 \mid \mathcal{F}_k\cup\{\xi_{i,k}\} \right].
\end{align*}
Next, we define $h_i(\bar v) = {\tilde{f}_i}(x_{i,k}+ \eta \bar v, \xi_{i,k})$, where $\bar v$ is uniformly distributed on \fyy{the} unit sphere. \fyy{Let} $\bar v_{i,k}$ \fyy{denote a sample withdrawn} uniformly from \fyy{the} unit sphere. \fyy{Then,} by employing Lemma \ref{lem:Levy} and the fact that $\mathbb E_{v_{i,k}}\left[h_i(\tfrac{v_{i,k}}{\eta})\right]=\mathbb E_{\bar v_{i,k}}[h_i(\bar v_{i,k})]$, we have
\begin{align*}
&\mathbb{E}_{v_{i,k}}\left[ \large|{\tilde{f}_i}(x_{i,k}+v_{i,k},\xi_{i,k}) - \mathbb E_{\hat v}[{\tilde{f}_i}(x_{i,k}+\hat v,\xi_{i,k})] \large|^2 \fyy{\mid \mathcal{F}_k}\right]\\
&= \mathbb{E}_{\bar v_{i,k}}\left[ \large|{\tilde{f}_i}(x_{i,k}+\eta \bar v_{i,k},\xi_{i,k}) - \mathbb E_{\bar v}[{\tilde{f}_i}(x_{i,k}+\eta \bar v,\xi_{i,k})] \large|^2 \fyy{\mid \mathcal{F}_k}\right]\\
&= \int_0^{\fyy{\infty}}\mathbb P\left(\large|h_i(\bar v_{i,k}) - \mathbb E_{\bar v}[h_i(\bar v)] \large|^2\geq \alpha \fyy{\mid \mathcal{F}_k}\right) d\alpha \\
&= \int_0^{\fyy{\infty}}\mathbb P\left(\large|h_i(\bar v_{i,k}) - \mathbb E_{\bar v}[h_i(\bar v)]  \large|\geq \sqrt{\alpha} \fyy{\mid \mathcal{F}_k}\right) d\alpha \\
&\stackrel{\text{Lemma }\ref{lem:Levy}}{\leq} \int_0^{\fyy{\infty}} 2\sqrt{2\pi} e^{-\tfrac{n\alpha}{8(\eta L_0(\xi_{i,k}))^2}}d\alpha = \tfrac{16 \sqrt{2\pi} (\eta L_0(\xi_{i,k}))^2}{n}.
\end{align*}

Then we obtain
\begin{align*}
\mathbb{E}\left[\left\|g^{\eta}_{i,k}\right\|^2 \mid  \mathcal{F}_k\cup\{v_{i,k}\} \right] 
\leq \tfrac{n^2}{\eta^2}\tfrac{16 \sqrt{2\pi} (L_0(\xi_{i,k})\eta)^2}{n} \stackrel{\text{Assumption }\ref{assum:Heterogeneous}\text{(i)}}{\leq} 16\sqrt{2\pi}L_0^2n.
\end{align*}
Using the preceding inequality and the fact that $\|a+b\|^2\leq 2\|a\|^2+2\|b\|^2$ \fyy{for any $a,b \in \mathbb{R}^n$}, we obtain
$\mathbb{E}\left[\|g^{\eta}_{i,k}+ \nabla^{\eta,d}_{i,k}\|^2 \mid \mathcal{F}_k\right] \leq  32\sqrt{2\pi}L_0^2n+2\|\nabla^{\eta,d}_{i,k}\|^2$.
}

\noindent (iii) Invoking \yqro{(i) and (ii)}, and using equation~\eqref{prob:main_unconstrained}, we have
\begin{align*}
&\mathbb{E}\left[\|(g^{\eta}_{i,k}+ \nabla^{\eta,d}_{i,k}) -  \nabla {\bf f}_i^{\eta}(x_{i, k})\|^2 \mid \mathcal{F}_k\right]  =\mathbb{E}\left[\|g^{\eta}_{i,k}+\nabla^{\eta,d}_{i,k}  -  \nabla f^{\eta}_{i}(x_{i, k})-\nabla^{\eta,d}_{i,k} \|^2 \mid \mathcal{F}_k\right]  \\
& = \mathbb{E}\left[\|g^{\eta}_{i,k}\|^2 +  \|\nabla f_{i}^\eta(x_{i, k})\|^2\mid \mathcal{F}_k\right]  - 2\mathbb{E}\left[  {g^{\eta}_{i,k}}^{\fyy{\top}} \nabla f_{i}^\eta(x_{i, k})   \mid \mathcal{F}_k\right] \\
&= \mathbb{E}\left[\|g^{\eta}_{i,k}\|^2\mid \mathcal{F}_k\right] +  \mathbb{E}\left[\| \nabla f_{i}^\eta(x_{i, k})\|^2\mid \mathcal{F}_k \right]  - 2\mathbb{E}\left[\| \nabla f_{i}^\eta(x_{i, k})\|^2 \mid \mathcal{F}_k\right] \leq \yqro{16\sqrt{2\pi}L_0^2n}.
\end{align*}

\begin{lemma}[Aggregated zeroth-order gradient]\label{lem:gbar}  \em                                                          
Let Assumption~\ref{assum:Heterogeneous} hold. Consider Algorithm~\ref{alg:PZO_LSGD}. Then, for all $k\geq 0$, {the following results hold}:
\begin{flalign*}
&\mbox{(i)}\ \   \mathbb{E}\left[\|\bar g_k^{\eta} +\bar{\nabla}_{k}^{\eta,d}\|^2 \mid \mathcal{F}_k \right] 
{\leq}\yqro{16\sqrt{2\pi}L_0^2n}  +2\left\|\nabla {\bf f}^\eta (\bar x_{k})\right\|^2  +\yqro{\tfrac{2(L_0\sqrt{n}+1)^2}{\eta^2}} \bar e_k.&& \\
&\mbox{(ii) } \mathbb{E}\left[ \nabla {\bf f}^\eta(\bar x_k)^{\fyy{\top}}(\bar{g}_k^{\eta}+\bar{\nabla}_{k}^{\eta,d})\mid \mathcal{F}_k \right]  \geq \tfrac{1}{2}\|\nabla {\bf f}^\eta(\bar x_k)\|^2 - \yqro{\tfrac{(L_0\sqrt{n}+1)^2}{2m\eta^2}}\textstyle\sum_{i=1}^{m} \|x_{i, k}-\bar x_{k} \|^2 .&&\end{flalign*}
\end{lemma}
\proof{Proof.} 
\noindent (i)  Using the definition of $\bar g^\eta_{k}$, we can write
\begin{align}                                                                                     
&\mathbb{E}\left[\|\bar g_k^{\eta}+\bar{\nabla}_{k}^{\eta,d}\|^2\mid  \mathcal{F}_k  \right] = \mathbb{E}\left[\left\|\tfrac{1}{m}\textstyle\sum_{i=1}^{m}(g_{i, k}^{\eta}+\nabla_{i,k}^{\eta,d})\right\|^2\mid  \mathcal{F}_k \right] \\
&{=} 
\mathbb{E}\left[\left\|\tfrac{1}{m}\textstyle\sum_{i=1}^{m}\left(g_{i, k}^{\eta}+\nabla_{i,k}^{\eta,d}-\nabla {\bf f}_i^\eta (x_{i,k})+\nabla {\bf f}_i^\eta (x_{i,k}) \right)\right\|^2\mid  \mathcal{F}_k \right] \nonumber \\
&\stackrel{\tiny \mbox{Lemma\mbox{ }\ref{hetero:local:gradient} {(i)}}}{=}\mathbb{E} \left[\left\|\tfrac{1}{m}\textstyle\sum_{i=1}^{m}\left(g_{i, k}^{\eta}+\nabla_{i,k}^{\eta,d}-\nabla {\bf f}_i^\eta (x_{i,k}) \right)\right\|^2 \mid  \mathcal{F}_k\right]+ \left\|\tfrac{1}{m}\textstyle\sum_{i=1}^{m}\nabla {\bf f}_i^\eta (x_{i,k})\right\|^2 \nonumber \\
&\stackrel{\tiny \mbox{Lemma\mbox{ }\ref{hetero:local:gradient} {(i)}}}{=} \tfrac{1}{m^2}\sum_{i=1}^{m}\mathbb{E}\left[\left\|g_{i, k}^{\eta}+\nabla_{i,k}^{\eta,d}-\nabla {\bf f}_i^\eta (x_{i,k})\right\|^2 \mid  \mathcal{F}_k \right]\notag\\
&+ \left\|\nabla {\bf f}^\eta (\bar x_{k})+\tfrac{1}{m}\textstyle\sum_{i=1}^{m}\left(\nabla {\bf f}_i^\eta (x_{i,k}) -\nabla {\bf f}_i^\eta (\bar x_k)\right) \right\|^2  \nonumber\\
&\stackrel{\tiny \mbox{Lemma\ \ref{hetero:local:gradient} (iii)}}{\leq}  \yqro{\tfrac{16\sqrt{2\pi}L_0^2n}{m}} +2\left\|\nabla {\bf f}^\eta (\bar x_{k})\right\|^2  +2 \left\|\tfrac{1}{m}\textstyle\sum_{i=1}^{m}\left(\nabla {\bf f}_i^\eta (x_{i,k}) -\nabla {\bf f}_i^\eta (\bar x_k)\right)\right\|^2. \nonumber
\end{align}
Note that given vectors $y_i \in \mathbb{R}^n$ for $i \in [m]$, we have $\| \tfrac{1}{m}\textstyle\sum_{i=1}^{m}y_i\|^2\leq \tfrac{1}{m}\textstyle\sum_{i=1}^{m}\|y_i\|^2$. Utilizing this inequality together with Lemma~\ref{SphericalSmooth} (iv), we obtain
\begin{align*}  
&\mathbb{E}\left[\|\bar g_k^{\eta}+\bar{\nabla}_{k}^{\eta,d}\|^2\mid  \mathcal{F}_k  \right]{\leq}\yqro{\tfrac{16\sqrt{2\pi}L_0^2n}{m}}+2\left\|\nabla {\bf f}^\eta (\bar x_{k})\right\|^2    + \yqro{\tfrac{2(L_0\sqrt{n}+1)^2}{\eta^2m}}  \textstyle\sum_{i=1}^{m}\left\|x_{i,k} -\bar x_k\right\|^2.
\end{align*}
Recalling the definition of $\bar e_k$, we obtain the {required} bound.

\noindent (ii) We have
\begin{align}
&\mathbb{E}\left[ \nabla {\bf f}^\eta(\bar x_k)^{\fyy{\top}}(\bar{g}_k^{\eta}+\bar{\nabla}_{k}^{\eta,d})\mid \mathcal{F}_k \right] =\nabla {\bf f}^\eta(\bar x_k)^{\fyy{\top}} \mathbb{E}\left[ \tfrac{1}{m}\textstyle\sum_{i=1}^{m}(g_{i, k}^{\eta}+\bar{\nabla}_{i,k}^{\eta,d})\mid \mathcal{F}_k \right] \nonumber \\
&\stackrel{\tiny \mbox{Lemma\mbox{ }\ref{hetero:local:gradient} {(i)}}}{=}
\nabla {\bf f}^\eta(\bar x_k)^{\fyy{\top}} \left(\tfrac{1}{m}\textstyle\sum_{i=1}^{m}\nabla {\bf f}_i^\eta (x_{i, k})  \right)  \nonumber \\
&=  \nabla {\bf f}^\eta(\bar x_k)^{\fyy{\top}} \tfrac{1}{m}\textstyle\sum_{i=1}^{m}\left(\nabla {\bf f}_i^\eta (x_{i, k}) -\nabla {\bf f}_i^\eta (\bar x_{k}) +\nabla {\bf f}_i^\eta (\bar x_{k}) \right)  \nonumber \\
&= \nabla {\bf f}^\eta(\bar x_k)^{\fyy{\top}} \tfrac{1}{m}\textstyle\sum_{i=1}^{m}\left(\nabla {\bf f}_i^\eta (x_{i, k})-\nabla {\bf f}_i^\eta (\bar x_{k}) \right)  +\|\nabla {\bf f}^\eta(\bar x_{k})\|^2 \nonumber \\
& {\geq}-\tfrac{1}{2}\|\nabla {\bf f}^\eta(\bar x_k)\|^2 - \tfrac{1}{2} \left\|\tfrac{1}{m}\textstyle\sum_{i=1}^{m}\left(\nabla {\bf f}_i^\eta (x_{i, k})-\nabla {\bf f}_i^\eta (\bar x_{k})\right) \right\|^2  +\|\nabla {\bf f}^\eta(\bar x_{k})\|^2  \nonumber \\
& {\geq}\tfrac{1}{2}\|\nabla {\bf f}^\eta(\bar x_k)\|^2 - \tfrac{1}{2m} \textstyle\sum_{i=1}^{m}\left\| \nabla {\bf f}_i^\eta (x_{i, k})-\nabla {\bf f}_i^\eta (\bar x_{k})  \right\|^2  \nonumber \\
&\stackrel{\tiny \mbox{Lemma\ \ref{SphericalSmooth} (iv)}}{\geq}\tfrac{1}{2}\|\nabla {\bf f}^\eta(\bar x_k)\|^2 - \yqro{\tfrac{(L_0\sqrt{n}+1)^2}{2m\eta^2}}\textstyle\sum_{i=1}^{m} \|x_{i, k}-\bar x_{k} \|^2  \nonumber .
\end{align}
The bound is obtained by recalling the definition of $\bar e_k$.

\begin{lemma}\label{lem:bound_on_dist2_ave} \em 
Let Assumption~\ref{assum:Heterogeneous} hold. Consider Algorithm~\ref{alg:PZO_LSGD}. {For any $k$}, we have
 \begin{align*}
\tfrac{1}{m}\textstyle\sum_{i=1}^{m}\|\nabla^{\eta,d}_{i,k}\|^2 \leq
 8\bar{e}_k+\tfrac{2B_1^2}{\eta^2} + 4B_2^2 \left\|\nabla {\bf f}^{\eta}(\bar{x}_k)\right\|^2+ \yqro{64B_2^2\sqrt{2\pi}L_0^2n}.
\end{align*}
\end{lemma} 

\proof{Proof.} Invoking Definition~\ref{def:basic_terms}, we can write 
\begin{align*}
\tfrac{1}{m}\textstyle\sum_{i=1}^{m}\|\nabla^{\eta,d}_{i,k}\|^2 &=\tfrac{1}{m}\textstyle\sum_{i=1}^{m}\|\nabla^{\eta,d}_{i,k}-\nabla d_i^\eta(\bar x_k)+\nabla d_i^\eta(\bar x_k)\|^2\\
& \leq \tfrac{2}{m}\textstyle\sum_{i=1}^{m}\|\nabla^{\eta,d}_{i,k}-\nabla d_i^\eta(\bar x_k)\|^2+\tfrac{2}{m}\textstyle\sum_{i=1}^{m}\|\nabla d_i^\eta(\bar x_k)\|^2\\
& =\tfrac{2}{m}\textstyle\sum_{i=1}^{m}\|(x_{i,k}-\bar{x}_k)-(\mathcal{P}_{X_i}(x_{i,k}) -\mathcal{P}_{X_i}(\bar x_k))\|^2+\tfrac{2}{m\eta^2}\textstyle\sum_{i=1}^{m} \mbox{dist}^2(\bar{x}_k,X_i) \\
& \leq 8\bar{e}_k+\tfrac{2}{m\eta^2}\textstyle\sum_{i=1}^{m} \mbox{dist}^2(\bar{x}_k,X_i) ,
\end{align*}
where the last inequality is obtained using the nonexpansiveness of the projection operator. Employing {Assumption \ref{assum:Heterogeneous} (iii)}, we obtain
\begin{align*}
\tfrac{1}{m}\textstyle\sum_{i=1}^{m}\|\nabla^{\eta,d}_{i,k}\|^2 & \leq 8\bar{e}_k+\tfrac{2B_1^2}{\eta^2}+\tfrac{2B_2^2}{\eta^2}\left\|\bar{x}_k-\tfrac{1}{m}\textstyle\sum_{i=1}^m\mathcal{P}_{X_i}(\bar{x}_k)\right\|^2\\
& \leq 8\bar{e}_k+\tfrac{2B_1^2}{\eta^2} + 2B_2^2 \left\|\tfrac{1}{m}\textstyle\sum_{i=1}^m\tfrac{1}{\eta}\left(\bar{x}_k-\mathcal{P}_{X_i}(\bar{x}_k)\right)\right\|^2\\
& = 8\bar{e}_k+\tfrac{2B_1^2}{\eta^2} + 2B_2^2 \left\|\tfrac{1}{m}\textstyle\sum_{i=1}^m\nabla d_i^\eta(\bar{x}_k) \right\|^2\\
& = 8\bar{e}_k+\tfrac{2B_1^2}{\eta^2} + 2B_2^2 \left\|\tfrac{1}{m}\textstyle\sum_{i=1}^m\left(\nabla {\bf f}_i^{\eta}(\bar{x}_k)- \nabla f_i^{\eta}(\bar{x}_k)\right)\right\|^2\\
& \leq 8\bar{e}_k+\tfrac{2B_1^2}{\eta^2} + 4B_2^2 \left\|\nabla {\bf f}^{\eta}(\bar{x}_k)\right\|^2+ 4B_2^2 \left\|\nabla f^{\eta}(\bar{x}_k) \right\|^2.
\end{align*}
Next, we find an upper bound on $\left\|\nabla f^{\eta}(\bar{x}_k) \right\|^2$. Using the definition of $f(x)$ and invoking Lemma~\ref{SphericalSmooth}, for any $x$ we have
 \begin{align*}
 \left\|\nabla f^{\eta}(x) \right\|^2 &\leq \tfrac{1}{m}\textstyle\sum_{i=1}^m\left\|\nabla f_i^{\eta}(x) \right\|^2 \\
 & =\left(\tfrac{n^2}{\eta^2}\right)  \tfrac{1}{m}\textstyle\sum_{i=1}^m\left\|\mathbb{E} [(  f_i(x+\yqro{v_i})-f_i(x\yqro{-v_i}))  \tfrac{v_i}{\|\yqro{2}v_i\|}  ]\right\|^2\\
&\stackrel{\tiny \mbox{Jensen's ineq.}}{\leq}\left(\tfrac{n^2}{\yqro{4}\eta^2}\right)  \tfrac{1}{m}\textstyle\sum_{i=1}^m\mathbb{E} [|f_i(x+\yqro{v_i})-f_i(x\yqro{-v_i})) |^2 ].
\end{align*}
\yqro{Next, we note that $f_i$ is $L_0$-Lipschitz, then we utilize Lemma \ref{lem:Levy} by following the same approach in Lemma~\ref{hetero:local:gradient}, we obtain
\begin{align}\label{eqn:bound:f:eta}
\left\|\nabla f^{\eta}(x) \right\|^2 \leq 16\sqrt{2\pi}L_0^2n.
\end{align}}
From the two preceding inequalities, we obtain the result.
\begin{lemma}\label{lem:aggr_a_k_extension}\em
Suppose for $T_r+1\leq k \leq T_{r+1}$, the nonnegative sequences $\{a_k\}$ and $\{\theta_k\}$ satisfy a recursive relation of the form $a_{k}\leq (k-T_r)\gamma^2\textstyle\sum_{t=T_r}^{k-1}(\beta a_t+\theta_t),$ where $a_{T_r}=0$ and $\beta>0$. Then, $a_k \leq H\gamma^2 \sum_{t=T_r}^{k-1}(\beta H\gamma^2+1)^{k-t-1}\theta_{t}$ for all $T_r+1\leq k \leq T_{r+1}$. Moreover, if $0<\gamma \leq \frac{1}{\sqrt{\beta}H}$, then $a_k\leq 3H\gamma^2 \sum_{t=T_r}^{k-1}\theta_t.$
\end{lemma}
\proof{Proof.}
This can be shown using induction by unrolling $a_{k}\leq (k-T_r)\gamma^2\textstyle\sum_{t=T_r}^{k-1}(\beta a_t+\theta_t)$ recursively. The proof is omitted.

\begin{lemma}[Bound on average consensus violation]\label{noncvx:hetero:deviation}\em                                        
Consider Algorithm~\ref{alg:PZO_LSGD}. Let Assumption~\ref{assum:Heterogeneous} hold and let $H\geq 1$ be given by Definition~\ref{def:H}. Then, for any communication round $r>0$, for all $T_r \leq k \leq T_{r+1}-1$ we have
\begin{align*}
\mathbb{E}\left[  \bar e_k\right]\  & \leq  \gamma^2(k-T_r)\textstyle\sum_{t=T_r}^{k-1}  \left(	\yqro{32\sqrt{2\pi}L_0^2n}+	\frac{4B_1^2}{\eta^2}+\yqro{128B_2^2\sqrt{2\pi}L_0^2n} +16\mathbb{E}\left[  \bar e_t\right] + 8B_2^2 \mathbb{E}\left[\left\|\nabla {\bf f}^{\eta}(\bar{x}_t)\right\|^2\right] \right).
\end{align*}
Moreover, if $0<\gamma \leq \frac{1}{4 H}$, then 
\begin{align*}\mathbb{E}\left[  \bar e_k\right] \leq  3H^2\gamma^2\left(\yqro{32\sqrt{2\pi}L_0^2n}+	\tfrac{4B_1^2}{\eta^2}+\yqro{128B_2^2\sqrt{2\pi}L_0^2n} \right) + 24B_2^2H\gamma^2 \textstyle\sum_{t=T_r}^{k-1}   \mathbb{E}\left[\left\|\nabla {\bf f}^{\eta}(\bar{x}_t)\right\|^2\right].\end{align*}
\end{lemma}
\proof{Proof.}
In view of Algorithm~\ref{alg:PZO_LSGD}, for any $i$ at any communication round $r>0$, for all $T_r \leq k \leq T_{r+1}-1$ we have $x_{i,k+1} =  x_{i,k} - \gamma  (g^{\eta}_{i,k}+\nabla^{\eta,d}_{i,k}).$
Equivalently, we can write $x_{i,k} =  x_{i,k-1} - \gamma  (g^{\eta}_{i,k-1}+\nabla^{\eta,d}_{i,k-1})$ for all $T_r+1 \leq k \leq T_{r+1}.$ This implies that 
\begin{align}\label{eqn:nonrecursive_local}
x_{i,k} =  x_{i,T_r} - \gamma  \textstyle\sum_{t=T_r}^{k-1}(g^{\eta}_{i,t}+\nabla^{\eta,d}_{i,t}),\quad \hbox{for all } T_r+1 \leq k \leq T_{r+1}.
\end{align} 
Again from Algorithm~\ref{alg:PZO_LSGD}, we have $\hat{x}_r =x_{i,T_r}$. From the definition of $\bar{x}_k$, we can write $\bar{x}_{T_r} = \hat{x}_r$. This implies that $\bar{x}_{T_r} =x_{i,T_r}$ for all $i$ and $r$. In view of Lemma~\ref{lem:bar_x_recursion}, we have
\begin{align}\label{eqn:nonrecursive_ave}
\bar x_k =  x_{i,T_r} - \gamma  \textstyle\sum_{t=T_r}^{k-1}(\bar{g}^{\eta}_{t}+\bar{\nabla}^{\eta,d}_{t}),\quad \hbox{for all } T_r+1 \leq k \leq T_{r+1}.
\end{align}
Utilizing \eqref{eqn:nonrecursive_local} and \eqref{eqn:nonrecursive_ave}, for all $T_r+1 \leq k \leq T_{r+1}$ we have
\begin{align*}
\mathbb{E}\left[  \bar e_k  \mid \mathcal{F}_{T_r}\right]  
&=\tfrac{1}{m}\textstyle\sum_{i=1}^{m}\mathbb{E}\left[ \|x_{i, k}-\bar x_k \|^2  \mid \mathcal{F}_{T_r} \right]   \\
&=\tfrac{1}{m}\textstyle\sum_{i=1}^{m}\mathbb{E}\left[ \left\|\gamma\textstyle\sum_{t=T_r}^{k-1} (g_{i, t}^{\eta} +\nabla^{\eta,d}_{i,t})-\gamma\textstyle\sum_{t=T_r}^{k-1} (\bar{g}_{t}^{\eta}+ \bar{\nabla}^{\eta,d}_{t}) \right\|^2 \mid \mathcal{F}_{T_r}  \right]   \\
&\leq \tfrac{\gamma^2(k-T_r)}{m}\textstyle\sum_{i=1}^{m}\textstyle\sum_{t=T_r}^{k-1} \mathbb{E}\left[ \left\| (g_{i, t}^{\eta} +\nabla^{\eta,d}_{i,t})- (\bar{g}_{t}^{\eta}+ \bar{\nabla}^{\eta,d}_{t}) \right\|^2 \mid \mathcal{F}_{T_r}  \right],
\end{align*}
where the preceding relation is implied {by} the inequality $\|  \textstyle\sum_{t=1}^{T}y_t\|^2\leq T\textstyle\sum_{t=1}^{T}\|y_t\|^2$ for any $y_t \in \mathbb{R}^n$ for $t \in [T]$. We have
\begin{align*}
\mathbb{E}\left[  \bar e_k  \mid \mathcal{F}_{T_r}\right]  
&\leq \tfrac{\gamma^2(k-T_r)}{m}\textstyle\sum_{t=T_r}^{k-1}\textstyle\sum_{i=1}^{m} \mathbb{E}\left[ \left\| (g_{i, t}^{\eta} +\nabla^{\eta,d}_{i,t})- (\bar{g}_{t}^{\eta}+ \bar{\nabla}^{\eta,d}_{t}) \right\|^2 \mid \mathcal{F}_{T_r}  \right]\\
& = \tfrac{\gamma^2(k-T_r)}{m}\textstyle\sum_{t=T_r}^{k-1}\sum_{i=1}^{m}\left(\mathbb{E}\left[ \left\|  g_{i, t}^{\eta} +\nabla^{\eta,d}_{i,t} \right\|^2 \mid \mathcal{F}_{T_r}  \right] +\mathbb{E}\left[ \left\| \bar{g}_{t}^{\eta}+ \bar{\nabla}^{\eta,d}_{t}  \right\|^2 \mid \mathcal{F}_{T_r}  \right] \right)\\
& -2\tfrac{\gamma^2(k-T_r)}{m}\textstyle\sum_{t=T_r}^{k-1}\sum_{i=1}^{m}\mathbb{E}\left[ (  g_{i, t}^{\eta} +\nabla^{\eta,d}_{i,t} )^{\fyy{\top}}(\bar{g}_{t}^{\eta}+ \bar{\nabla}^{\eta,d}_{t}  ) \mid \mathcal{F}_{T_r}  \right],
\end{align*}
Observing that 
$$\tfrac{1}{m}\textstyle\sum_{i=1}^{m}\mathbb{E}\left[ (  g_{i, t}^{\eta} +\nabla^{\eta,d}_{i,t} )^{\fyy{\top}}(\bar{g}_{t}^{\eta}+ \bar{\nabla}^{\eta,d}_{t}  ) \mid \mathcal{F}_{T_r}  \right]= \mathbb{E}\left[\left\|\bar{g}_{t}^{\eta}+ \bar{\nabla}^{\eta,d}_{t}  \right\|^2 \mid \mathcal{F}_{T_r}  \right],$$
we obtain 
\begin{align*}
\mathbb{E}\left[  \bar e_k  \mid \mathcal{F}_{T_r}\right]  
&\leq  \tfrac{\gamma^2(k-T_r)}{m}\textstyle\sum_{t=T_r}^{k-1}\sum_{i=1}^{m} \mathbb{E}\left[ \left\|  g_{i, t}^{\eta} +\nabla^{\eta,d}_{i,t} \right\|^2 \mid \mathcal{F}_{T_r}  \right].
\end{align*}
From the law of total expectation, for any $T_r\leq  t\leq k-1 $, we can write
\begin{align*}
\mathbb{E}\left[ \left\| g_{i, t}^{\eta} +\nabla^{\eta,d}_{i,t} \right\|^2  \mid \mathcal{F}_{T_r}  \right]&=\mathbb{E}\left[\mathbb{E}\left[ \left\| g_{i, t}^{\eta}  +\nabla^{\eta,d}_{i,t}\right\|^2  \mid \mathcal{F}_{T_r}\cup\left( \cup_{i=1}^m\cup_{t=T_r}^{t-1} \{ \xi_{i, t}, v_{i, t} \} \right)  \right] \right]\\
& =\mathbb{E}\left[\mathbb{E}\left[ \left\| g_{i, t}^{\eta}  +\nabla^{\eta,d}_{i,t}\right\|^2  \mid \mathcal{F}_t  \right] \right]\\
& \stackrel{\tiny\mbox {Lemma \ref{hetero:local:gradient} {(ii)}}}{\leq }\yqro{32\sqrt{2\pi}L_0^2n} +2\|\nabla^{\eta,d}_{i,t}\|^2.
 \end{align*}
From the two preceding relations, we obtain
\begin{align*}
\mathbb{E}\left[  \bar e_k  \mid \mathcal{F}_{T_r}\right]\  & \leq \tfrac{\gamma^2(k-T_r)}{m}\textstyle\sum_{t=T_r}^{k-1}\sum_{i=1}^{m} \left(\yqro{32\sqrt{2\pi}L_0^2n} +2\|\nabla^{\eta,d}_{i,t}\|^2\right).
\end{align*}
Invoking Lemma~\ref{lem:bound_on_dist2_ave}, we obtain 
\begin{align*}
\mathbb{E}\left[  \bar e_k  \mid \mathcal{F}_{T_r}\right]\  & \leq  \gamma^2(k-T_r)\textstyle\sum_{t=T_r}^{k-1}  \left(\yqro{32\sqrt{2\pi}L_0^2n}+	\frac{4B_1^2}{\eta^2}+\yqro{128B_2^2\sqrt{2\pi}L_0^2n} +16\bar{e}_t + 8B_2^2 \left\|\nabla {\bf f}^{\eta}(\bar{x}_t)\right\|^2 \right).
\end{align*}
Taking expectations on both sides, we obtain the first result.  The second bound is obtained by invoking Lemma~\ref{lem:aggr_a_k_extension}.

\subsection{Proof of Proposition \ref{thm:main_bound}}
\proof{Proof.} 
\noindent (i) Recall from Lemma~\ref{SphericalSmooth} that each of the local functions $f_i^{\eta}$ is \yqro{$\tfrac{L_0\sqrt{n}}{\eta}$}-smooth. Also, $d_i^{\eta}$ is $\frac{1}{\eta}$-smooth. As such, ${\bf f}^\eta$ is \yqro{$\left(\frac{L_0\sqrt{n}+1}{\eta}\right)$}-smooth. Invoking Lemma~\ref{lem:bar_x_recursion}, we {may obtain}
\begin{align*}
  {\bf f}^{\eta}(\bar x_{k+1})    \leq   {\bf f}^{\eta}(\bar x_{k}) -\gamma   \nabla {\bf f}^{\eta}(\bar x_k)^{\fyy{\top}} \left(\bar{g}_k^{\eta} +\bar{\nabla}_{i,k}^{\eta,d}\right)  +\yqro{\tfrac{L_0\sqrt{n}}{2\eta}}\gamma^2 \|\bar{g}_k^{\eta}+\bar{\nabla}_{i,k}^{\eta,d}\|^2.
\end{align*}
Taking expectations on both sides, we obtain 
\begin{align*}
 \mathbb{E}\left[ {\bf f}^{\eta}(\bar x_{k+1})  \right]  \leq   \mathbb{E}\left[ {\bf f}^{\eta}(\bar x_{k})\right] -\gamma    \mathbb{E}\left[\nabla {\bf f}^{\eta}(\bar x_k)^{\fyy{\top}} \left(\bar{g}_k^{\eta} +\bar{\nabla}_{i,k}^{\eta,d}\right) \right] +\yqro{\tfrac{L_0\sqrt{n}}{2\eta}}\gamma^2  \mathbb{E}\left[\|\bar{g}_k^{\eta}+\bar{\nabla}_{i,k}^{\eta,d}\|^2\right].
\end{align*}
Invoking Lemma~\ref{lem:gbar}, we obtain 
\begin{align*}
 \mathbb{E}\left[ {\bf f}^{\eta}(\bar x_{k+1})  \right]  &\leq   \mathbb{E}\left[ {\bf f}^{\eta}(\bar x_{k})\right] -\gamma    \left(\tfrac{1}{2} \mathbb{E}\left[ \|\nabla {\bf f}^\eta(\bar x_k)\|^2\right] - \yqro{\tfrac{(L_0\sqrt{n}+1)^2}{2 \eta^2}} \mathbb{E}\left[ \bar e_k\right]\right) \\
&  +\yqro{\tfrac{L_0\sqrt{n}}{2\eta}}\gamma^2\left(\tfrac{\yqro{16\sqrt{2\pi}L_0^2n}}{m} +2 \mathbb{E}\left[ \left\|\nabla {\bf f}^\eta (\bar x_{k})\right\|^2\right]  +\yqro{\tfrac{2(L_0\sqrt{n}+1)^2}{\eta^2}}  \mathbb{E}\left[ \bar e_k\right]\right) .
\end{align*}
Using $\gamma \leq \yqro{\tfrac{\eta}{4L_0\sqrt{n}}}$ and rearranging the terms, we have
\begin{align*}
 \tfrac{\gamma}{4} \mathbb{E}\left[\|\nabla {\bf f}^\eta(\bar x_k)\|^2\right]  &\leq   \mathbb{E}\left[ {\bf f}^{\eta}(\bar x_{k})\right] - \mathbb{E}\left[ {\bf f}^{\eta}(\bar x_{k+1})  \right]   + \yqro{\tfrac{\gamma^2 16\sqrt{2\pi}L_0^3n^{1.5}}{2\eta m}}    +\yqro{\tfrac{3\gamma(L_0\sqrt{n}+1)^2}{4\eta^2}}  \mathbb{E}\left[ \bar e_k\right]  .
\end{align*}
Summing both sides over $k=1=0,\ldots,K$, then dividing both sides by {$\tfrac{\gamma(K+1)}{4}$}, and using the definition of $k^*$, we obtain
\begin{align*}
  \mathbb{E}\left[\|\nabla {\bf f}^\eta(\bar x_{k^*})\|^2\right]  &\leq  \tfrac{4( \mathbb{E}\left[ {\bf f}^{\eta}(\bar x_{0})\right] - \mathbb{E}\left[ {\bf f}^{\eta}(\bar x_{K+1})  \right])}{\gamma(K+1)}   + \yqro{\tfrac{\gamma 32\sqrt{2\pi} L_0^3n^{1.5}}{\eta m}} + \yqro{\tfrac{3(L_0\sqrt{n}+1)^2}{\eta^2(K+1)}} \textstyle\sum_{k=0}^K \mathbb{E}\left[ \bar e_k\right]  .
\end{align*}
From Lemma~\ref{noncvx:hetero:deviation} and the definition of $k^*$, we obtain 
 \begin{align*}
 \tfrac{1}{K+1}\textstyle\sum_{k=0}^K \mathbb{E}\left[  \bar e_k\right]& \leq 3H^2\gamma^2\left(\yqro{32\sqrt{2\pi}L_0^2n+\tfrac{4B_1^2}{\eta^2}+128B_2^2\sqrt{2\pi}L_0^2n }\right) \\
 &+ 24B_2^2H\gamma^2 \tfrac{1}{K+1}\textstyle\sum_{k=0}^K\textstyle\sum_{t=T_r}^{k-1}   \mathbb{E}\left[\left\|\nabla {\bf f}^{\eta}(\bar{x}_t)\right\|^2\right]\\
 &\leq 3H^2\gamma^2\left(\yqro{32\sqrt{2\pi}L_0^2n+\tfrac{4B_1^2}{\eta^2}+128B_2^2\sqrt{2\pi}L_0^2n } \right) + 24B_2^2H^2\gamma^2    \mathbb{E}\left[\left\|\nabla {\bf f}^{\eta}(\bar{x}_{k^*})\right\|^2\right],
 \end{align*}
 where in the preceding relation,  we used 
 $$\textstyle\sum_{k=0}^K\textstyle\sum_{t=T_r}^{k-1}   \mathbb{E}\left[\left\|\nabla {\bf f}^{\eta}(\bar{x}_t)\right\|^2\right] \leq H\textstyle\sum_{k=0}^K    \mathbb{E}\left[\left\|\nabla {\bf f}^{\eta}(\bar{x}_k)\right\|^2\right].$$
Thus, invoking $\gamma \leq \yqro{\tfrac{\eta}{12\sqrt{3} B_2(L_0\sqrt{n}+1)H}}$, from the preceding relations, we obtain
 \begin{align*}
  \mathbb{E}\left[\|\nabla {\bf f}^\eta(\bar x_{k^*})\|^2\right]  &\leq  \tfrac{4( \mathbb{E}\left[ {\bf f}^{\eta}(\bar x_{0})\right] - \mathbb{E}\left[ {\bf f}^{\eta,*} \right])}{\gamma(K+1)}   + \yqro{\tfrac{\gamma 32\sqrt{2\pi} L_0^3n^{1.5}}{\eta m}} \\
  &  +\yqro{\tfrac{9H^2\gamma^2(L_0\sqrt{n}+1)^2}{\eta^2}}\left(\yqro{32\sqrt{2\pi}L_0^2n+\tfrac{4B_1^2}{\eta^2}+128B_2^2\sqrt{2\pi}L_0^2n } \right)   +0.5 \mathbb{E}\left[\|\nabla {\bf f}^\eta(\bar x_{k^*})\|^2\right].
\end{align*}
Rearranging the terms, we obtain the inequality in part (i).

\noindent (ii) Substituting $\gamma:=\yqro{\sqrt{\tfrac{m\eta}{n^{1.5}L_0^3K}}}$ and $H:= \sqrt[4]{\tfrac{K}{m^3}}$ in the preceding bound, we obtain
 \begin{align*}\mathbb{E}\left[\|\nabla {\bf f}^\eta(\bar x_{k^*})\|^2\right]   
 &\leq    \tfrac{8( \mathbb{E}\left[ {\bf f}^{\eta}(\bar x_{0})\right]-{\bf f}^{\eta,*})\yqro{\tfrac{n^{0.75}L_0^{1.5}}{\eta^{0.5}}}+\yqro{\tfrac{64\sqrt{2\pi} L_0^{1.5}n^{0.75}}{\eta^{0.5}}}+\yqro{(n^{-1.5})\tfrac{72 (L_0\sqrt{n}+1)^2}{L_0^3\eta}}\left(\yqro{\tfrac{B_1^2}{\eta^2}+8\sqrt{2\pi}(4B_2^2+1)L_0^2n } \right) }{\sqrt{ m K}}.
\end{align*}
This leads to iteration complexity of \yqro{$\mathcal{O}\left(\left(\tfrac{n^{0.75}L_0^{1.5}}{\eta^{0.5}}+\tfrac{L_0^{1.5}n^{0.75}}{\eta^{0.5}}+\tfrac{B_1^2L_0^{-1}n^{-0.5}}{\eta^3}+\tfrac{L_0 \sqrt n}{\eta}+\tfrac{B_2^2L_0 \sqrt n}{\eta} \right)^2\tfrac{1}{m\epsilon^2}\right)$, by combining terms and assuming $\tfrac{B_1^2L_0^{-1}n^{-0.5}}{\eta^3}\leq 1$, the iteration complexity can be written as $$\mathcal{O}\left(\left(\tfrac{L_0^{1.5}n^{0.75}}{\eta^{0.5}}+\tfrac{B_2^2L_0 \sqrt n}{\eta} \right)^2\tfrac{1}{m\epsilon^2}\right).$$}


\noindent (iii) From the choice of $H$ is (ii), we obtain $R=\mathcal{O}(\tfrac{K}{H}) = \mathcal{O}\left(\tfrac{K}{ \sqrt[4]{{K}/{m^3}} }\right) =\mathcal{O}\left((mK)^{3/4}\right)$.


\subsection{Proof of Proposition~\ref{proposition:1}}
\proof{Proof.}
(i) The proof for this part follows from \yqro{Theorem 3 in \cite{lin2022gradient}.}

(ii) From $\nabla \mathbf{f}^{\eta}(x) =0$, we have that $\nabla f^{\eta}(x) +\frac{1}{\eta}(x-\mathcal{P}_{X}(x)) =0$. This implies that $\|x-\mathcal{P}_{X}(x)\| \leq \eta \|\nabla f^{\eta}(x)\|$. \yqro{From \eqref{eqn:bound:f:eta}, we have $\|\nabla f^{\eta}(x)\| \leq 4\sqrt[4]{2\pi}L_0\sqrt{n}$}
Thus, the infeasibility of $x$ is bounded as $\|x-\mathcal{P}_X(x)\| \leq  \yqro{4\sqrt[4]{2\pi}\eta L_0\sqrt{n}}$. Recall that the $\delta$-Clarke generalized gradient of ${\mathbb{I}_X}$ at $x$  is defined as $\partial_{\delta} {\mathbb{I}_X}(x) \triangleq \mbox{conv}\left\{ \zeta: \zeta \in \mathcal{N}_X(y), \|x-y\| \leq \delta\right\},$ 
where $\mathcal{N}_X(\bullet)$ denotes the normal cone of $X$. In view of $\|x-\mathcal{P}_X(x)\| \leq  \yqro{4\sqrt[4]{2\pi}\eta L_0\sqrt{n}}$, for $y:=\mathcal{P}_X(x) $ and $\eta \leq \frac{\delta}{\max\{\yqro{1,4\sqrt[4]{2\pi}L_0\sqrt{n}}\}}$, we have $\|x-y\| \leq \delta$. Next we show that for  $\zeta:= \frac{1}{\eta}(x-\mathcal{P}_X(x))$ we have $\zeta \in \mathcal{N}_X(y)$. From the projection theorem, we may write $(x-\mathcal{P}_X(x))^{\fyy{\top}}(\mathcal{P}_X(x)-z) \geq 0,$ for all $z \in X. $ This implies that $\zeta^{\fyy{\top}}(y-z) \geq 0$ for all $z \in X$. Thus, we have $\zeta \in \mathcal{N}_X(y)$ which implies $\frac{1}{\eta}(x-\mathcal{P}_X(x)) \in \partial_{\delta}\mathbb{I}(x)$. From (i) and that $2\eta \leq \delta$, we have $\nabla f^{\eta}(x) \in \partial_{\delta} f(x)$. Adding the preceding relations and invoking $\nabla \mathbf{f}^{\eta}(x) =0$, we obtain $ 0 \in \partial_{\delta} \left(  f + \mathbb{I}_X \right)(x)$.  


{\begin{remark}   We note that the approximate Clarke stationary point is also referred to as Goldstein stationary point. (e.g. \cite{lin2022gradient}) \end{remark}}

\subsection{Proof of Lemma~\ref{lem:implicit_props}}
\proof{Proof.}
\noindent (i) From the strong convexity of $h(x,\bullet)$, we can write for any $x\in\mathbb{R}^n$ and all $y_1,y_2 \in \mathbb{R}^{\tilde{n}}$,
\begin{align}\label{eqn:s_monotone_h}
\mu_h\|y_1-y_2\|^2 \leq (y_1-y_2)^{\fyy{\top}}\left(\nabla_y h(x,y_1)-\nabla_y h(x,y_2)\right).
\end{align}
Substituting $x:=x_1$, $y_1:=y(x_1)$, and $y_2:=y(x_1)$ in \eqref{eqn:s_monotone_h}, we obtain 
$$\mu_h\|y(x_1)-y(x_2)\|^2 \leq (y(x_1)-y(x_2))^{\fyy{\top}}\left(\nabla_y h(x_1,y(x_1))-\nabla_y h(x_1,y(x_2))\right).$$
Similarly, substituting $x:=x_2$, $y_1:=y(x_1)$, and $y_2:=y(x_1)$ in \eqref{eqn:s_monotone_h}, we obtain 
$$\mu_h\|y(x_1)-y(x_2)\|^2 \leq (y(x_1)-y(x_2))^{\fyy{\top}}\left(\nabla_y h(x_2,y(x_1))-\nabla_y h(x_2,y(x_2))\right).$$
Adding the preceding two inequalities together, we have
\begin{align*}
2\mu_h\|y(x_1)-y(x_2)\|^2 &\leq (y(x_1)-y(x_2))^{\fyy{\top}}\left(\nabla_y h(x_1,y(x_1))-\nabla_y h(x_1,y(x_2))\right.\\&\left.+\nabla_y h(x_2,y(x_1))-\nabla_y h(x_2,y(x_2))\right).
\end{align*} 
Note that from the definition of $y(\bullet)$, we have $\nabla_y h(x_1,y(x_1))=\nabla_y h(x_2,y(x_2))=0$. As such, from the preceding inequality we have
\begin{align*}
2\mu_h\|y(x_1)-y(x_2)\|^2 &\leq (y(x_1)-y(x_2))^{\fyy{\top}}\left(-\nabla_y h(x_1,y(x_1))-\nabla_y h(x_1,y(x_2))\right.\\&\left.+\nabla_y h(x_2,y(x_1))+\nabla_y h(x_2,y(x_2))\right).
\end{align*} 
Using the Cauchy-Schwarz inequality and the triangle inequality, we obtain 
\begin{align*}
2\mu_h\|y(x_1)-y(x_2)\|^2 &\leq \|y(x_1)-y(x_2)\|\left\|\nabla_y h(x_2,y(x_1))-\nabla_y h(x_1,y(x_1))\right\|\\&+\|y(x_1)-y(x_2)\|\left\|\nabla_y h(x_2,y(x_2))-\nabla_y h(x_1,y(x_2))\right\|.
\end{align*} 
If $x_1=x_2$, the relation in (i) holds. Suppose $x_1\neq x_2$. Thus, $y(x_1)\neq y(x_2)$. We obtain 
\begin{align*}
2\mu_h\|y(x_1)-y(x_2)\| &\leq  \left\|\nabla_y h(x_2,y(x_1))-\nabla_y h(x_1,y(x_1))\right\| + \left\|\nabla_y h(x_2,y(x_2))-\nabla_y h(x_1,y(x_2))\right\|.
\end{align*} 
From Assumption~\ref{assum:bilevel}, we obtain 
$
2\mu_h\|y(x_1)-y(x_2)\|  \leq  L_{0,x}^{\nabla h}\left\|x_1-x_2\right\| + L_{0,x}^{\nabla h}\left\|x_1-x_2\right\|.
$ This implies the bound in (i).

{
\noindent(ii) Let $L_0^{\text{imp}}(\xi_i)$ denote the Lipschitz constant of ${\tilde{f}_i}(x,y(x),\xi_i)$. We have
\begin{align*}
 & |{\tilde{f}_i}(x_1,y(x_1),\xi_i) -{\tilde{f}_i}(x_2,y(x_2),\xi_i) |\\
 &= |{\tilde{f}_i}(x_1,y(x_1),\xi_i) -{\tilde{f}_i}(x_1,y(x_2),\xi_i)+{\tilde{f}_i}(x_1,y(x_2),\xi_i)-{\tilde{f}_i}(x_2,y(x_2),\xi_i) |\\
 & \leq |{\tilde{f}_i}(x_1,y(x_1),\xi_i) -{\tilde{f}_i}(x_1,y(x_2),\xi_i)|+| {\tilde{f}_i}(x_1,y(x_2),\xi_i)-{\tilde{f}_i}(x_2,y(x_2),\xi_i)|\\
 & \leq L_{0,y}^f(\xi_i)\|y(x_1)-y(x_2)\|+L_{0,x}^f(\xi_i)\|x_1-x_2\|\\
 & \leq \left( \tfrac{L_{0,y}^f(\xi_i) L_{0,x}^{\nabla h}}{\mu_h}+ L_{0,x}^f(\xi_i) \right)\|x_1-x_2\| = L_0^{\text{imp}}(\xi_i)\|x_1-x_2\|.
\end{align*}
}

\noindent (iii) 
{First, we show that for any $x_1,x_2,y$, we have $|f_i(x_1,y)-f_i(x_2,y)|\leq L_{0,x}^f\|x_1-x_2\|$. Also, for any $x,y_1,y_2$, we have $|f_i(x,y_1)-f_i(x,y_2)|\leq L_{0,y}^f\|y_1-y_2\|$. From Assumption~\ref{assum:bilevel}, we have
\begin{align*}
|f_i(x_1,y)-f_i(x_2,y)| &= |\mathbb{E}_{\xi_i} [{\tilde{f}_i}(x_1,y,\xi_i)-{\tilde{f}_i}(x_2,y,\xi_i)] |  \leq \mathbb{E}_{\xi_i}[| {\tilde{f}_i}(x_1,y,\xi_i)-{\tilde{f}_i}(x_2,y,\xi_i) |]\\  &\leq \mathbb{E}_{\xi_i}[L_{0,x}^f(\xi_i)]\|x_1-x_2 \|.
\end{align*}
From the Jensen's inequality, we have $(\mathbb{E}_{\xi_i}[L_{0,x}^f(\xi_i)])^2 \leq \mathbb{E}_{\xi_i}[{L_{0,x}^f}^2(\xi_i)].$ Therefore,
$$\mathbb{E}_{\xi_i}[L_{0,x}^f(\xi_i)] \leq |\mathbb{E}_{\xi_i}[L_{0,x}^f(\xi_i)]| \leq \sqrt{\mathbb{E}_{\xi_i}[({L_{0,x}^f}(\xi_i))^2]} {\leq} L_{0,x}^f.$$
From the preceding two relations, we obtain $
|f_i(x_1,y)-f_i(x_2,y)|\leq L_{0,x}^f\|x_1-x_2\|.$
Similarly, we obtain $
|f_i(x,y_1)-f_i(x,y_2)|\leq L_{0,y}^f\|y_1-y_2\|.$
} Next,  \fy{from the last two inequalities and} Assumption~\ref{assum:bilevel}, we obtain 
\begin{align*}
 | f (x_1)-  f (x_2)|& = |f(x_1,y(x_1)) -f(x_2,y(x_2)) |\\
 & =|f(x_1,y(x_1)) -f(x_1,y(x_2))+f(x_1,y(x_2))-f(x_2,y(x_2)) |\\
 & \leq |f(x_1,y(x_1)) -f(x_1,y(x_2))|+|f(x_1,y(x_2))-f(x_2,y(x_2)) |\\
 & {\leq \tfrac{1}{m}\sum_{i=1}^m|f_i(x_1,y(x_1)) -f_i(x_1,y(x_2))|+\tfrac{1}{m}\sum_{i=1}^m|f_i(x_1,y(x_2))-f_i(x_2,y(x_2))| }\\
 &  {\leq} L_{0,y}^f\|y(x_1)-y(x_2)\|+L_{0,x}^f\|x_1-x_2\|.
\end{align*}
The bound in (iii) is obtained by invoking the bound in (i).


\subsection{Preliminaries for proof of Theorem \ref{thm:bilevel}}
In this subsection, we provide preliminary results that will be utilized to prove Theorem~\ref{thm:bilevel} and provide some preliminary lemmas and their proofs.


\begin{definition}\label{def:basic_terms_bilevel}   
Let $x_{i,k}$ be given by Algorithm~\ref{alg:FedRiZO_upper} for all $i \in [m]$ and $k\geq 0$.  \fy{Let $d_i(x), \nabla_{i,k}^{\eta,d}, \bar{\nabla}_{k}^{\eta,d}, \bar{x}_k , $ and $\bar{e}_k$ be given by Definition~\ref{def:basic_terms}}.  Let us define  \fy{an average delay term as} $ \hat{e}_k \triangleq  \frac{\sum_{i=1}^m \|x_{i,k}-\hat{x}_k\|^2}{m}.$
\end{definition}
\begin{definition}\label{def:s_terms_bilevel} Let us define the following terms. 
\begin{align*}
& g^{\eta}_{i,k} \triangleq\tfrac{n}{2\eta^2} \left(\tilde{f}_i(x_{i,k}+v_{T_r},y(x_{i,k}+v_{T_r}),\xi_{i,k}) - \tilde{f}_i(x_{i,k}\yqro{-v_{T_r}},y(x_{i,k}\yqro{-v_{T_r}}),\xi_{i,k})\right)v_{T_r},\\
&{\hat g}^{\eta}_{i,k} \triangleq  \tfrac{n}{2\eta^2} \left( \tilde{f}_i(x_{i,k}+v_{T_r},y(\hat x_r+v_{T_r}),\xi_{i,k})- \tilde{f}_i(x_{i,k}\yqro{-v_{T_r}},y(\hat x_r\yqro{-v_{T_r}}),\xi_{i,k})\right)v_{T_r},\\
&g^{\eta,\varepsilon_r}_{i,k} \triangleq \tfrac{n}{2\eta^2} \left( \tilde{f}_i(x_{i,k}+v_{T_r},y_{\varepsilon_r}(\hat x_r+v_{T_r}),\xi_{i,k})- \tilde{f}_i(x_{i,k}\yqro{-v_{T_r}},y_{\varepsilon_r}(\hat x_r\yqro{-v_{T_r}}),\xi_{i,k})\right)v_{T_r},\\
&{\omega_{i,k}^\eta}  \triangleq {\hat g}^{\eta}_{i,k} - g^{\eta}_{i,k} , \qquad  {w_{i,k}^\eta}  \triangleq g^{\eta,\varepsilon_r}_{i,k} - {\hat g}^{\eta}_{i,k}, \qquad \fy{\bar{g}^{\eta}_{k} \triangleq \tfrac{1}{m}\textstyle\sum_{i=1}^m g^{\eta}_{i,k}, } \qquad  \bar{g}^{\eta,\varepsilon_r}_{k} \triangleq \tfrac{1}{m}\textstyle\sum_{i=1}^m g^{\eta,\varepsilon_r}_{i,k} \\
& {\bar{\omega}^\eta}_k \triangleq \tfrac{1}{m}\textstyle\sum_{i=1}^m {\omega_{i,k}^\eta}, \qquad  {\bar{w}^\eta}_k \triangleq \tfrac{1}{m}\textstyle\sum_{i=1}^m {w_{i,k}^\eta}.
\end{align*}

\end{definition}

\begin{remark} 
In view of Definition~\ref{def:s_terms_bilevel}, we have $g^{\eta,\varepsilon_r}_{i,k} =g^{\eta}_{i,k}+{\omega_{i,k}^\eta} +{w_{i,k}^\eta}$.  The term $g^{\eta}_{i,k}$ denotes a zeroth-order stochastic gradient of the local implicit objective of {client} $i$ at iteration $k$,  {$\hat g^{\eta}_{i,k}$} denotes a variant of $g^{\eta}_{i,k}$ where $y(\bullet)$ are obtained at delayed updates,  and $g^{\eta,\varepsilon_r}_{i,k}$ denotes the inexact variant of $g^{\eta}_{i,k}$ where $y(\bullet)$ is only inexactly evaluated with prescribed accuracy $\varepsilon_r$.  While in Algorithm~\ref{alg:FedRiZO_upper}, only  $g^{\eta,\varepsilon_r}_{i,k}$ is employed at the local steps,  we utilize the equations $g^{\eta,\varepsilon_r}_{i,k} =g^{\eta}_{i,k}+{\omega_{i,k}^\eta} +{w_{i,k}^\eta}$ and $\bar{g}^{\eta,\varepsilon_r}_{k} =\bar{g}^{\eta}_{k}+{\bar{\omega}^\eta}_{k} +{\bar{w}^\eta}_{k}$ to analyze the method.  
\end{remark}

We define the history of Algorithm~\ref{alg:FedRiZO_upper}, for $T_r \leq k \leq T_{r+1}-1$ and $r\geq 1$ as  
\begin{align*}
\mathcal{F}_k \triangleq\left(\cup_{i=1}^{m}\cup_{t=0}^{k-1} \{ \xi_{i, t}\}\right) \cup \left(\cup_{j=0}^{r-1} \{v_{T_j}\}\right) \cup \left(  \cup_{j=0}^{r}\mathcal{F}^2_{j}\right)\cup\{\fy{\hat{x}_0}\}, 
\end{align*}
and for $1 \leq k \leq T_{1}-1$ as $
\mathcal{F}_k \triangleq\left(\cup_{i=1}^{m}\cup_{t=0}^{k-1} \{ \xi_{i, t}\}\right) \fy{\cup  \{v_{T_0}\}} \cup  \mathcal{F}^2_{0}\cup\{\fy{\hat{x}_0}\}, $
and $\mathcal{F}_0 \triangleq  \{\fy{\hat{x}_0}\}$. Here, $\mathcal{F}_j^2$ denotes the collection of all random variables generated in the two calls to {the lower-level FL method (e.g., FedAvg)} during the $j$th round of Algorithm~\ref{alg:FedRiZO_upper}.

In the following, we analyze the error terms in Definition~\ref{def:s_terms_bilevel}. 
\begin{lemma}\label{lemma:bounds_bilevel_6terms}\em
Consider Algorithm~\ref{alg:FedRiZO_upper} \fy{where $\{\varepsilon_r\}$ is such that $\mathbb{E}[\|y_{\varepsilon_{r}}(\bullet)-y(\bullet)\|^2 ]\leq \varepsilon_{r}$ for all $r\geq 0$. } Let Assumption~\ref{assum:bilevel} hold.  Then, the following statements hold {in an almost-sure sense} for all $k\geq 0$ and $i \in [m]$.
\begin{flalign*}
&\mbox{(i)}\ \   \mathbb{E}\left[g^{\eta}_{i,k} \mid \mathcal{F}_k  \right]  = \nabla f_i^{\eta}(x_{i,k}).&&\\
&\mbox{(ii)}\ \   \mathbb{E}\left[\|g^{\eta}_{i,k}\|^2  \right] \leq \yqro{16\sqrt{2\pi}(L_0^{\text{imp}})^2n.}&&\\
&\mbox{(iii)}\ \ \mathbb{E}\left[ \|{\omega_{i,k}^\eta}\|^2 \right]\leq \yqro{\tfrac{n^2}{\eta^2}} \left( \tfrac{ L_{0,y}^{f}L_{0,x}^{\nabla h}}{\mu_h}\right)^2\mathbb{E}\left[\|x_{i,k}-\hat x_r\|^2\right]. &&\\
&\mbox{(iv)}\ \  \mathbb{E}\left[ \|{w_{i,k}^\eta}\|^2\right] \leq \yqro{\tfrac{n^2}{\eta^2}} \left(L_{0,y}^{f}\right)^2\varepsilon_r. &&\\
&\mbox{(v)}\ \   \mathbb{E}\left[\|\bar g_k^{\eta}+\bar{\nabla}_{k}^{\eta,d}\|^2 \right] {\leq} \tfrac{\yqro{16\sqrt{2\pi}(L_0^{\text{imp}})^2n}}{m} +2\mathbb{E}\left[\left\|\nabla {\bf f}^\eta (\bar x_{k})\right\|^2\right]   + \yqro{\tfrac{2(L_0^{\text{imp}}\sqrt{n}+1)^2}{\eta^2}} \mathbb{E}\left[\bar e_k\right].&&\\
&\mbox{(vi)}\ \ \mathbb{E}\left[ \nabla {\bf f}^\eta(\bar x_k)^{\fyy{\top}}(\bar{g}_k^{\eta}+\bar{\nabla}_{k}^{\eta,d}) \right] \geq \tfrac{1}{2}\mathbb{E}\left[\|\nabla {\bf f}^\eta(\bar x_k)\|^2\right] - \yqro{\tfrac{(L_0^{\text{imp}}\sqrt{n}+1)^2}{2\eta^2}} \mathbb{E}\left[ \bar e_k\right].&&\\
&\mbox{(vii)}\ \ \mathbb{E}\left[ \|{\bar{\omega}^\eta}_{k}\|^2 \right] \leq  \yqro{\tfrac{n^2}{\eta^2}}\left( \tfrac{L_{0,x}^{\nabla h}}{\mu_h}\right)^2 (L_{0,y}^f)^2\ \mathbb{E}[\hat{e}_k ] .&&\\
&\mbox{(viii)}\ \ \mathbb{E}\left[ \|{\bar{w}^\eta}_{k}\|^2 \right] \leq \yqro{\tfrac{n^2}{\eta^2}}(L_{0,y}^f)^2 \varepsilon_r  .&&
\end{flalign*}
\end{lemma} 
\proof{Proof.}
\noindent (i) From Definition~\ref{def:s_terms_bilevel}, we can write
\begin{align*}
&\mathbb{E}\left[g^{\eta}_{i,k} \mid \mathcal{F}_k\right]\\  &= \mathbb{E}_{v_{T_r}}\left[\mathbb{E}_{\xi_{i, k}}\left[\tfrac{n}{\yqro{2}\eta^2} \left({\tilde{f}_i}(x_{i,k}+v_{T_r},y(x_{i,k}+v_{T_r}),\xi_{i,k}) - {\tilde{f}_i}(x_{i,k}\yqro{-v_{T_r}},y(x_{i,k}\yqro{-v_{T_r}}),\xi_{i,k})\right)v_{T_r} \right]\right]\\
& =\mathbb{E}_{v_{T_r}}\left[\tfrac{n}{\yqro{2}\eta^2}  f_i(x_{i,k}+v_{T_r},y(x_{i,k}+v_{T_r}))   v_{T_r} \mid \mathcal{F}_k\right] \yqro{+\mathbb{E}_{v_{T_r}}\left[\tfrac{n}{2\eta^2}  f_i(x_{i,k}-v_{T_r},y(x_{i,k}-v_{T_r})) (-v_{T_r}) \mid \mathcal{F}_k\right]} \\
&\stackrel{\tiny \mbox{Lemma~\ref{SphericalSmooth} (i)}}{=} \nabla f_{i}^\eta(x_{i, k}).
\end{align*}

{
\noindent (ii) From the definition of $g^{\eta}_{i,k}$ and that $\|v_{T_r}\| =\eta$, we have
\begin{align*}
&\mathbb{E}_{v_{T_r}}\left[\left\|g^{\eta}_{i,k}\right\|^2 \mid \mathcal{F}_k {\cup \{ \xi_{i, k}\}} \right] \\
&= \mathbb{E}_{v_{T_r}}\left[\left\|\tfrac{n}{\yqro{2}\eta^2} \left({\tilde{f}_i}(x_{i,k}+v_{T_r},y(x_{i,k}+v_{T_r}),\xi_{i,k}) - {\tilde{f}_i}(x_{i,k}\yqro{-v_{T_r}},y(x_{i,k}\yqro{-v_{T_r}}),\xi_{i,k})\right)v_{T_r}\right\|^2 \right].
\end{align*}
\yqro{From Lemma \ref{lem:implicit_props} (ii), we know that the random implicit function ${\tilde{f}_i}(\bullet,y(\bullet),\xi_{i})$ is $L_0^{\text{imp}}(\xi_i)$-Lipschitz. Then by adding and subtracting $\mathbb{E}_{\hat v_r}[{\tilde{f}_i}(x_{i,k}+\hat v_{r},y(x_{i,k}+\hat v_{r}),\xi_{i,k})]$, where $\hat v_r$ is uniformly distributed on $\eta \mathbb S$, we can utilize Lemma \ref{lem:Levy} by following the same approach in Lemma \ref{hetero:local:gradient} (ii) and obtain
\begin{align*}
\mathbb{E}_{v_{T_r}}\left[\left\|g^{\eta}_{i,k}\right\|^2 \right] &\leq 16\sqrt{2\pi}(L_0^{\text{imp}}(\xi_i))^2n.
\end{align*}
}

Taking expectations with respect to $\xi_{i,k}$ on both sides and invoking $L_0^{\text{imp}} := \max_{i=1,\ldots,m}\sqrt{\mathbb{E}[(L_0^{\text{imp}}(\xi_{i}))^2]}$, we get $\mathbb{E}\left[\|g^{\eta}_{i,k}\|^2 \mid \mathcal{F}_k\right] \leq \yqro{16\sqrt{2\pi}(L_0^{\text{imp}})^2n}$. Then, taking  expectations on both sides of the preceding inequality, we obtain $\mathbb{E}\left[\|g^{\eta}_{i,k}\|^2 \right] \leq \yqro{16\sqrt{2\pi}(L_0^{\text{imp}})^2n}$.
}

\noindent (iii) Consider the definition of ${\omega_{i,k}^\eta}$. We can write
\begin{align}
&\mathbb{E}\left[ \|{\omega_{i,k}^\eta}\|^2\mid \mathcal{F}_k {\cup \{ \xi_{i, k}\}}\right] = \tfrac{n^2}{\yqro{4}\eta^4} \mathbb{E}\left[ \|\left( \tilde{f}_i(x_{i,k}+v_{T_r},y(\hat x_r+v_{T_r}),\xi_{i,k})- \tilde{f}_i(x_{i,k}\yqro{-v_{T_r}},y(\hat x_r\yqro{-v_{T_r}}),\xi_{i,k}) \right.\right.\notag\\
&\left.\left.-\tilde{f}_i(x_{i,k}+v_{T_r},y(x_{i,k}+v_{T_r}),\xi_{i,k}) + \tilde{f}_i(x_{i,k}\yqro{-v_{T_r}},y(x_{i,k}\yqro{-v_{T_r}}),\xi_{i,k})\right)v_{T_r} \|^2\mid \mathcal{F}_k {\cup \{ \xi_{i, k}\}}\right]\notag\\
&\leq \tfrac{n^2}{\yqro{2}\eta^4} \mathbb{E}\left[ \| (\tilde{f}_i(x_{i,k}+v_{T_r},y(\hat x_r+v_{T_r}),\xi_{i,k})-\tilde{f}_i(x_{i,k}+v_{T_r},y(x_{i,k}+v_{T_r}),\xi_{i,k}))v_{T_r}  \|^2 \mid \mathcal{F}_k {\cup \{ \xi_{i, k}\}}\right]\notag\\
&+\tfrac{n^2}{\yqro{2}\eta^4} \mathbb{E}\left[  \| (-\tilde{f}_i(x_{i,k}\yqro{-v_{T_r}},y(\hat x_r\yqro{-v_{T_r}}),\xi_{i,k})+ \tilde{f}_i(x_{i,k}\yqro{-v_{T_r}},y(x_{i,k}\yqro{-v_{T_r}}),\xi_{i,k}))v_{T_r} \|^2\mid \mathcal{F}_k {\cup \{ \xi_{i, k}\}}\right].\label{ineq:lemma9_split}
\end{align}
Now consider the second term in the preceding relation.  \fy{Invoking Assumption~\ref{assum:bilevel} (i)} we have
\begin{align*}
&\mathbb{E}\left[\| (-\tilde{f}_i(x_{i,k}\yqro{-v_{T_r}},y(\hat x_r\yqro{-v_{T_r}}),\xi_{i,k})+ \tilde{f}_i(x_{i,k}\yqro{-v_{T_r}},y(x_{i,k}\yqro{-v_{T_r}}),\xi_{i,k}))v_{T_r} \|^2\mid \mathcal{F}_k {\cup \{ \xi_{i, k}\}}\right]\\
& {\leq}(L_{0,y}^{f}(\xi_{i,k}))^2\int_{\eta\mathbb{S}} \left\|y(x_{i,k}\yqro{-v_{T_r}})-y(\hat x_r\yqro{-v_{T_r}})\right\|^2\|v_{T_r}\|^2p_v(v_{T_r}) \,dv_{T_r} \\
& \stackrel{\text{Lemma } \ref{lem:implicit_props} \text{ (i)}}{\leq} \left( \tfrac{\eta L_{0,y}^{f}(\xi_{i,k})L_{0,x}^{\nabla h}}{\mu_h}\right)^2\|x_{i,k}-\hat x_r\|^2.
\end{align*}
We can obtain a similar bound on the first term in \eqref{ineq:lemma9_split}. We obtain 
\begin{align*}
\mathbb{E}\left[ \|{\omega_{i,k}^\eta}\|^2\mid \mathcal{F}_k{\cup \{ \xi_{i, k}\}}\right]  &\leq \tfrac{n^2}{\yqro{2}\eta^2} \left( \tfrac{ L_{0,y}^{f}(\xi_{i,k})L_{0,x}^{\nabla h}}{\mu_h}\right)^2\|x_{i,k}-\hat x_r\|^2+\tfrac{n^2}{\yqro{2}\eta^2} \left( \tfrac{ L_{0,y}^{f}(\xi_{i,k})L_{0,x}^{\nabla h}}{\mu_h}\right)^2\|x_{i,k}-\hat x_r\|^2\\
&=\yqro{\tfrac{n^2}{\eta^2}} \left( \tfrac{ L_{0,y}^{f}(\xi_{i,k})L_{0,x}^{\nabla h}}{\mu_h}\right)^2\|x_{i,k}-\hat x_r\|^2.
\end{align*}
Taking expectations with respect to $\xi_{i,k}$ on the both sides, we get $\mathbb{E}\left[\|{\omega_{i,k}^\eta}\|^2 \mid \mathcal{F}_k\right] \leq \yqro{\tfrac{n^2}{\eta^2}} \left( \tfrac{ L_{0,y}^{f}L_{0,x}^{\nabla h}}{\mu_h}\right)^2\mathbb{E}_{\xi_{i,k}}\left[\|x_{i,k}-\hat x_r\|^2\right]$.  Then, taking expectations on the both sides of the preceding relation, we obtain $\mathbb{E}\left[\|{\omega_{i,k}^\eta}\|^2\right] \leq \yqro{\tfrac{n^2}{\eta^2}} \left( \tfrac{ L_{0,y}^{f}L_{0,x}^{\nabla h}}{\mu_h}\right)^2\mathbb{E}\left[\|x_{i,k}-\hat x_r\|^2\right]$.

{
\noindent (iv) Consider the definition of ${w_{i,k}^\eta}$.  Then we may bound $\mathbb{E}\left[ \|{w_{i,k}^\eta}\|^2\right]$ as follows.
\begin{align*}
&\mathbb{E}\left[ \|{w_{i,k}^\eta}\|^2\mid \mathcal{F}_k{\cup \{ \xi_{i, k}\}}\right] =  \tfrac{n^2}{\yqro{4}\eta^4} \mathbb{E}\left[ \|\left( \tilde{f}_i(x_{i,k}+v_{T_r},y_{\varepsilon_r}(\hat x_r+v_{T_r}),\xi_{i,k})- \tilde{f}_i(x_{i,k}\yqro{-v_{T_r}},y_{\varepsilon_r}(\hat x_r\yqro{-v_{T_r}}),\xi_{i,k}) \right.\right.\\
&\left.\left. -\tilde{f}_i(x_{i,k}+v_{T_r},y(\hat x_r+v_{T_r}),\xi_{i,k})+ \tilde{f}_i(x_{i,k}\yqro{-v_{T_r}},y(\hat x_r\yqro{-v_{T_r}}),\xi_{i,k})\right)v_{T_r} \|^2\mid \mathcal{F}_k{\cup \{ \xi_{i, k}\}}\right]\\
&\leq \tfrac{n^2}{\yqro{2}\eta^4} \mathbb{E}\left[ \| (\tilde{f}_i(x_{i,k}+v_{T_r},y_{\varepsilon_r}(\hat x_r+v_{T_r}),\xi_{i,k})-\tilde{f}_i(x_{i,k}+v_{T_r},y(\hat x_r+v_{T_r}),\xi_{i,k}))v_{T_r} \|^2 \mid \mathcal{F}_k{\cup \{ \xi_{i, k}\}} \right]  \\
&+ \tfrac{n^2}{\yqro{2}\eta^4} \mathbb{E}\left[\| (- \tilde{f}_i(x_{i,k}\yqro{-v_{T_r}},y_{\varepsilon_r}(\hat x_r\yqro{-v_{T_r}}),\xi_{i,k})+ \tilde{f}_i(x_{i,k}\yqro{-v_{T_r}},y(\hat x_r\yqro{-v_{T_r}}),\xi_{i,k}))v_{T_r}\|^2\mid \mathcal{F}_k{\cup \{ \xi_{i, k}\}} \right]\\
&\leq \yqro{\tfrac{n^2}{\eta^2}} \left(L_{0,y}^{f}(\xi_{i,k})\right)^2\varepsilon_r.
\end{align*}
The last inequality is obtained by following steps similar to those in (iii) and by invoking $\mathbb{E}\left[\|{y}_{\varepsilon_r}(x) -y(x) \|^2 \mid x\right] \leq \varepsilon_r$. Next, we take  expectations with respect to $\xi_{i,k}$ on both sides of the preceding relation. Then, we take expectations w.r.t. $\mathcal{F}_k$ on both sides of the resulting inequality and obtain the desired bound.

}

\noindent (v) Using the definition of $\bar g^\eta_{k}$ in Definition~\ref{def:s_terms_bilevel}, we can write
\begin{align}                                                                                     
&\mathbb{E}\left[\|\bar g_k^{\eta}+\bar{\nabla}_{k}^{\eta,d}\|^2\mid  \mathcal{F}_k  \right] = \mathbb{E}\left[\left\|\tfrac{1}{m}\textstyle\sum_{i=1}^{m}(g_{i, k}^{\eta}+\nabla_{i,k}^{\eta,d})\right\|^2\mid  \mathcal{F}_k \right] \nonumber\\
&{=} 
\mathbb{E}\left[\left\|\tfrac{1}{m}\textstyle\sum_{i=1}^{m}\left(g_{i, k}^{\eta}+\nabla_{i,k}^{\eta,d}-\nabla {\bf f}_i^\eta (x_{i,k})+\nabla {\bf f}_i^\eta (x_{i,k}) \right)\right\|^2\mid  \mathcal{F}_k \right] \nonumber \\
&{=}\mathbb{E} \left[\left\|\tfrac{1}{m}\textstyle\sum_{i=1}^{m}\left(g_{i, k}^{\eta}+\nabla_{i,k}^{\eta,d}-(\nabla {f}_i^\eta (x_{i,k})+\nabla_{i,k}^{\eta,d}) \right)\right\|^2 \mid  \mathcal{F}_k\right]+ \left\|\tfrac{1}{m}\textstyle\sum_{i=1}^{m}\nabla {\bf f}_i^\eta (x_{i,k})\right\|^2 \nonumber \\
&\stackrel{\tiny \mbox{Lemma\mbox{ }\ref{lemma:bounds_bilevel_6terms} {(i)}}}{=} \tfrac{1}{m^2}\sum_{i=1}^{m}\mathbb{E}\left[\left\|g_{i, k}^{\eta}  -\nabla {f}_i^\eta (x_{i,k}) \right\|^2 \mid  \mathcal{F}_k \right] + \left\|\nabla {\bf f}^\eta (\bar x_{k})+\tfrac{1}{m}\textstyle\sum_{i=1}^{m}\left(\nabla {\bf f}_i^\eta (x_{i,k}) -\nabla {\bf f}_i^\eta (\bar x_k)\right) \right\|^2  \nonumber\\
&{\leq}  \tfrac{1}{m^2}\sum_{i=1}^{m}\mathbb{E}\left[\left\|g_{i, k}^{\eta}  \right\|^2 \mid  \mathcal{F}_k \right] +2\left\|\nabla {\bf f}^\eta (\bar x_{k})\right\|^2  +2 \left\|\tfrac{1}{m}\textstyle\sum_{i=1}^{m}\left(\nabla {\bf f}_i^\eta (x_{i,k}) -\nabla {\bf f}_i^\eta (\bar x_k)\right)\right\|^2. \nonumber\\
&\stackrel{\tiny \mbox{Lemma\mbox{ }\ref{lemma:bounds_bilevel_6terms} {(ii)}}}{\leq}  \tfrac{\yqro{16\sqrt{2\pi}(L_0^{\text{imp}})^2n}}{m} +2\left\|\nabla {\bf f}^\eta (\bar x_{k})\right\|^2  +2 \left\|\tfrac{1}{m}\textstyle\sum_{i=1}^{m}\left(\nabla {\bf f}_i^\eta (x_{i,k}) -\nabla {\bf f}_i^\eta (\bar x_k)\right)\right\|^2. \nonumber
\end{align}
Note that given vectors $y_i \in \mathbb{R}^n$ for $i \in [m]$, we have $\| \tfrac{1}{m}\textstyle\sum_{i=1}^{m}y_i\|^2\leq \tfrac{1}{m}\textstyle\sum_{i=1}^{m}\|y_i\|^2$. Utilizing this inequality, Lemma~\ref{SphericalSmooth} (iv), and that $d_i^{\eta}(\bullet)$ is $\tfrac{1}{\eta}$-smooth, we obtain
\begin{align*}  
&\mathbb{E}\left[\|\bar g_k^{\eta}+\bar{\nabla}_{k}^{\eta,d}\|^2\mid  \mathcal{F}_k  \right]{\leq} \tfrac{\yqro{16\sqrt{2\pi}(L_0^{\text{imp}})^2n}}{m} +2\left\|\nabla {\bf f}^\eta (\bar x_{k})\right\|^2    + \yqro{\tfrac{2(L_0^{\text{imp}}\sqrt{n}+1)^2}{\eta^2m}} \textstyle\sum_{i=1}^{m}\left\|x_{i,k} -\bar x_k\right\|^2.
\end{align*}
Recalling the definition of $\bar e_k$ and taking expectations on both sides, we obtain the bound in (i).

\noindent (vi) We can write
\begin{align}
&\mathbb{E}\left[ \nabla {\bf f}^\eta(\bar x_k)^{\fyy{\top}}(\bar{g}_k^{\eta}+\bar{\nabla}_{k}^{\eta,d})\mid \mathcal{F}_k \right] =\nabla {\bf f}^\eta(\bar x_k)^{\fyy{\top}} \mathbb{E}\left[ \tfrac{1}{m}\textstyle\sum_{i=1}^{m}(g_{i, k}^{\eta}+\bar{\nabla}_{i,k}^{\eta,d})\mid \mathcal{F}_k \right] \nonumber \\
&\stackrel{\tiny \mbox{Lemma\mbox{ }\ref{lemma:bounds_bilevel_6terms} {(i)}}}{=}
\nabla {\bf f}^\eta(\bar x_k)^{\fyy{\top}} \left(\tfrac{1}{m}\textstyle\sum_{i=1}^{m}\nabla {\bf f}_i^\eta (x_{i, k})  \right)  \nonumber \\
&=  \nabla {\bf f}^\eta(\bar x_k)^{\fyy{\top}} \tfrac{1}{m}\textstyle\sum_{i=1}^{m}\left(\nabla {\bf f}_i^\eta (x_{i, k}) -\nabla {\bf f}_i^\eta (\bar x_{k}) +\nabla {\bf f}_i^\eta (\bar x_{k}) \right)  \nonumber \\
&= \nabla {\bf f}^\eta(\bar x_k)^{\fyy{\top}} \tfrac{1}{m}\textstyle\sum_{i=1}^{m}\left(\nabla {\bf f}_i^\eta (x_{i, k})-\nabla {\bf f}_i^\eta (\bar x_{k}) \right)  +\|\nabla {\bf f}^\eta(\bar x_{k})\|^2 \nonumber \\
& {\geq}-\tfrac{1}{2}\|\nabla {\bf f}^\eta(\bar x_k)\|^2 - \tfrac{1}{2} \left\|\tfrac{1}{m}\textstyle\sum_{i=1}^{m}\left(\nabla {\bf f}_i^\eta (x_{i, k})-\nabla {\bf f}_i^\eta (\bar x_{k})\right) \right\|^2  +\|\nabla {\bf f}^\eta(\bar x_{k})\|^2  \nonumber \\
& {\geq}\tfrac{1}{2}\|\nabla {\bf f}^\eta(\bar x_k)\|^2 - \tfrac{1}{2m} \textstyle\sum_{i=1}^{m}\left\| \nabla {\bf f}_i^\eta (x_{i, k})-\nabla {\bf f}_i^\eta (\bar x_{k})  \right\|^2  \nonumber \\
&\stackrel{\tiny \mbox{Lemma\ \ref{SphericalSmooth} (iv)}}{\geq}\tfrac{1}{2}\|\nabla {\bf f}^\eta(\bar x_k)\|^2 - \yqro{\tfrac{(L_0^{\text{imp}}\sqrt{n}+1)^2}{2m\eta^2}}\textstyle\sum_{i=1}^{m} \|x_{i, k}-\bar x_{k} \|^2  \nonumber .
\end{align}
The bound is obtained by recalling the definition of $\bar e_k$ and taking expectations on both sides.

\noindent (vii) We have
\begin{align*}
\mathbb{E}\left[ \|{\bar{\omega}^\eta}_{k}\|^2 \right] &\leq  \tfrac{1}{m^2}\mathbb{E}\left[ \|\textstyle\sum_{i=1}^m {\omega}_{i,k}\|^2 \right] \leq \tfrac{1}{m }\mathbb{E}\left[ \textstyle\sum_{i=1}^m \|{\omega}_{i,k}\|^2 \right]=\tfrac{1}{m }\textstyle\sum_{i=1}^m \mathbb{E}\left[ \|{\omega}_{i,k}\|^2 \right]\\
& \stackrel{\tiny\mbox{(iii)}}{\leq} \yqro{\tfrac{n^2}{\eta^2}} \left( \tfrac{ L_{0,y}^{f}L_{0,x}^{\nabla h}}{\mu_h}\right)^2\textstyle\sum_{i=1}^m\frac{\mathbb{E}[\|x_{i,k}-\hat{x}_r\|^2]}{m}.
\end{align*}

\noindent (viii) We can write
\begin{align*}
\mathbb{E}\left[ \|{\bar{w}^\eta}_{k}\|^2 \right] &\leq  \tfrac{1}{m^2}\mathbb{E}\left[ \|\textstyle\sum_{i=1}^m {w_{i,k}^\eta}\|^2 \right] \leq  \tfrac{1}{m}\mathbb{E}\left[ \textstyle\sum_{i=1}^m \|{w_{i,k}^\eta}\|^2 \right]=\tfrac{1}{m} \textstyle\sum_{i=1}^m \mathbb{E}\left[\|{w_{i,k}^\eta}\|^2 \right]  \stackrel{\tiny\mbox{(iv)}}{\leq}  \yqro{\tfrac{n^2}{\eta^2}} \left(L_{0,y}^{f}\right)^2\varepsilon_r.
\end{align*}
Next, we derive an upper bound on the average delay and average consensus violation, after performing $k$ local steps by a {client} in Algorithm~\ref{alg:FedRiZO_upper}.  We make use of the following result.

\begin{lemma}\label{lem:aggr_a_k_bilevel}\em
Suppose for $T_r+1\leq k \leq T_{r+1}$, for any $r\geq 0$, where $T_0:=0$, the nonnegative sequences $\{a_k\}$, $\{b_k\}$, and $\{\theta_k\}$ satisfy   
 $\max\{a_{k},b_k\}\leq (k-T_r)\gamma^2 \sum_{t=T_r}^{k-1}(\beta_1 a_t +\beta_2 b_t+\theta_t),$ 
 where  $a_{T_r}=b_{T_r}=0$ and $\beta_1,\beta_2 >0$ . If $0<\gamma \leq \frac{1}{\sqrt{2\max\{\beta_1,\beta_2\}}H}$, then for $T_r+1\leq k \leq T_{r+1}$, for any $r\geq 0$,
 $a_k+b_k\leq  {6}H\gamma^2 \sum_{t=T_r}^{k-1}\theta_t.$ 

{
\proof{Proof.}
From the given recursive relation, we have for any $k\geq 1$,  $a_{k}\leq (k-T_r)\gamma^2 \sum_{t=T_r}^{k-1}(\beta_1 a_t +\beta_2 b_t+\theta_t)$ 
and
 $b_k\leq (k-T_r)\gamma^2 \sum_{t=T_r}^{k-1}(\beta_1 a_t +\beta_2 b_t+\theta_t).$ 
Summing the preceding inequalities, we have
\begin{align*}
a_{k}+b_k &\leq 2(k-T_r)\gamma^2 \sum_{t=T_r}^{k-1}(\beta_1 a_t +\beta_2 b_t+\theta_t) \leq (k-T_r)\gamma^2 \sum_{t=T_r}^{k-1}(2\max\{\beta_1,\beta_2\} (a_t + b_t)+2\theta_t).
\end{align*}
If $0<\gamma \leq \frac{1}{{\sqrt{2\max\{\beta_1,\beta_2\}}}H}$, by invoking Lemma \ref{lem:aggr_a_k_extension}, from the preceding relation we have
\begin{align*}
a_{k}+b_k &\leq 3H\gamma^2 \sum_{t=T_r}^{k-1}2\theta_t=6H\gamma^2 \sum_{t=T_r}^{k-1}\theta_t.
\end{align*}
}
\end{lemma}

\begin{lemma}[Bounds on average delay and average consensus violation]\label{lem:delay_cons_bound}  \em                                                         
Let Assumption~\ref{assum:bilevel} hold. Consider Algorithm~\ref{alg:FedRiZO_upper}.  The following holds.

\noindent (i)  [Recursive bound] For any communication round $r>0$, for all $T_r \leq k \leq T_{r+1}-1$, we have
 \begin{align*}
\max \left\{\mathbb{E}\left[  \bar e_k  \right] , \mathbb{E}\left[  \hat e_k  \right]  \right \}
& \leq  4\gamma^2(k-T_r) \textstyle\sum_{t=T_r}^{k-1} \left( 8\mathbb{E}\left[\bar{e}_t\right]+ {\left( \tfrac{L_{0,x}^{\nabla h}}{\mu_h}\right)^2\yqro{\tfrac{n^2}{\eta^2}}(L_{0,y}^f)^2} \mathbb{E}[\hat{e}_k]\right.\\
&\left. + 4B_2^2 \mathbb{E}\left[\left\|\nabla {\bf f}^{\eta}(\bar{x}_t)\right\|^2\right]+\yqro{16\sqrt{2\pi}(L_0^{\text{imp}})^2n} +\tfrac{2B_1^2}{\eta^2}+ \yqro{64B_2^2\sqrt{2\pi}(L_0^{\text{imp}})^2n}  +\yqro{\tfrac{n^2}{\eta^2}}(L_{0,y}^f)^2\varepsilon_r\right).
\end{align*}

\noindent (ii)  [Non-recursive bound] Let {$0<\gamma \leq \left(\yqro{2}H\max\left\{\yqro{4},{ \left( \tfrac{L_{0,x}^{\nabla h}}{\mu_h}\right) \tfrac{\sqrt{2}n}\eta L_{0,y}^f  }\right\}\right)^{-1}$.} Then, for any $r>0$, for all $T_r \leq k \leq T_{r+1}-1$, we have
 \begin{align*}
\mathbb{E}\left[  \bar e_k  \right]+\mathbb{E}\left[  \hat e_k  \right] & \leq {96B_2^2}H\gamma^2\textstyle\sum_{t=T_r}^{k-1}\mathbb{E}\left[\left\|\nabla {\bf f}^{\eta}(\bar{x}_t)\right\|^2\right]\\
    &+{24}H^2\gamma^2\left(\yqro{(4B^2_2+1)16\sqrt{2\pi}(L_0^{\text{imp}})^2n} +\tfrac{2B_1^2}{\eta^2}+  +\yqro{\tfrac{n^2}{\eta^2}}(L_{0,y}^f)^2\varepsilon_r \right).
\end{align*}
\end{lemma}
\proof{Proof.}  
\noindent (i) This proof is comprised of two steps. In the first step we show that the required bound holds for $\mathbb{E}\left[  \bar e_k  \right] $. Then, in the second step we show that it holds for $\mathbb{E}\left[  \hat e_k  \right] $ as well.

\noindent (Step 1)  In view of Algorithm~\ref{alg:FedRiZO_upper}, for any $i$ at any communication round $r>0$, for all $T_r \leq k \leq T_{r+1}-1$ we have $x_{i,k+1} =  x_{i,k} - \gamma  (g^{\eta,\varepsilon_r}_{i,k}+\nabla^{\eta,d}_{i,k}).$ This implies that $x_{i,k} =  x_{i,k-1} - \gamma  (g^{\eta,\varepsilon_r}_{i,k-1}+\nabla^{\eta,d}_{i,k-1})$ for all $ T_r+1 \leq k \leq T_{r+1}.$ Unrolling the preceding relation recursively, we obtain
\begin{align}\label{eqn:nonrecursive_local_bilevel}
x_{i,k} =  x_{i,T_r} - \gamma  \textstyle\sum_{t=T_r}^{k-1}(g^{\eta,\varepsilon_r}_{i,t}+\nabla^{\eta,d}_{i,t}),\quad \hbox{for all } T_r+1 \leq k \leq T_{r+1}.
\end{align} 
From Algorithm~\ref{alg:FedRiZO_upper}, we have $\hat{x}_r =x_{i,T_r}$. Invoking the definition of $\bar{x}_k$, we have $\bar{x}_{T_r} = \hat{x}_r$. This implies that $\bar{x}_{T_r} =x_{i,T_r}$ for all $i$ and $r$. Extending Lemma~\ref{lem:bar_x_recursion}, we can write
\begin{align}\label{eqn:nonrecursive_ave_bilevel}
\bar x_k =  x_{i,T_r} - \gamma  \textstyle\sum_{t=T_r}^{k-1}(\bar{g}^{\eta,\varepsilon_r}_{t}+\bar{\nabla}^{\eta,d}_{t}),\quad \hbox{for all } T_r+1 \leq k \leq T_{r+1}.
\end{align}
Using \eqref{eqn:nonrecursive_local_bilevel} and \eqref{eqn:nonrecursive_ave_bilevel}, for all $T_r+1 \leq k \leq T_{r+1}$ we have
\begin{align*}
\mathbb{E}\left[  \bar e_k  \right]  
&=\tfrac{1}{m}\textstyle\sum_{i=1}^{m}\mathbb{E}\left[ \|x_{i, k}-\bar x_k \|^2   \right]   \\
&=\tfrac{1}{m}\textstyle\sum_{i=1}^{m}\mathbb{E}\left[ \left\|\gamma\textstyle\sum_{t=T_r}^{k-1} (g_{i, t}^{\eta,\varepsilon_r} +\nabla^{\eta,d}_{i,t})-\gamma\textstyle\sum_{t=T_r}^{k-1} (\bar{g}_{t}^{\eta,\varepsilon_r}+ \bar{\nabla}^{\eta,d}_{t}) \right\|^2  \right]   \\
&\leq \tfrac{\gamma^2(k-T_r)}{m}\textstyle\sum_{i=1}^{m}\textstyle\sum_{t=T_r}^{k-1} \mathbb{E}\left[ \left\| (g_{i, t}^{\eta,\varepsilon_r} +\nabla^{\eta,d}_{i,t})- (\bar{g}_{t}^{\eta,\varepsilon_r}+ \bar{\nabla}^{\eta,d}_{t}) \right\|^2  \right],
\end{align*}
where the preceding relation follows from the inequality $\|  \textstyle\sum_{t=1}^{T}y_t\|^2\leq T\textstyle\sum_{t=1}^{T}\|y_t\|^2$ for any $y_t \in \mathbb{R}^n$ for $t \in [T]$. Consequently,
\begin{align*}
\mathbb{E}\left[  \bar e_k   \right]  
&\leq \tfrac{\gamma^2(k-T_r)}{m}\textstyle\sum_{t=T_r}^{k-1}\textstyle\sum_{i=1}^{m} \mathbb{E}\left[ \left\| (g_{i, t}^{\eta} +\nabla^{\eta,d}_{i,t})- (\bar{g}_{t}^{\eta,\varepsilon_r}+ \bar{\nabla}^{\eta,d}_{t}) \right\|^2   \right]\\
& = \tfrac{\gamma^2(k-T_r)}{m}\textstyle\sum_{t=T_r}^{k-1}\sum_{i=1}^{m}\left(\mathbb{E}\left[ \left\|  g_{i, t}^{\eta,\varepsilon_r} +\nabla^{\eta,d}_{i,t} \right\|^2   \right] +\mathbb{E}\left[ \left\| \bar{g}_{t}^{\eta,\varepsilon_r}+ \bar{\nabla}^{\eta,d}_{t}  \right\|^2   \right] \right)\\
& -2\tfrac{\gamma^2(k-T_r)}{m}\textstyle\sum_{t=T_r}^{k-1}\sum_{i=1}^{m}\mathbb{E}\left[ (  g_{i, t}^{\eta,\varepsilon_r} +\nabla^{\eta,d}_{i,t} )^{\fyy{\top}}(\bar{g}_{t}^{\eta,\varepsilon_r}+ \bar{\nabla}^{\eta,d}_{t}  )   \right],
\end{align*}
Observing that 
 $\tfrac{1}{m}\textstyle\sum_{i=1}^{m}\mathbb{E}\left[ (  g_{i, t}^{\eta,\varepsilon_r} +\nabla^{\eta,d}_{i,t} )^{\fyy{\top}}(\bar{g}_{t}^{\eta,\varepsilon_r}+ \bar{\nabla}^{\eta,d}_{t}  ) \mid \mathcal{F}_{T_r}  \right]= \mathbb{E}\left[\left\|\bar{g}_{t}^{\eta,\varepsilon_r}+ \bar{\nabla}^{\eta,d}_{t}  \right\|^2  \right],$ 
we obtain 
\begin{align*}
\mathbb{E}\left[  \bar e_k  \right]  
&\leq  \tfrac{\gamma^2(k-T_r)}{m}\textstyle\sum_{t=T_r}^{k-1}\sum_{i=1}^{m} \mathbb{E}\left[ \left\|  g_{i, t}^{\eta,\varepsilon_r} +\nabla^{\eta,d}_{i,t} \right\|^2  \right]\\
& \leq \tfrac{4\gamma^2(k-T_r)}{m}\textstyle\sum_{t=T_r}^{k-1}\sum_{i=1}^{m} \left(\mathbb{E}\left[  \|  g_{i, t}^{\eta}\|^2 +\|\nabla^{\eta,d}_{i,t}  \|^2 +\|{\omega_{i,k}^\eta}\|^2+\|{w_{i,k}^\eta}\|^2 \right]\right)\\
&\stackrel{\tiny \mbox{Lemma\mbox{ }\ref{lemma:bounds_bilevel_6terms}}}{\leq} \tfrac{4\gamma^2(k-T_r)}{m}\textstyle\sum_{t=T_r}^{k-1}\sum_{i=1}^{m} \left(\yqro{16\sqrt{2\pi}(L_0^{\text{imp}})^2n}+\|\nabla^{\eta,d}_{i,t}  \|^2\right.\\
&{\left. + \left( \tfrac{L_{0,x}^{\nabla h}}{\mu_h}\right)^2\yqro{\tfrac{n^2}{\eta^2}}(L_{0,y}^f)^2\mathbb{E}[\|x_{i,t}-\hat{x}_r\|^2]+\yqro{\tfrac{n^2}{\eta^2}}(L_{0,y}^f)^2\varepsilon_r\right)}.
\end{align*}
Invoking Lemma~\ref{lem:bound_on_dist2_ave}, the following bound emerges. 
\begin{align*}
\mathbb{E}\left[  \bar e_k  \right]   
& \leq  4\gamma^2(k-T_r) \textstyle\sum_{t=T_r}^{k-1} \left(\yqro{16\sqrt{2\pi}(L_0^{\text{imp}})^2n}+ 8\mathbb{E}\left[\bar{e}_t\right]+\tfrac{2B_1^2}{\eta^2} + 4B_2^2 \mathbb{E}\left[\left\|\nabla {\bf f}^{\eta}(\bar{x}_t)\right\|^2\right]\right.\\
&\left.+ \yqro{64B_2^2\sqrt{2\pi}(L_0^{\text{imp}})^2n}  + {\left( \tfrac{L_{0,x}^{\nabla h}}{\mu_h}\right)^2\yqro{\tfrac{n^2}{\eta^2}}(L_{0,y}^f)^2} \mathbb{E}[\hat{e}_k]+{\yqro{\tfrac{n^2}{\eta^2}}(L_{0,y}^f)^2\varepsilon_r}\right).
\end{align*}

\noindent (Step 2) Consider \eqref{eqn:nonrecursive_local_bilevel}.  From Algorithm~\ref{alg:FedRiZO_upper}, we have $\hat{x}_r =x_{i,T_r}$. Thus, we obtain
\begin{align*}
\mathbb{E}\left[\|x_{i,k}-\hat{x}_r\|^2\right] &\leq \gamma^2 \textstyle\sum_{t=T_r}^{k-1} \mathbb{E}\left[ \left\|  g_{i, t}^{\eta,\varepsilon_r} +\nabla^{\eta,d}_{i,t} \right\|^2  \right].
\end{align*}
Thus, for all $T_r+1 \leq k \leq T_{r+1}$ we have
\begin{align*}
\mathbb{E}\left[  \hat e_k  \right]  
&=\tfrac{1}{m}\textstyle\sum_{i=1}^{m}\mathbb{E}\left[ \|x_{i, k}-\hat x_k \|^2   \right]    = \tfrac{\gamma^2(k-T_r)}{m}\textstyle\sum_{t=T_r}^{k-1}\sum_{i=1}^{m} \mathbb{E}\left[ \left\|  g_{i, t}^{\eta,\varepsilon_r} +\nabla^{\eta,d}_{i,t} \right\|^2  \right].
\end{align*}
Following the steps in (Step 1), we obtain the same bound on $\mathbb{E}\left[  \hat e_k  \right]  $. 

\noindent (ii) The bound in (ii) is obtained by applying Lemma~\ref{lem:aggr_a_k_bilevel} on the inequality in part (i).
\subsection{Proof of Theorem \ref{thm:bilevel}}
\proof{Proof.} (i)
From the $L$-smoothness property of the implicit function, where $L=\yqro{\tfrac{L_0^{\text{imp}}\sqrt{n}+1}{\eta}}$, we have
\begin{align*}
 {\bf f}^{\eta}(\bar x_{k+1})    \leq   {\bf f}^{\eta}(\bar x_{k}) +  \nabla {\bf f}^{\eta}(\bar x_k)^{\fyy{\top}} (\bar x_{k+1}-\bar x_k)   +\yqro{\tfrac{L_0^{\text{imp}}\sqrt{n}+1}{2\eta}} \|\bar x_{k+1}-\bar x_k\|^2.
\end{align*}
In view of the recursion $\bar{x}_{k+1}=\bar{x}_k-\gamma(\bar{g}_k^{\eta,\varepsilon_r}+\bar{\nabla}^{\eta,d}_{k})$, we have
\begin{align*}
 {\bf f}^{\eta}(\bar x_{k+1})    \leq   {\bf f}^{\eta}(\bar x_{k}) -\gamma   \nabla {\bf f}^{\eta}(\bar x_k)^{\fyy{\top}} (\bar{g}_k^{\eta,\varepsilon_r} +\nabla^{\eta,d}_{i,k})  +\yqro{\tfrac{L_0^{\text{imp}}\sqrt{n}+1}{2\eta}}\gamma^2 \|\bar{g}_k^{\eta,\varepsilon_r}+\bar{\nabla}^{\eta,d}_{k}\|^2.
\end{align*}
In view of $\bar{g}^{\eta,\varepsilon_r}_{k} =\bar{g}^{\eta}_{k}+{\bar{\omega}^\eta}_{k} +{\bar{w}^\eta}_{k}$, we have
\begin{align*}
 {\bf f}^{\eta}(\bar x_{k+1})   & \leq   {\bf f}^{\eta}(\bar x_{k}) -\gamma   \nabla {\bf f}^{\eta}(\bar x_k)^{\fyy{\top}}(\bar{g}_k^{\eta}+\bar{\nabla}_{k}^{\eta,d}) -\gamma   \nabla {\bf f}^{\eta}(\bar x_k)^{\fyy{\top}}{\bar{\omega}^\eta}_{k} -\gamma   \nabla {\bf f}^{\eta}(\bar x_k)^{\fyy{\top}}{\bar{w}^\eta}_{k} \\ 
 &  +\yqro{\tfrac{L_0^{\text{imp}}\sqrt{n}+1}{2\eta}}\gamma^2 \|\bar{g}^{\eta}_{k}+\bar{\nabla}_{k}^{\eta,d}+{\bar{\omega}^\eta}_{k} +{\bar{w}^\eta}_{k}\|^2\\
 &   \leq   {\bf f}^{\eta}(\bar x_{k}) -\gamma  \nabla {\bf f}^{\eta}(\bar x_k)^{\fyy{\top}}(\bar{g}_k^{\eta}+\bar{\nabla}_{k}^{\eta,d})+\tfrac{\gamma}{8}\|\nabla {\bf f}^{\eta}(\bar x_k)\|^2+4\gamma\|{\bar{\omega}^\eta}_{k}\|^2+4\gamma\|{\bar{w}^\eta}_{k}\|^2  \\
 & +\yqro{\tfrac{3(L_0^{\text{imp}}\sqrt{n}+1)}{2\eta}}\gamma^2\left( \|\bar{g}^{\eta}_{k}+\bar{\nabla}_{k}^{\eta,d}\|^2+\|{\bar{\omega}^\eta}_{k}\|^2 +\|{\bar{w}^\eta}_{k}\|^2\right),
\end{align*}
where in the second inequality, we utilized the following relation twice: $|u_1^{\fyy{\top}}u_2| \leq 0.5\left(\tfrac{1}{8}\|u_1\|^2 +8\|u_2\|^2\right)$ for any $u_1,u_2 \in \mathbb{R}^n$. Taking expectations on both sides and requiring that $\gamma \leq \yqro{\tfrac{\eta}{24(L_0^{\text{imp}}\sqrt{n}+1)}}$, we obtain 
\begin{align}\label{eqn:thm_bilevel_term0}
 \mathbb{E}\left[{\bf f}^{\eta}(\bar x_{k+1}) \right]      & \leq   \mathbb{E}\left[{\bf f}^{\eta}(\bar x_{k}) \right]-\gamma  \mathbb{E}\left[ \nabla {\bf f}^{\eta}(\bar x_k)^{\fyy{\top}}   (\bar{g}_k^{\eta}+\bar{\nabla}_{k}^{\eta,d})\right]+\tfrac{\gamma}{8}\mathbb{E}\left[\|\nabla {\bf f}^{\eta}(\bar x_k)\|^2\right]\notag\\
&    +\yqro{\tfrac{3(L_0^{\text{imp}}\sqrt{n}+1)}{2\eta}}\gamma^2   \mathbb{E}\left[\|\bar{g}^{\eta}_{k}+\bar{\nabla}_{k}^{\eta,d}\|^2\right]+5\gamma \mathbb{E}\left[\|{\bar{\omega}^\eta}_{k}\|^2 +\|{\bar{w}^\eta}_{k}\|^2\right] .
\end{align}
From Lemma~\ref{lemma:bounds_bilevel_6terms} (vi) and (v), we have
\begin{align}\label{eqn:thm_bilevel_term1}
-\gamma\mathbb{E}\left[ \nabla {\bf f}^\eta(\bar x_k)^{\fyy{\top}}(\bar{g}_k^{\eta}+\bar{\nabla}_{k}^{\eta,d}) \right] \leq  -\tfrac{\gamma}{2}\mathbb{E}\left[\|\nabla {\bf f}^\eta(\bar x_k)\|^2\right] +\gamma\yqro{\tfrac{(L_0^{\text{imp}}\sqrt{n}+1)^2}{2\eta^2}}\mathbb{E}\left[ \bar e_k\right],
\end{align}
and $
\mathbb{E}\left[\|\bar g_k^{\eta}+\bar{\nabla}_{k}^{\eta,d}\|^2 \right] {\leq}\yqro{\tfrac{16\sqrt{2\pi}(L_0^{\text{imp}})^2n}{m}} +2\mathbb{E}\left[\left\|\nabla {\bf f}^\eta (\bar x_{k})\right\|^2\right]   +\yqro{\tfrac{2(L_0^{\text{imp}}\sqrt{n}+1)^2}{\eta^2}} \mathbb{E}\left[\bar e_k\right],$
{respectively.} In view of $\gamma \leq \yqro{\tfrac{\eta}{24(L_0^{\text{imp}}\sqrt{n}+1)}}$, {multiplying both sides of the preceding inequality by $\yqro{\tfrac{3(L_0^{\text{imp}}\sqrt{n}+1)}{2\eta}}\gamma^2$,} we obtain 
\begin{align}\label{eqn:thm_bilevel_term2}
\yqro{\tfrac{3(L_0^{\text{imp}}\sqrt{n}+1)}{2\eta}}\gamma^2\mathbb{E}\left[\|\bar g_k^{\eta}+\bar{\nabla}_{k}^{\eta,d}\|^2   \right]& \leq \yqro{\gamma^2\tfrac{48\sqrt{2\pi}(L_0^{\text{imp}}\sqrt{n}+1)(L_0^{\text{imp}})^2n}{2m\eta}}+\tfrac{\gamma}{8}\mathbb{E}\left[\|\nabla {\bf f}^\eta (\bar x_{k})\|^2\right] +\gamma\yqro{\tfrac{(L_0^{\text{imp}}\sqrt{n}+1)^2}{8\eta^2}}\mathbb{E}\left[\bar e_k\right].
\end{align}
Also,  from Lemma~\ref{lemma:bounds_bilevel_6terms} {(vii) and (viii)}, we have  
\begin{align}\label{eqn:thm_bilevel_term3}
5\gamma\mathbb{E}\left[\|{\bar{\omega}^\eta}_{k}\|^2 +\|{\bar{w}^\eta}_{k}\|^2\right] \leq { \left( \tfrac{L_{0,x}^{\nabla h}}{\mu_h}\right)^2\tfrac{\yqro{5} \gamma n^2}{\eta^2}(L_{0,y}^f)^2}\left( \mathbb{E}[\hat{e}_k ]  +\varepsilon_r \right).
\end{align}
From \eqref{eqn:thm_bilevel_term0}--\eqref{eqn:thm_bilevel_term3}, we obtain 
\begin{align*} 
 \mathbb{E}\left[{\bf f}^{\eta}(\bar x_{k+1}) \right]  & \leq   \mathbb{E}\left[{\bf f}^{\eta}(\bar x_{k}) \right] -\tfrac{\gamma}{4}\mathbb{E}\left[\|\nabla {\bf f}^{\eta}(\bar x_k)\|^2\right] +\yqro{\gamma^2\tfrac{48\sqrt{2\pi}(L_0^{\text{imp}}\sqrt{n}+1)(L_0^{\text{imp}})^2n}{2m\eta}}\\
 & +\yqro{\gamma\tfrac{5(L_0^{\text{imp}}\sqrt{n}+1)^2}{8\eta^2}}\mathbb{E}\left[\bar e_k\right]+ {\left( \tfrac{L_{0,x}^{\nabla h}}{\mu_h}\right)^2\tfrac{\yqro{5} \gamma n^2}{\eta^2}(L_{0,y}^f)^2}\left( \mathbb{E}[\hat{e}_k ]  +\varepsilon_r \right).
\end{align*}

This implies that 
\begin{align*} 
 \mathbb{E}\left[{\bf f}^{\eta}(\bar x_{k+1}) \right]      & \leq   \mathbb{E}\left[{\bf f}^{\eta}(\bar x_{k}) \right] -\tfrac{\gamma}{4}\mathbb{E}\left[\|\nabla {\bf f}^{\eta}(\bar x_k)\|^2\right] +\tfrac{\gamma^2 \Theta_1}{m}    +\gamma\max\{\Theta_2,\Theta_3\}(\mathbb{E}\left[\bar e_k\right]+\mathbb{E}[\hat{e}_k ] )+ \gamma\Theta_3\varepsilon_r .
\end{align*}
where $\Theta_1 :=\yqro{\tfrac{48\sqrt{2\pi}(L_0^{\text{imp}}\sqrt{n}+1)(L_0^{\text{imp}})^2n}{2m\eta}}$, $\Theta_2:=\yqro{\tfrac{5(L_0^{\text{imp}}\sqrt{n}+1)^2}{8\eta^2}}$, and ${\Theta_3:= \left( \tfrac{L_{0,x}^{\nabla h}}{\mu_h}\right)^2\tfrac{\yqro{5}  n^2}{\eta^2}(L_{0,y}^f)^2}$.
From Lemma~\ref{lem:delay_cons_bound}, for {$\gamma \leq \tfrac{\max\left\{\yqro{4}, \sqrt{\yqro{0.4}\Theta_3}\right\} ^{-1}}{\yqro{2}H}$}, we have
\begin{align*} 
\mathbb{E}\left[  \bar e_k  \right]+\mathbb{E}\left[  \hat e_k  \right] & \leq {96B_2^2}H\gamma^2\textstyle\sum_{t=T_r}^{k-1}\mathbb{E}\left[\left\|\nabla {\bf f}^{\eta}(\bar{x}_t)\right\|^2\right]+H^2\gamma^2\Theta_4\varepsilon_r  +H^2\gamma^2\Theta_5,
\end{align*}
where $\Theta_4:=\tfrac{{\yqro{24}}n^2}{\eta^2}(L_{0,y}^f)^2$ and $\Theta_5:= \tfrac{48 B_1^2}{\eta^2}+ \yqro{384\sqrt{2\pi}(4B^2_2+1)(L_0^{\text{imp}})^2n}$.
Since $\gamma \leq \tfrac{\left(\sqrt{48}B_2\max\{\sqrt{\Theta_2},\sqrt{\Theta_3}\}\right)^{-1}}{4H}$,  we have ${96B_2^2}H\gamma^2\max\{\Theta_2,\Theta_3\} \leq  {\tfrac{1}{8H}}$.
Combining the two preceding inequalities, we obtain
\begin{align*} 
 \mathbb{E}\left[{\bf f}^{\eta}(\bar x_{k+1}) \right]      & \leq   \mathbb{E}\left[{\bf f}^{\eta}(\bar x_{k}) \right] -\tfrac{\gamma}{4}\mathbb{E}\left[\|\nabla {\bf f}^{\eta}(\bar x_k)\|^2\right] +{\tfrac{\gamma }{8H}}\textstyle\sum_{t=T_r}^{k-1}\mathbb{E}\left[\left\|\nabla {\bf f}^{\eta}(\bar{x}_t)\right\|^2\right]+\tfrac{\gamma^2 \Theta_1}{m}  \\
 & +\gamma\left(H^2\gamma^2\max\{\Theta_2,\Theta_3\}\Theta_4+\Theta_3\right)\varepsilon_r  +H^2\gamma^3\max\{\Theta_2,\Theta_3\}\Theta_5   .
\end{align*}
Summing the preceding relation from $k:=0,\ldots,K:=T_{R}-1$ for some $R\geq 1$, we have
 \begin{align*} 
 \mathbb{E}\left[{\bf f}^{\eta}(\bar x_{T_R}) \right]      & \leq   \mathbb{E}\left[{\bf f}^{\eta}(\bar x_{0}) \right] -\textstyle\sum_{k=0}^{T_R-1}\tfrac{\gamma}{4}\mathbb{E}\left[\|\nabla {\bf f}^{\eta}(\bar x_k)\|^2\right] +{\tfrac{\gamma}{8H}}\sum_{k=0}^{T_R-1}\textstyle\sum_{t=T_r}^{k-1}\mathbb{E}\left[\left\|\nabla {\bf f}^{\eta}(\bar{x}_t)\right\|^2\right]\\
 &+\tfrac{T_R\gamma^2 \Theta_1}{m}   +\gamma\left(H^2\gamma^2\max\{\Theta_2,\Theta_3\}\Theta_4+\Theta_3\right)\textstyle\sum_{r=0}^{R-1}(T_{r+1}-T_r-1)\varepsilon_r  \\
 &+T_RH^2\gamma^3\max\{\Theta_2,\Theta_3\}\Theta_5   .
\end{align*}
We obtain 
  \begin{align*} 
  \textstyle\sum_{k=0}^{T_R-1}\tfrac{\gamma}{{8}}\mathbb{E}\left[\|\nabla {\bf f}^{\eta}(\bar x_k)\|^2\right]   & \leq   \mathbb{E}\left[{\bf f}^{\eta}(\bar x_{0}) \right] - \mathbb{E}\left[{\bf f}^{\eta}(\bar x_{T_R}) \right]   +\tfrac{T_R\gamma^2 \Theta_1}{m}  \\
 & +\gamma\left(H^2\gamma^2\max\{\Theta_2,\Theta_3\}\Theta_4+\Theta_3\right){H}\textstyle\sum_{r=0}^{R-1}\varepsilon_r    +T_RH^2\gamma^3\max\{\Theta_2,\Theta_3\}\Theta_5   .
\end{align*}
Multiplying both sides by $\tfrac{8}{\gamma T_R}$ and using the definition of $k^*$, we obtain  
\begin{align*} 
 \mathbb{E}\left[\|\nabla {\bf f}^{\eta}(\bar x_{k^*})\|^2\right]   & \leq   8(\gamma T_R)^{-1}(\mathbb{E}\left[{\bf f}^{\eta}(x_{0}) \right] -  {\bf f}^{\eta,*}  )   +\tfrac{8\gamma \Theta_1}{m}  \\
 & + 8\left(H^2\gamma^2\max\{\Theta_2,\Theta_3\}\Theta_4+\Theta_3\right){H}\tfrac{\textstyle\sum_{r=0}^{R-1}\varepsilon_r}{T_R}    + 8H^2\gamma^2\max\{\Theta_2,\Theta_3\}\Theta_5   .
\end{align*}

{
\noindent (ii) Substituting $\gamma:=\yqro{\sqrt{\tfrac{m\eta}{Kn^{1.5}(L_0^{\text{imp}})^3}}}$ and $H:=  \sqrt[4]{\tfrac{K}{m^3}}$ in the error bound in (i), we obtain  
\begin{align}\label{eqn:ineq_proof_thm1_ii}
 \mathbb{E}\left[\|\nabla {\bf f}^{\eta}(\bar x_{k^*})\|^2\right]   & \leq   \tfrac{8(\mathbb{E}\left[{\bf f}^{\eta}(x_{0}) \right]-{\bf f}^{\eta,*})\yqro{\tfrac{n^{0.75}(L_0^{\text{imp}})^{1.5}}{\eta^{0.5}}}}{\sqrt{mK}}   +\tfrac{8 \Theta_1\yqro{\tfrac{\eta^{0.5}}{n^{0.75}(L_0^{\text{imp}})^{1.5}}}}{\sqrt{mK}} \nonumber\\
&+ 8\left(\tfrac{\max\{\Theta_2,\Theta_3\}\Theta_4\yqro{\tfrac{\eta}{n^{1.5}(L_0^{\text{imp}})^{3}}}}{\sqrt{mK}}+\Theta_3\right)\tfrac{\textstyle\sum_{r=0}^{R-1}\varepsilon_r}{(mK)^{\tfrac{3}{4}}}    + 8\tfrac{\max\{\Theta_2,\Theta_3\}\Theta_5\yqro{\tfrac{\eta}{n^{1.5}(L_0^{\text{imp}})^{3}}}}{\sqrt{mK}}   .
\end{align}
Invoking $\varepsilon_r:=\tilde{\mathcal{O}}(\tfrac{1}{m \tilde{T}_{\tilde{R}_r}})$ where $\tilde{T}_{\tilde{R}_r}:=\tilde{\mathcal{O}}\left(m^{-1}(r+1)^{\frac{2}{3}}\yqro{n}\right)$,  we have $\sum_{r=0}^{R-1}\varepsilon_r =\tilde{\mathcal{O}}(\yqro{n^{-1}}R^{\frac{1}{3}}) = \tilde{\mathcal{O}}\left(\yqro{n^{-1}}\tfrac{K^{\frac{1}{3}}}{H^{\frac{1}{3}}}\right)$.  Substituting $H:=  \sqrt[4]{\tfrac{K}{m^3}}$, we obtain $\sum_{r=0}^{R-1}\varepsilon_r = \tilde{\mathcal{O}}\left( \yqro{n^{-1}}(m K)^{\frac{1}{4}}\right)$.  Substituting this bound in \eqref{eqn:ineq_proof_thm1_ii}, \yqro{leads to iteration complexity of $$\mathcal{O}\left(\left(\tfrac{n^{0.75}(L_0^{\text{imp}})^{1.5}}{\eta^{0.5}}+\tfrac{n^{1.5}(L_0^{\text{imp}})^{-3}}{\eta^3}+\tfrac{n}{\eta^2}+\tfrac{\sqrt n(L_0^{\text{imp}})^{-3}}{\eta^3}+\tfrac{n^{1.5}(L_0^{\text{imp}})^{-1}}{\eta} \right)^2\tfrac{1}{m\epsilon^2}\right).$$}


\noindent (iii)  From (ii), we obtain $R=\mathcal{O}(\tfrac{K}{H}) = \mathcal{O}\left(\tfrac{K}{ \sqrt[4]{{K}/{m^3}} }\right) =\mathcal{O}\left((mK)^{3/4}\right)$.

}
 
\subsection{Preliminaries for the proof of Theorem~\ref{thm:fed2smpec}}
In this subsection, we provide preliminary results to prove Theorem~\ref{thm:fed2smpec}.  Some results \us{are similar to} those presented for the analysis of Theorem~\ref{thm:bilevel}. However,  as noted in Remark~\ref{rem:bl_vs_2s},  there are some major algorithmic distinctions, e.g.,  unlike {FedRZO$_{\texttt {bl}}$}, {FedRZO$_{\texttt {2s}}$} does not leverage delays and instead, employs calls to the lower-level deterministic VI solver during the local steps.  Accordingly,  in the following results,  we provide the detailed proofs as needed.  We first introduce some terms in the following.
\begin{definition}\label{def:s_terms_2smpec} Consider the terms given in Definition~\ref{def:basic_terms_bilevel}. Further, consider the following.
\begin{align*}
& {g^{\eta}_{i,k}\triangleq \tfrac{n}{\yqro{2}\eta^2} \left( \tilde{f}_i(x_{i,k}+v_{i,k},y(x_{i,k}+v_{i,k}, \xi_{i,k}),\xi_{i,k})- \tilde{f}_i(x_{i,k}\yqro{-v_{i,k}},y(x_{i,k}\yqro{-v_{i,k}}, \xi_{i,k}),\xi_{i,k})\right)v_{i,k}},\\
& {g^{\eta, \varepsilon_k}_{i,k}\triangleq\tfrac{n}{\yqro{2}\eta^2} \left( \tilde{f}_i(x_{i,k}+v_{i,k},y_{\varepsilon_k}(x_{i,k}+v_{i,k}, \xi_{i,k}),\xi_{i,k})- \tilde{f}_i(x_{i,k}\yqro{-v_{i,k}},y_{\varepsilon_k}(x_{i,k}\yqro{-v_{i,k}}, \xi_{i,k}),\xi_{i,k})\right)v_{i,k}},\\
&  \bar{g}^{\eta,\varepsilon_k}_{k} \triangleq \tfrac{1}{m}\textstyle\sum_{i=1}^m g^{\eta,\varepsilon_k}_{i,k} ,  \qquad {\bar{g}^{\eta}_{k} \triangleq \tfrac{1}{m}\textstyle\sum_{i=1}^m g^{\eta}_{i,k}} ,  \qquad  {{w_{i,k}^\eta}  \triangleq g^{\eta,\varepsilon_k}_{i,k} - {g}^{\eta}_{i,k}},\qquad {\bar{w}^\eta}_k \triangleq \tfrac{1}{m}\textstyle\sum_{i=1}^m w_{i,k}^\eta. 
\end{align*}
\end{definition}

\begin{remark}\label{rem:thm2_bar_terms}
In view of Definition~\ref{def:s_terms_2smpec},   we have $g^{\eta,\varepsilon_k}_{i,k} =g^{\eta}_{i,k} +{w_{i,k}^\eta}$ and $\bar{g}^{\eta,\varepsilon_k}_{k} =\bar{g}^{\eta}_{k} +{\bar{w}^\eta}_{k}$. $\hfill$ $\Box$
\end{remark}
Let us define the history of Algorithm~\ref{alg:fed2:upper} for $k\geq 1$ as  $
\mathcal{F}_k \triangleq {\left(\cup_{i=1}^{m}\cup_{t=0}^{k-1} \{ \xi_{i, t},v_{i,t}\}\right)} \cup\{\hat{x}_0\}, $ where  $\mathcal{F}_0 \triangleq  \{\hat{x}_0\}$. 
\begin{lemma}\label{lemma:bounds_2smpec_6terms}\em
Consider Algorithm~\ref{alg:fed2:upper} where $\{\varepsilon_k\}$ is such that $\|y_{\varepsilon_{k}}(\bullet, \bullet)-y(\bullet, \bullet)\|^2\leq \varepsilon_{k}$ for all $k\geq 0$.   Let Assumption~\ref{assum:fed2smpecs} hold. Then, the following statements hold in an almost-sure sense for all $k\geq 0$ and $i \in [m]$.
\begin{flalign*}
&\mbox{(i)}\ \   \mathbb{E}\left[g^{\eta}_{i,k} \mid \mathcal{F}_k  \right]  = \nabla f_i^{\eta}(x_{i,k}).&&\\
&\mbox{(ii)}\ \   \mathbb{E}\left[\|g^{\eta}_{i,k}\|^2  \right] \leq \yqro{16\sqrt{2\pi}(L_0^{\text{imp}})^2n}.&&\\
&\mbox{(iii)}\ \  \mathbb{E}\left[ \|{w_{i,k}^\eta}\|^2\right] \leq \yqro{\tfrac{n^2}{\eta^2}} \left(L_{0,y}^{f}\right)^2{\varepsilon_k}. &&\\
&\mbox{(iv)}\ \   \mathbb{E}\left[\|\bar g_k^{\eta}+\bar{\nabla}_{k}^{\eta,d}\|^2 \right] \leq \tfrac{\yqro{16\sqrt{2\pi}(L_0^{\text{imp}})^2n}}{m} +2\mathbb{E}\left[\left\|\nabla {\bf f}^\eta (\bar x_{k})\right\|^2\right]   +\yqro{\tfrac{2(L_0^{\text{imp}}\sqrt n+1)^2}{\eta^2}}\mathbb{E}\left[\bar e_k\right].&&\\
&\mbox{(v)}\ \ \mathbb{E}\left[ \nabla {\bf f}^\eta(\bar x_k)^{\fyy{\top}}(\bar{g}_k^{\eta}+\bar{\nabla}_{k}^{\eta,d}) \right] \geq \tfrac{1}{2}\mathbb{E}\left[\|\nabla {\bf f}^\eta(\bar x_k)\|^2\right] -\yqro{\tfrac{(L_0^{\text{imp}}\sqrt n+1)^2}{2\eta^2}}\mathbb{E}\left[ \bar e_k\right].&&\\
&\mbox{(vi)}\ \ \mathbb{E}\left[ \|{\bar{w}^\eta}_{k}\|^2 \right] \leq \yqro{\tfrac{n^2}{\eta^2}}(L_{0,y}^f)^2 {\varepsilon_k}  .&&
\end{flalign*}
\end{lemma} 
\proof{Proof.}
\noindent (i) From Definition~\ref{def:s_terms_2smpec}, we can write
\begin{align*}
&\mathbb{E}\left[g^{\eta}_{i,k} \mid \mathcal{F}_k\right]\\  
&\yqro{= \mathbb{E}_{v_{i,k}}\left[\mathbb{E}_{\xi_{i, k}}\left[\tfrac{n}{\yqro{2}\eta^2} \left({\tilde{f}_i}(x_{i,k}+v_{i,k},y(x_{i,k}+v_{i,k},\xi_{i,k}),\xi_{i,k})\right)v_{i,k} \mid \mathcal{F}_k\cup\{v_{i,k}\}\right]\right]}\\
&\yqro{+ \mathbb{E}_{v_{i,k}}\left[\mathbb{E}_{\xi_{i, k}}\left[\tfrac{n}{\yqro{2}\eta^2} \left({\tilde{f}_i}(x_{i,k}-v_{i,k},y(x_{i,k}-v_{i,k},\xi_{i,k}),\xi_{i,k})\right)(-v_{i,k}) \mid \mathcal{F}_k\cup\{v_{i,k}\}\right]\right]}\\
& \ {=\mathbb{E}_{v_{i,k}}\left[\tfrac{n}{\eta^2}  f_i(x_{i,k}+v_{i,k},y(x_{i,k}+v_{i,k}))v_{i,k} \mid \mathcal{F}_k\right]}  \stackrel{\tiny \mbox{Lemma~\ref{SphericalSmooth} (i)}}{=} \nabla f_{i}^\eta(x_{i, k}).
\end{align*}

\noindent (ii) \yqro{From Assumption \ref{assum:fed2smpecs} (i), we know that the random implicit function ${\tilde{f}_i}(\bullet,y(\bullet,\xi_i),\xi_{i})$ is $L_0^{\text{imp}}(\xi_i)$-Lipschitz. Then by adding and subtracting $\mathbb{E}_{\hat v_{i,k}}[{\tilde{f}_i}(x_{i,k}+\hat v_{i,k},y(x_{i,k}+\hat v_{i,k},\xi_{i,k}),\xi_{i,k})]$, where $\hat v_{i,k}$ is uniformly distributed on $\eta \mathbb S$, we can utilize Lemma \ref{lem:Levy} by following the same approach in Lemma \ref{hetero:local:gradient} (ii) and obtain
\begin{align*}
\mathbb{E}_{v_{i,k}}\left[\left\|g^{\eta}_{i,k}\right\|^2 \mid \mathcal{F}_k \cup \{ \xi_{i, k}\} \right]  &\leq 16\sqrt{2\pi}(L_0^{\text{imp}}(\xi_i))^2n.
\end{align*}}
Taking expectations with respect to $\xi_{i,k}$ on both sides and recalling $L_0^{\text{imp}} \triangleq {\displaystyle \max_{i=1,\ldots,m}}\sqrt{\mathbb{E}[(L_0^{\text{imp}}(\xi_{i}))^2]}$, we get $\mathbb{E}\left[\|g^{\eta}_{i,k}\|^2 \mid \mathcal{F}_k\right] \leq \yqro{16\sqrt{2\pi}(L_0^{\text{imp}})^2n}$.  The result follows by taking expectations on both sides of the last relation.

\noindent (iii) Consider the definition of ${w_{i,k}^\eta}$.  Then we may bound $\mathbb{E}\left[ \|{w_{i,k}^\eta}\|^2\right]$ as follows.
\begin{align*}
&\mathbb{E}\left[ \|{w_{i,k}^\eta}\|^2\mid \mathcal{F}_k\cup \{ \xi_{i, k}\}\right] \\
&=  \tfrac{n^2}{\yqro{4}\eta^4} \mathbb{E}\left[ \|\left( \tilde{f}_i(x_{i,k}+{v_{i,k},y_{\varepsilon_k}(x_{i,k}+v_{i,k}, \xi_{i,k})},\xi_{i,k})- \tilde{f}_i(x_{i,k},{y_{\varepsilon_k}(x_{i,k}, \xi_{i,k})},\xi_{i,k}) \right.\right.\\
&\left.\left. -\tilde{f}_i(x_{i,k}+{v_{i,k},y(x_{i,k}+v_{i,k}, \xi_{i,k})},\xi_{i,k})+ \tilde{f}_i(x_{i,k},{y(x_{i,k}, \xi_{i,k})},\xi_{i,k})\right){v_{i,k}} \|^2\mid \mathcal{F}_k\cup \{ \xi_{i, k}\}\right]\\
&\leq \tfrac{n^2}{\yqro{2}\eta^4} \mathbb{E}\left[ \| (\tilde{f}_i(x_{i,k}+{v_{i,k},y_{\varepsilon_k}(x_{i,k}+v_{i,k}}, \xi_{i,k}),\xi_{i,k})-\tilde{f}_i(x_{i,k}+{v_{i,k},y(x_{i,k}+v_{i,k}, \xi_{i,k})},\xi_{i,k})){v_{i,k}} \|^2\right. \\ 
&\left. \mid \mathcal{F}_k\cup \{ \xi_{i, k}\} \right]   + \tfrac{n^2}{\yqro{2}\eta^4} \mathbb{E}\left[\| ( - \tilde{f}_i(x_{i,k},{y_{\varepsilon_k}(x_{i,k}}, \xi_{i,k}),\xi_{i,k})+ \tilde{f}_i(x_{i,k},{y(x_{i,k}, \xi_{i,k})},\xi_{i,k}))v_{i,k}\|^2\mid \mathcal{F}_k\cup \{ \xi_{i, k}\} \right]\\
&\leq \yqro{\tfrac{n^2}{\eta^2}} \left(L_{0,y}^{f}(\xi_{i,k})\right)^2 {\varepsilon_k}.
\end{align*}
{The last inequality is obtained by invoking Assumption~\ref{assum:fed2smpecs} and following similar steps in Lemma \ref{lemma:bounds_bilevel_6terms} (iv) and then, by invoking the definition of $ \varepsilon_k$.} Next, we take  expectations with respect to $\xi_{i,k}$ on both sides of the preceding relation.  The result follows by taking expectations on both sides of the resulting inequality.

\noindent \us{(iv,v,vi)} The proofs \us{are} similar to that of Lemma \ref{lemma:bounds_bilevel_6terms} (v,vi,vii) \us{respectively}  and \us{are} omitted.



Next, we derive an upper bound on the average consensus violation, after performing $k$ local steps by a client in Algorithm~\ref{alg:fed2:upper}.
\begin{lemma}[Bounds on average consensus violation]\label{lem:delay_cons_bound_fed2smpec}  \em                                                         
    Let Assumption~\ref{assum:fed2smpecs} hold. Consider Algorithm~\ref{alg:fed2:upper}.  The following \us{hold}.

\noindent (i)  [Recursive bound] For any communication round $r>0$, for all $T_r \leq k \leq T_{r+1}-1$, we have
 \begin{align*}
 \mathbb{E}\left[  \bar e_k  \right] 
& \leq   {3}\gamma^2(k-T_r) \displaystyle\sum_{t=T_r}^{k-1} \left( 8\mathbb{E}\left[\bar{e}_t\right] + 4B_2^2 \mathbb{E}\left[\left\|\nabla {\bf f}^{\eta}(\bar{x}_t)\right\|^2\right]+\tfrac{2B_1^2}{\eta^2}\right. \\
&\left.+ (4B_2^2+1)\yqro{16\sqrt{2\pi}(L_0^{\text{imp}})^2n} +\yqro{\tfrac{n^2}{\eta^2}}(L_{0,y}^f)^2 {\varepsilon_t}\right).
\end{align*}
\noindent (ii)  [Non-recursive bound] Let  {$0<\gamma \leq \left( \sqrt{24}H\right)^{-1}$.} Then, for any $r>0$, for all $T_r \leq k \leq T_{r+1}-1$, we have
\begin{align*}
\mathbb{E}\left[  \bar e_k  \right] 
& \leq \yqro{9}H\gamma^2\sum_{t=T_r}^{k-1}\left(\yqro{4}B_2^2\mathbb{E}[\|\nabla {\bf f}^{\eta}(\bar x_t)\|^2] +\tfrac{n^2}{\eta^2}(L_{0,y}^f)^2 \varepsilon_t \right)+ \yqro{18}H^2\gamma^2 \left(\yqro{\tfrac{B_1^2}{\eta^2}} + (4B_2^2+1)\yqro{8\sqrt{2\pi}(L_0^{\text{imp}})^2n} \right). 
\end{align*}
\end{lemma}
\proof{Proof.}  
\noindent {(i) 
In view of Algorithm~\ref{alg:fed2:upper}, for any $i$ at any communication round $r>0$, for all $T_r \leq k \leq T_{r+1}-1$, we have $x_{i,k+1} =  x_{i,k} - \gamma  ({g^{\eta,\varepsilon_k}_{i,k}}+\nabla^{\eta,d}_{i,k}).$ This implies that $x_{i,k} =  x_{i,k-1} - \gamma  ({g^{\eta,\varepsilon_k}_{i,k-1}}+\nabla^{\eta,d}_{i,k-1})$ for all $T_r+1 \leq k \leq T_{r+1}.$ Unrolling the preceding relation recursively, we obtain
\begin{align}\label{eqn:nonrecursive_local_fed2smpec}
x_{i,k} =  x_{i,T_r} - \gamma  \textstyle\sum_{t=T_r}^{k-1}({g^{\eta,\varepsilon_t}_{i,t}}+\nabla^{\eta,d}_{i,t}),\quad \hbox{for all } T_r+1 \leq k \leq T_{r+1}.
\end{align} 
From Algorithm~\ref{alg:fed2:upper}, we have $\hat{x}_r =x_{i,T_r}$. Invoking the definition of $\bar{x}_k$, we have $\bar{x}_{T_r} = \hat{x}_r$. This implies that $\bar{x}_{T_r} =x_{i,T_r}$ for all $i$ and $r$. Extending Lemma~\ref{lem:bar_x_recursion}, we can write
\begin{align}\label{eqn:nonrecursive_ave_fed2smpec}
\bar x_k =  x_{i,T_r} - \gamma  \textstyle\sum_{t=T_r}^{k-1}({\bar{g}^{\eta,\varepsilon_t}_{t}}+\bar{\nabla}^{\eta,d}_{t}),\quad \hbox{for all } T_r+1 \leq k \leq T_{r+1}.
\end{align}
Using \eqref{eqn:nonrecursive_local_fed2smpec} and \eqref{eqn:nonrecursive_ave_fed2smpec}, for all $T_r+1 \leq k \leq T_{r+1}$ we have
\begin{align*}
\mathbb{E}\left[  \bar e_k  \right]  
&=\tfrac{1}{m}\textstyle\sum_{i=1}^{m}\mathbb{E}\left[ \|x_{i, k}-\bar x_k \|^2   \right]   \\
&=\tfrac{1}{m}\textstyle\sum_{i=1}^{m}\mathbb{E}\left[ \left\|\gamma\textstyle\sum_{t=T_r}^{k-1} ({g_{i, t}^{\eta,\varepsilon_t}} +\nabla^{\eta,d}_{i,t})-\gamma\textstyle\sum_{t=T_r}^{k-1} ({\bar{g}_{t}^{\eta,\varepsilon_t}}+ \bar{\nabla}^{\eta,d}_{t}) \right\|^2  \right]   \\
&\leq \tfrac{\gamma^2(k-T_r)}{m}\textstyle\sum_{i=1}^{m}\textstyle\sum_{t=T_r}^{k-1} \mathbb{E}\left[ \left\| ({g_{i, t}^{\eta,\varepsilon_t}} +\nabla^{\eta,d}_{i,t})- ({\bar{g}_{t}^{\eta,\varepsilon_t}} + \bar{\nabla}^{\eta,d}_{t}) \right\|^2  \right],
\end{align*}
where the preceding relation follows from the inequality $\|  \textstyle\sum_{t=1}^{T}y_t\|^2\leq T\textstyle\sum_{t=1}^{T}\|y_t\|^2$ for any $y_t \in \mathbb{R}^n$ for $t \in [T]$. Consequently,
\begin{align*}
\mathbb{E}\left[  \bar e_k   \right]  
&\leq \tfrac{\gamma^2(k-T_r)}{m}\textstyle\sum_{t=T_r}^{k-1}\textstyle\sum_{i=1}^{m} \mathbb{E}\left[ \left\| ({g_{i, t}^{\eta,\varepsilon_t}} +\nabla^{\eta,d}_{i,t})- ({\bar{g}_{t}^{\eta,\varepsilon_t}} + \bar{\nabla}^{\eta,d}_{t}) \right\|^2   \right]\\
& = \tfrac{\gamma^2(k-T_r)}{m}\textstyle\sum_{t=T_r}^{k-1}\sum_{i=1}^{m}\left(\mathbb{E}\left[ \left\|  {g_{i, t}^{\eta,\varepsilon_t}} +\nabla^{\eta,d}_{i,t} \right\|^2   \right] +\mathbb{E}\left[ \left\| {\bar{g}_{t}^{\eta,\varepsilon_t}} + \bar{\nabla}^{\eta,d}_{t}  \right\|^2   \right] \right)\\
& -\tfrac{2\gamma^2(k-T_r)}{m}\textstyle\sum_{t=T_r}^{k-1}\sum_{i=1}^{m}\mathbb{E}\left[ (  {g_{i, t}^{\eta,\varepsilon_t}} +\nabla^{\eta,d}_{i,t} )^{\fyy{\top}}({\bar{g}_{t}^{\eta,\varepsilon_t}} + \bar{\nabla}^{\eta,d}_{t}  )   \right],
\end{align*}
Observing that  $\tfrac{1}{m}\textstyle\sum_{i=1}^{m}\mathbb{E}\left[ (  {g_{i, t}^{\eta,\varepsilon_t}} +\nabla^{\eta,d}_{i,t} )^{\fyy{\top}}({\bar{g}_{t}^{\eta,\varepsilon_t}} + \bar{\nabla}^{\eta,d}_{t}  ) \mid \mathcal{F}_{T_r}  \right]= {\left\|{\bar{g}_{t}^{\eta,\varepsilon_t}}+ \bar{\nabla}^{\eta,d}_{t}  \right\|^2},$ 
we obtain 
\begin{align*}
\mathbb{E}\left[  \bar e_k  \right]  
&\leq  \tfrac{\gamma^2(k-T_r)}{m}\textstyle\sum_{t=T_r}^{k-1}\sum_{i=1}^{m} \mathbb{E}\left[ \left\|  {g_{i, t}^{\eta,\varepsilon_t}} +\nabla^{\eta,d}_{i,t} \right\|^2  \right]\\
& \leq \tfrac{{3}\gamma^2(k-T_r)}{m}\textstyle\sum_{t=T_r}^{k-1}\sum_{i=1}^{m} \left(\mathbb{E}\left[  \|  {g_{i, t}^{\eta}}\|^2 +\|{w_{i,t}^\eta}\|^2+\|\nabla^{\eta,d}_{i,t}  \|^2  \right]\right)\\
& \stackrel{\tiny {\mbox{Lemma\mbox{ }\ref{lemma:bounds_2smpec_6terms}}}}{\leq} \tfrac{{3}\gamma^2(k-T_r)}{m}\textstyle\sum_{t=T_r}^{k-1}\sum_{i=1}^{m} \left(\yqro{16\sqrt{2\pi}(L_0^{\text{imp}})^2n}+\|\nabla^{\eta,d}_{i,t}  \|^2+\yqro{\tfrac{n^2}{\eta^2}}(L_{0,y}^f)^2{\varepsilon_t} \right).
\end{align*}
Invoking Lemma~\ref{lem:bound_on_dist2_ave},  we obtain the inequality in (i).
 
\noindent (ii) This result is obtained by applying Lemma~\ref{lem:aggr_a_k_extension} on the inequality in part (i).}

\subsection{Proof of Theorem~\ref{thm:fed2smpec}}
\,

\proof{Proof.} 
\noindent (i) The proof can be done in a similar vein to that of~\cite[Theorem 4 (a)]{cui2022complexity} and is omitted.

\noindent (ii)
From the $L$-smoothness property of the implicit function, where $L=\yqro{\tfrac{L_0^{\text{imp}}\sqrt n+1}{\eta}}$, we have
\begin{align*}
 {\bf f}^{\eta}(\bar x_{k+1})    \leq   {\bf f}^{\eta}(\bar x_{k}) +  \nabla {\bf f}^{\eta}(\bar x_k)^{\fyy{\top}} (\bar x_{k+1}-\bar x_k)   +\tfrac{({L_0^{\tiny \mbox{imp}}}n+1)}{2\eta} \|\bar x_{k+1}-\bar x_k\|^2.
\end{align*}
In view of the recursion $\bar{x}_{k+1}=\bar{x}_k-\gamma({\bar{g}_k^{\eta,\varepsilon_k}}+\bar{\nabla}^{\eta,d}_{k})$, we have
\begin{align*}
{\bf f}^{\eta}(\bar x_{k+1})    \leq   {\bf f}^{\eta}(\bar x_{k}) -\gamma   \nabla {\bf f}^{\eta}(\bar x_k)^{\fyy{\top}} ({\bar{g}_k^{\eta,\varepsilon_k}} +\bar\nabla^{\eta,d}_{k})  +\yqro{\tfrac{L_0^{\text{imp}}\sqrt n+1}{2\eta}}\gamma^2 \|{\bar{g}_k^{\eta,\varepsilon_r}}+\bar{\nabla}^{\eta,d}_{k}\|^2.
\end{align*}
In view of {$\bar{g}^{\eta,\varepsilon_k}_{k} =\bar{g}^{\eta}_{k} +{\bar{w}^\eta}_{k}$}, we have
\begin{align*}
 {\bf f}^{\eta}(\bar x_{k+1})   & \leq   {\bf f}^{\eta}(\bar x_{k}) -\gamma   \nabla {\bf f}^{\eta}(\bar x_k)^{\fyy{\top}}(\bar{g}_k^{\eta}+\bar{\nabla}_{k}^{\eta,d})  -\gamma   \nabla {\bf f}^{\eta}(\bar x_k)^{\fyy{\top}}{\bar{w}^\eta}_{k} \\ 
 &  +\yqro{\tfrac{L_0^{\text{imp}}\sqrt n+1}{2\eta}}\gamma^2 \|{\bar{g}^{\eta}_{k}+\bar{\nabla}_{k}^{\eta,d} +{\bar{w}^\eta}_{k}}\|^2\\
 &   \leq   {\bf f}^{\eta}(\bar x_{k}) -\gamma  \nabla {\bf f}^{\eta}(\bar x_k)^{\fyy{\top}}(\bar{g}_k^{\eta}+\bar{\nabla}_{k}^{\eta,d})+\tfrac{\gamma}{{16}}\|\nabla {\bf f}^{\eta}(\bar x_k)\|^2+4\gamma\|{\bar{w}^\eta}_{k}\|^2  \\
 & +\yqro{\tfrac{L_0^{\text{imp}}\sqrt n+1}{2\eta}}\gamma^2\left( \|\bar{g}^{\eta}_{k}+\bar{\nabla}_{k}^{\eta,d}\|^2 +\|{\bar{w}^\eta}_{k}\|^2\right),
\end{align*}
where in the second inequality, we utilized the following relation: $|u_1^{\fyy{\top}}u_2| \leq 0.5\left(\tfrac{1}{8}\|u_1\|^2 +8\|u_2\|^2\right)$ for any $u_1,u_2 \in \mathbb{R}^n$. Taking expectations on both sides and requiring that $\gamma \leq \yqro{\tfrac{\eta}{16(L_0^{\text{imp}}\sqrt n+1)}}$, we obtain 
\begin{align}\label{eqn:thm_fed2smpec_term0}
 \mathbb{E}\left[{\bf f}^{\eta}(\bar x_{k+1}) \right]      & \leq   \mathbb{E}\left[{\bf f}^{\eta}(\bar x_{k}) \right]-\gamma  \mathbb{E}\left[ \nabla {\bf f}^{\eta}(\bar x_k)^{\fyy{\top}}   (\bar{g}_k^{\eta}+\bar{\nabla}_{k}^{\eta,d})\right]+\tfrac{\gamma}{{16}}\mathbb{E}\left[\|\nabla {\bf f}^{\eta}(\bar x_k)\|^2\right]\notag\\
&    +\yqro{\tfrac{L_0^{\text{imp}}\sqrt n+1}{\eta}}\gamma^2   \mathbb{E}\left[\|\bar{g}^{\eta}_{k}+\bar{\nabla}_{k}^{\eta,d}\|^2\right]+5\gamma \mathbb{E}\left[\|{\bar{w}^\eta}_{k}\|^2\right] .
\end{align}
From {Lemma~\ref{lemma:bounds_2smpec_6terms} (iv) and (v)}, we have
\begin{align}\label{eqn:thm_fed2smpec_term1}
-\gamma\mathbb{E}\left[ \nabla {\bf f}^\eta(\bar x_k)^{\fyy{\top}}(\bar{g}_k^{\eta}+\bar{\nabla}_{k}^{\eta,d}) \right] \leq  -\tfrac{\gamma}{2}\mathbb{E}\left[\|\nabla {\bf f}^\eta(\bar x_k)\|^2\right] +\gamma\yqro{\tfrac{(L_0^{\text{imp}}\sqrt n+1)^2}{2\eta^2}}\mathbb{E}\left[ \bar e_k\right],
\end{align}
and $
\mathbb{E}\left[\|\bar g_k^{\eta}+\bar{\nabla}_{k}^{\eta,d}\|^2 \right] \leq \tfrac{\yqro{16\sqrt{2\pi}(L_0^{\text{imp}})^2n}}{m} +2\mathbb{E}\left[\left\|\nabla {\bf f}^\eta (\bar x_{k})\right\|^2\right]   +\yqro{\tfrac{2(L_0^{\text{imp}}\sqrt n+1)^2}{\eta^2}} \mathbb{E}\left[\bar e_k\right],$
respectively. In view of $\gamma \leq \yqro{\tfrac{\eta}{16(L_0^{\text{imp}}\sqrt n+1)}}$, multiplying both sides of the preceding inequality by $\yqro{\tfrac{L_0^{\text{imp}}\sqrt n+1}{\eta}}\gamma^2$, we obtain 
\begin{align}\label{eqn:thm_fed2smpec_term2}
\yqro{\tfrac{L_0^{\text{imp}}\sqrt n+1}{\eta}}\gamma^2\mathbb{E}\left[\|\bar g_k^{\eta}+\bar{\nabla}_{k}^{\eta,d}\|^2   \right]& \leq \yqro{\tfrac{16\sqrt{2\pi}(L_0^{\text{imp}}\sqrt n+1)(L_0^{\text{imp}})^2n}{m\eta}}\gamma^2+\tfrac{\gamma}{8}\mathbb{E}\left[\|\nabla {\bf f}^\eta (\bar x_{k})\|^2\right] +\gamma\yqro{\tfrac{(L_0^{\text{imp}}\sqrt n+1)^2}{8\eta^2}}\mathbb{E}\left[\bar e_k\right].
\end{align}
Also, from {Lemma~\ref{lemma:bounds_2smpec_6terms} (vi)}, we have $
5\gamma\mathbb{E}\left[\|{\bar{w}^\eta}_{k}\|^2\right] \leq \tfrac{\yqro{5} \gamma n^2}{\eta^2}(L_{0,y}^f)^2\varepsilon_k .$ From the preceding relation and {\eqref{eqn:thm_fed2smpec_term0}--\eqref{eqn:thm_fed2smpec_term2}}, we obtain 
\begin{align*} 
 \mathbb{E}\left[{\bf f}^{\eta}(\bar x_{k+1}) \right]      & \leq   \mathbb{E}\left[{\bf f}^{\eta}(\bar x_{k}) \right] -\tfrac{\gamma}{4}\mathbb{E}\left[\|\nabla {\bf f}^{\eta}(\bar x_k)\|^2\right] +\yqro{\tfrac{16\sqrt{2\pi}(L_0^{\text{imp}}\sqrt n+1)(L_0^{\text{imp}})^2n}{m\eta}}\gamma^2 \\
 & +\gamma\tfrac{5\yqro{(L_0^{\text{imp}}\sqrt n+1)^2}}{8\eta^2}\mathbb{E}\left[\bar e_k\right]+ {\tfrac{\yqro{5} \gamma n^2}{\eta^2}(L_{0,y}^f)^2\varepsilon_k}.
\end{align*}
This implies that  $\mathbb{E}\left[{\bf f}^{\eta}(\bar x_{k+1}) \right]       \leq   \mathbb{E}\left[{\bf f}^{\eta}(\bar x_{k}) \right] -\tfrac{\gamma}{4}\mathbb{E}\left[\|\nabla {\bf f}^{\eta}(\bar x_k)\|^2\right] +\tfrac{\gamma^2 \Theta_1}{m}    +{\gamma\Theta_2\mathbb{E}\left[\bar e_k\right]+ \gamma\Theta_3\varepsilon_k},$ where we define $\Theta_1 :=\yqro{\tfrac{16\sqrt{2\pi}(L_0^{\text{imp}}\sqrt n+1)(L_0^{\text{imp}})^2n}{m\eta}}$, $\Theta_2:=\tfrac{5\yqro{(L_0^{\text{imp}}\sqrt n+1)^2}}{8\eta^2}$,  and ${\Theta_3:= \tfrac{\yqro{5} n^2}{\eta^2}(L_{0,y}^f)^2}$.
From {Lemma~\ref{lem:delay_cons_bound_fed2smpec} (ii)},  for {$0<\gamma \leq \left( \sqrt{24}H\right)^{-1}$}, we have
\begin{align*}
\mathbb{E}\left[  \bar e_k  \right] 
& \leq 36B_2^2H\gamma^2\sum_{t=T_r}^{k-1}\mathbb{E}[\|\nabla {\bf f}^{\eta}(\bar x_t)\|^2] + H\gamma^2\Theta_4\sum_{t=T_r}^{k-1} \varepsilon_t + H^2\gamma^2\Theta_5, 
\end{align*}
where ${\Theta_4:=\tfrac{\yqro{9}n^2}{\eta^2}(L_{0,y}^f)^2}$ and $\Theta_5:={\yqro{18}\left(\yqro{\tfrac{B_1^2}{\eta^2}} + (4B_2^2+1)\yqro{8\sqrt{2\pi}(L_0^{\text{imp}})^2n} \right)}$.
Let $\gamma \leq \tfrac{1}{12\sqrt{6}HB_2\sqrt{\Theta_2}}$,  we have ${{36B_2^2}H\gamma^2\Theta_2} \leq  {\tfrac{1}{24H}}$.
Combining the two preceding inequalities, we obtain
\begin{align*} 
 \mathbb{E}\left[{\bf f}^{\eta}(\bar x_{k+1}) \right]      & \leq   \mathbb{E}\left[{\bf f}^{\eta}(\bar x_{k}) \right] -\tfrac{\gamma}{4}\mathbb{E}\left[\|\nabla {\bf f}^{\eta}(\bar x_k)\|^2\right] +{\tfrac{\gamma }{24H}}\textstyle\sum_{t=T_r}^{k-1}\mathbb{E}\left[\left\|\nabla {\bf f}^{\eta}(\bar{x}_t)\right\|^2\right]+\tfrac{\gamma^2 \Theta_1}{m}  \\
 & +{H\gamma^3\Theta_2\Theta_4\sum_{t=T_r}^{k-1}\varepsilon_t+ \gamma\Theta_3\varepsilon_k  +H^2\gamma^3\Theta_2\Theta_5}.
\end{align*} 
Summing the preceding relation from $k:=0,\ldots,K:=T_{R}-1$ for some $R\geq 1$, we have
\begin{align*} 
\mathbb{E}\left[{\bf f}^{\eta}(\bar x_{T_R}) \right]      
& \leq   \mathbb{E}\left[{\bf f}^{\eta}(\bar x_{0}) \right] -\textstyle\sum_{k=0}^{T_R-1}\tfrac{\gamma}{4}\mathbb{E}\left[\|\nabla {\bf f}^{\eta}(\bar x_k)\|^2\right] +{\tfrac{\gamma}{24}}\sum_{k=0}^{T_R-1}\mathbb{E}\left[\left\|\nabla {\bf f}^{\eta}(\bar{x}_k)\right\|^2\right]+\tfrac{T_R\gamma^2 \Theta_1}{m} \\
&  +H^2\gamma^3\Theta_2\Theta_4 \sum_{k=0}^{T_R-1} \varepsilon_k + \gamma\Theta_3\sum_{k=0}^{T_R-1}\varepsilon_k +T_RH^2\gamma^3\Theta_2\Theta_5 ,
\end{align*}
where we used $\sum_{k=0}^{T_R-1} \sum_{t=T_r}^{k-1}\varepsilon_t \leq \sum_{k=0}^{T_R-1} \sum_{t=T_r}^{T_{r+1}-1}\varepsilon_t \leq H\sum_{k=0}^{T_R-1}  \varepsilon_k$.  From $\tfrac{\gamma}{8} \leq \tfrac{\gamma}{4} -\tfrac{\gamma}{24} $, we obtain
  \begin{align*} 
  \textstyle\sum_{k=0}^{T_R-1}\tfrac{\gamma}{8}\mathbb{E}\left[\|\nabla {\bf f}^{\eta}(\bar x_k)\|^2\right]   & \leq   \mathbb{E}\left[{\bf f}^{\eta}(\bar x_{0}) \right] - \mathbb{E}\left[{\bf f}^{\eta}(\bar x_{T_R}) \right]   +\tfrac{T_R\gamma^2 \Theta_1}{m}  \\
 & {+H^2\gamma^3\Theta_2\Theta_4 \sum_{k=0}^{T_R-1}\ \varepsilon_k + \gamma\Theta_3\sum_{k=0}^{T_R-1}\varepsilon_k +T_RH^2\gamma^3\Theta_2\Theta_5}.
\end{align*}
Multiplying both sides by $\tfrac{8}{\gamma T_R}$ and using the definition of $k^*$, we obtain  
\begin{align*} 
\mathbb{E}\left[\|\nabla {\bf f}^{\eta}(\bar x_{k^*})\|^2\right]   & \leq   8(\gamma T_R)^{-1}(\mathbb{E}\left[{\bf f}^{\eta}(x_{0}) \right] -  {\bf f}^{\eta,*}  )   +\tfrac{8\gamma \Theta_1}{m}  {+8H^2\gamma^2\Theta_2\Theta_5}\\
 & {+8(T_R)^{-1}H^2\gamma^2\Theta_2\Theta_4 \sum_{k=0}^{T_R-1}\ \varepsilon_k + 8(T_R)^{-1}\Theta_3\sum_{k=0}^{T_R-1}\varepsilon_k }.
\end{align*}
 
\noindent (iii) In view of $t_k=\tau \text{ln}(k+1)$ where $\tau\geq\tfrac{-1}{\text{ln}(1-1/\kappa^2_F)}$, then we have 
$$(k+1)\varepsilon_k\leq (k+1)B(1-1/\kappa^2_F)^{\tau \text{ln}(k+1)}=B((1-1/\kappa^2_F)^{\tau}e)^{\text{ln}(k+1)}\leq B.$$
Therefore, $\varepsilon_k\leq\frac{B}{k+1}$. 
From~\cite[Lemma 10]{yousefian2017smoothing}, we have 
$$\sum_{k=0}^{T_R-1}\varepsilon_k \leq B\sum_{k=0}^{T_R-1}\tfrac{1}{k+1}\leq B(1+\text{ln}(K)) =B\text{ln}( K\exp(1)).$$ 
Finally, by substituting $\gamma:=\yqro{\sqrt{\tfrac{\eta m}{Kn^{1.5}(L_0^{\text{imp}})^3}}}$ and $H:=   \left\lceil\sqrt[4]{\tfrac{K}{m^3}}\right\rceil$ in the error bound in (ii), we obtain
\begin{align*} 
\mathbb{E}\left[\|\nabla {\bf f}^{\eta}(\bar x_{k^*})\|^2\right]   
& \leq    \frac{8(\mathbb{E}\left[{\bf f}^{\eta}(x_{0}) \right] -  {\bf f}^{\eta,*})\yqro{\tfrac{n^{0.75}(L_0^{\text{imp}})^{1.5}}{\eta^{0.5}}}+ 8\Theta_1\yqro{\tfrac{\eta^{0.5}}{n^{0.75}(L_0^{\text{imp}})^{1.5}}}  + 8\Theta_2\Theta_5\yqro{\tfrac{\eta}{n^{1.5}(L_0^{\text{imp}})^{3}}}}{\sqrt{mK}}  \\
&+ \frac{8B\Theta_2\Theta_4\yqro{\tfrac{\eta}{n^{1.5}(L_0^{\text{imp}})^{3}}}\text{ln}(K)}{\sqrt{mK^3}}  + \frac{8B\Theta_3\text{ln}(K\exp(1))}{K} .
\end{align*}
Therefore, the iteration complexity for $\mathbb{E}\left[\|\nabla {\bf f}^{\eta}(\bar x_{k^*})\|^2\right]\leq \epsilon$ is 
 $$K_\epsilon:=\mathcal{\tilde{O}}\left( \tfrac{\yqro{\tfrac{n^{1.5}(L_0^{\text{imp}})^{3}}{\eta}+\tfrac{n(L_0^{\text{imp}})^{-2}}{\eta^{2}}}}{m\epsilon^2} + \tfrac{\yqro{\tfrac{n(L_0^{\text{imp}})^{-2}}{\eta^{-2/3}}}}{m^{1/3}\epsilon^{2/3}} + \tfrac{\yqro{n^2}}{\epsilon}\right).$$  

\noindent  (iv)  From $H:=  \left\lceil\sqrt[4]{\tfrac{K}{m^3}}\right\rceil$,  we obtain $R=\mathcal{O}(\tfrac{K}{H}) = \mathcal{O}\left( (mK)^{\tfrac{3}{4}}\right) $.  Note that there \us{is} no \us{communication} among the clients during the implementations of the lower-level calls, implying the overall communication complexity of $\mathcal{O}\left( (mK)^{\tfrac{3}{4}}\right) $.   The total number of lower-level projection steps per client is $ \textstyle\sum_{k=1}^{K}  t_k = \tau\textstyle\sum_{k=1}^{K} \text{ln}(k+1)  = \mathcal{O}( \tau K\text{ln}(K)).$ Multiplying this by the number of clients yields $\tilde{\mathcal{O}}(\tfrac{mK}{-\text{ln}(1-1/\kappa^2_F)} ).$ 

 

\section{Appendix}\label{sec:appendix}
\yqro{This section presents \fyy{two supplementary} materials that are omitted in the main text.

First, we \fyy{prove the claim in \fyy{S}ection \ref{sec:hyper} stating} that the \fyy{upper-level} stochastic objective function \fyy{of the hyperparameter optimization problem} is Lipschitz with respect to the inner variable.
\begin{lemma}\em
Consider the hyperparameter optimization problem in \fyy{S}ection \ref{sec:hyper}. Let $\xi_i\in \mathcal D_i$ denote a random data sample of client $i$. We have ${\tilde{f}_i}(x,y,\xi_i)=\log\left( 1+\exp(-v_{i,\ell}U_{i,\ell}^{\fyy{\top}}y)\right)$ is $L_{0,y}^f(\xi_i)$-Lipschitz with respect to $y$, where $L_{0,y}^f(\xi_i)=\|U_{i,\ell}\|$, and $\displaystyle{\max_{i=1,\ldots,m}}\sqrt{\mathbb{E}_{\xi_i\in \mathcal D_i}[(L_{0,y}^f(\xi_i))^2]}<\infty.$
\end{lemma}
\proof{Proof.} 
From \cite[Theorem 3.61]{beck2017first}, we know that given a differentiable proper convex function $h(x)$, and $x\in X\subseteq \text{int dom}h$, we have $h(x)$ is $L$-Lipschitz continuous if $\|\nabla h(x)\|\leq L$. Therefore, it suffices to show that $\|\nabla_y {\tilde{f}_i}(x,y,\xi_i)\| \leq \|U_{i,\ell}\|$.

We have $\nabla_y {\tilde{f}_i}(x,y,\xi_i)= \tfrac{-v_{i,\ell}U_{i,\ell}^{\fyy{\top}} \, e^{-v_{i,\ell}U_{i,\ell}^{\fyy{\top}}y}}{1+e^{-v_{i,\ell}U_{i,\ell}^{\fyy{\top}}y}}$. 
Therefore, $\|\nabla_y {\tilde{f}_i}(x,y,\xi_i)\|= \|v_{i,\ell}U_{i,\ell}\|\tfrac{e^{-v_{i,\ell}U_{i,\ell}^{\fyy{\top}}y}}{1+e^{-v_{i,\ell}U_{i,\ell}^{\fyy{\top}}y}}\leq \|v_{i,\ell}U_{i,\ell}\|$. Since $v_{i,\ell}\in \{-1,1\}$, we have $\|\nabla_y {\tilde{f}_i}(x,y,\xi_i)\| \leq \|U_{i,\ell}\|$.  Furthermore, given that $\|U_{i,\ell}\|$ is finite (e.g., in MNIST, each of the coordinates of the input image ranges between $0$ to $1$ after data normalization), we have $\displaystyle{\max_{i=1,\ldots,m}}\sqrt{\mathbb{E}_{\xi_i\in \mathcal D_i}[\|U_{i,\ell}\|^2]}<\infty.$
\bigskip

Next, we \fyy{provide a simple example to} demonstrate that the implicit function \fyy{in bilevel optimization may be} nondifferentiable and nonconvex. Consider a constrained bilevel optimization problem given as
\begin{align}\label{eg:bl}
\begin{aligned}
&\min_{x\in X \subseteq \mathbb R^n} \, \tfrac{1}{2}\|x+{\bf 1}_n -y(x)\|^2 \\
&\,\, \text{s.t. }y(x)=\text{arg min}_{y\in \mathbb R_{+}^n} \|y-x\|^2.
\end{aligned}
\end{align}
Assume \fyy{that} set $X$ is nonempty, closed and convex. Let $f(\bullet)=\tfrac{1}{2}\|\bullet+{\bf 1}_n -y(\bullet)\|^2$ denote the implicit function. The lower-level problem can be expressed as $y(x)=\mathcal P_{\mathbb R_{+}^n}(x)$, where $\mathcal P_{\mathbb R_{+}^n}(\bullet)$ denotes the projection operator onto the n-dimensional nonnegative orthant. Note that $f(x)=\sum_{i=1}^n g_i(x_i)$, where $g_i(x_i)=\tfrac{1}{2}(x_i+1-\max\{0,x_i\})^2$, and $x_i\in \mathbb R$ is the $i$-th element of $x$. $g_i(x_i)$ can be written as
$$g_i(x_i)=
\begin{cases}
\tfrac{1}{2}, & \text{if } x_i \ge 0, \\
\tfrac{1}{2}(x_i+1)^2, & \text{if } x_i < 0.
\end{cases}$$
Let \fyy{$X=\prod_{i=1}^n[-a_i,a_i]$, where $a_i>0$ for all $i=1,\ldots,n$. \yqro{We first show that $g_i(x_i)$ is Lipschitz continuous in $x_i$. (Case 1) When $x_i\geq 0$, we obtain that $g_i(x_i)$ is $L_i$-Lipschitz, where $L_i=0$. (Case 2) When $x_i < 0$, we have $|\nabla g_i(x_i)|=|x_i+1|\leq \max\{1,|1-a_i|\}$. Then, we have $L_i=\max\{1,|1-a_i|\}$. (Case 3) Take two arbitrary points $y_i \in [-a_i,0)$, $z_i \in (0,a_i]$. We have $|g_i(y_i)-g_i(z_i)| = |0.5(y_i+1)^2-0.5| = 0.5|y_i(y_i+2)| = 0.5 |y_i||y_i+2|$. We also have $|y_i-z_i|>|y_i|$. Then, we obtain $\tfrac{|g_i(y_i)-g_i(z_i)|}{|y_i-z_i|} < \tfrac{0.5 |y_i||y_i+2|}{|y_i|}=0.5|y_i+2|\leq 0.5 \max\{2,|2-a_i|\}$. Therefore, by combining all cases, we can conclude that $g_i(x_i)$ is $L_i$-Lipschitz in $x_i$ over $[-a_i,a_i]$, where $L_i=\max\{1,|1-a_i|,|1-0.5a_i|\}$.} Next, we prove that $f(x)$ is Lipschitz continuous, nondifferentiable, and nonconvex in terms of $x$. Let us denote  $ L = \max_{i=1,\dots,n}\{L_i\}$}. We have
\begin{align*}
|f(x)-f(z)|^2 &= |\textstyle \sum_{i=1}^n g_i(x_i)-\textstyle \sum_{i=1}^n g_i(z_i)|^2\\
&=|(g_1(x_1)-g_1(z_1))+\dots+(g_n(x_n)-g_n(z_n))|^2\\
&\stackrel{\text{Cauchy-Schwarz inequality}}{\leq}n|(g_1(x_1)-g_1(z_1))|^2+\dots+n|(g_n(x_n)-g_n(z_n))|^2\\
&\leq nL_1^2|x_1-z_1|^2+\dots+nL_n^2|x_n-z_n|^2 \\
& = nL^2 \|x-z\|^2.
\end{align*}
Thus, the implicit function $f(x)$ is {\it $L\sqrt{n}$--Lipschitz continuous}. Next, we show that $f(x)$ is nonconvex. Recall that a function $f: X\to \mathbb R$ is convex over a closed convex set $X\subseteq \mathbb R^n$ if and only if for any $\lambda\in [0,1]$ and any $x,z \in X$, the following holds
$$f(\lambda x+ (1-\lambda)z)\leq \lambda f(x)+ (1-\lambda)f(z).$$
Let $\lambda=0.5$, $x={\bf 1}_n$, and $z=-{\bf 1}_n$, where ${\bf 1}_n$ denotes an n-dimensional vector with all elements equal to $1$. Then, we obtain $f(\lambda x+ (1-\lambda)z)= \tfrac{n}{2} > \tfrac{n}{4}=  \lambda f(x)+ (1-\lambda)f(z) $. Therefore, $f(x)$ is a {\it nonconvex} function over $X$. Furthermore, by noting that functions $g_i(x_i)$ for all $i=1,\dots,n$ are independent of each other in terms of the variable $x_i$, and $g_i(x_i)$ are nondifferentiable at $x_i=0$, we can conclude that $f(x)=\sum_{i=1}^n g_i(x_i)$ is {\it nondifferentiable} at $x$ when at least one element $x_i=0$. Figure \ref{fig:NDNCBiLv} shows the implicit function $f(\bullet)$ in \eqref{eg:bl} when $n=2$, visualizing its nondifferentiability and nonconvexity. 
\begin{figure}[htbp]
  \centering
  \includegraphics[width=0.6\textwidth]{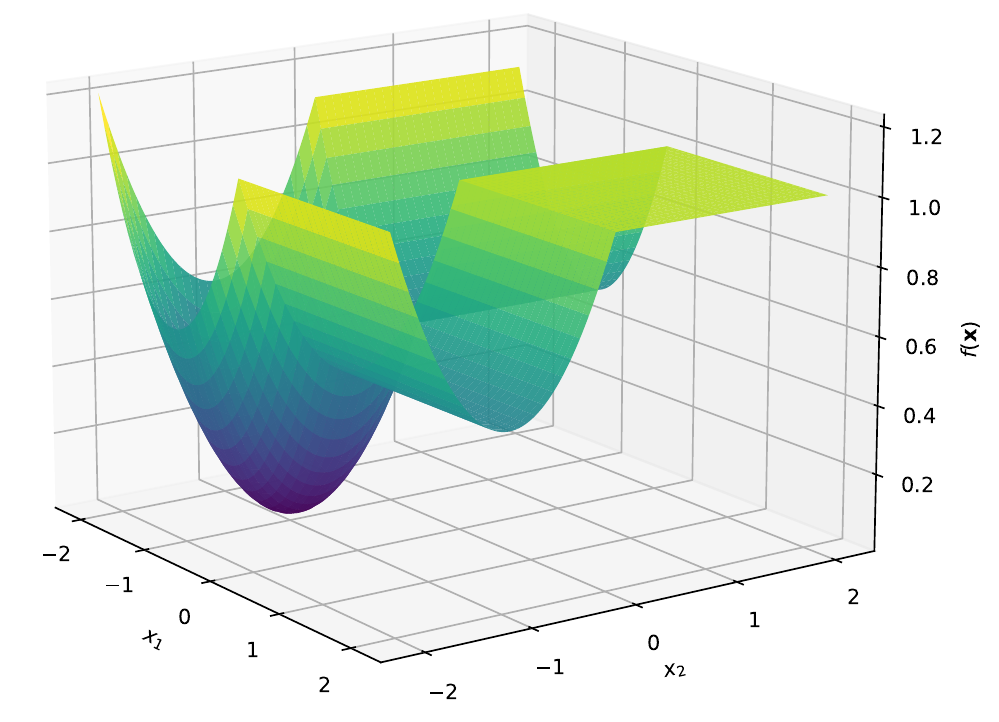}
  \caption{The implicit function $f(x)$ in \eqref{eg:bl} when $n=2$.}
  \label{fig:NDNCBiLv}
\end{figure}


}

\end{document}